\documentclass{article}
 \usepackage{amssymb,amsmath,amsfonts}
\newtheorem{lemma}{Lemma}[section]
 \newtheorem{theorem}{Theorem}[section]
 \newtheorem{proposition}{Proposition}[section]
 \newtheorem{corollary}{Corollary}[section]
 \newtheorem{remark}{Remark}[section]

\numberwithin{equation}{section}

\newcommand{\p}{\alpha}
\newcommand{\pe}{\psi_E}
\newcommand{\po}{\psi_0}
\newcommand{\R}{\mathbb{R}}
\newcommand{\C}{\mathbb{C}}
\renewcommand{\le}{\leqslant}
\renewcommand{\ge}{\geqslant}
\newcommand{\<}{\langle}
\renewcommand{\>}{\rangle}
\newcommand{\Dg}{Dg|_{\psi_E}}
\newcommand{\xs}{\langle x \rangle^{\sigma}}
\newcommand{\xsn}{\langle x \rangle^{-\sigma}}
\newcommand{\eht}{e^{-iH(t-\tau)}}
\newcommand{\ehs}{e^{-iH(t-s)}}
\newcommand{\ehts}{e^{-iH(\tau-s)}}
\newcommand{\xst}{\langle x \rangle^{2\sigma}}
\newcommand{\Ls}{L_{\sigma}^2}
\newcommand{\Lsn}{L_{-\sigma}^2}
\newcommand{\om}{\Omega}
\newcommand{\ipt}{(\frac{1}{2}-\frac{1}{p_2})}
\newcommand{\ipo}{(\frac{1}{2}-\frac{1}{p_1})}
\newcommand{\ip}{(\frac{1}{2}-\frac{1}{p})}

\begin{document}
\title{On the stability of ground states in 4D and 5D nonlinear Schr\" odinger equation including subcritical cases}
\author{E. Kirr and \" O. M\i zrak \thanks{Department of Mathematics, University of Illinois at Urbana-Champaign}}
\maketitle

\begin{abstract}
We consider a class of nonlinear Schr\"{o}dinger equation in four
and five space dimensions with an attractive potential. The
nonlinearity is local but rather general encompassing for the first
time both subcritical and supercritical (in $L^2$) nonlinearities.
We show that the center manifold formed by localized in space
periodic in time solutions (bound states) is an attractor for all
solutions with a small initial data. The proof hinges on dispersive
estimates that we obtain for the time dependent, Hamiltonian,
linearized dynamics around a one parameter family of bound states
that ``shadows" the nonlinear evolution of the system. The methods
we employ are an extension to higher dimensions, hence different
linear dispersive behavior, and to rougher nonlinearities of our
previous results \cite{kz:as2d,kz:as2d2,km:asd3}.
 \end{abstract}

\section{Introduction}
In this paper we study the long time behavior of solutions of the
nonlinear Schr\" odinger equation (NLS) with potential in four and five
space dimensions (4-d and 5-d):
\begin{align}
i\partial_t u(t,x)&=[-\Delta_x+V(x)]u+g(u), \quad t\in\R ,\quad x\in\R^N,\quad N=4,5 \label{u} \\
u(0,x)&=u_0(x) \label{ic}
\end{align}
where the local nonlinearity is constructed from the real valued,
odd, $C^1$ function $g:\mathbb{R}\mapsto\mathbb{R}$ which is twice differentiable except maybe at zero and satisfies:
\begin{equation}\label{gest}
g(0)=0,\ g'(0)=0, \quad |g''(s)|\leq
C(|s|^{\alpha_1-1}+|s|^{\p_2-1}),\quad \textrm{for
$s\not=0$},\end{equation} with
\begin{equation}
\alpha_0(N)=\frac{2-N+\sqrt{N^2+12N+4}}{2N}<\alpha_1\leq\p_2<\frac{4}{N-2}\quad\textrm{for
N=4, 5.}\label{alpha12}\end{equation} $g$ is then extended to a
complex function via the gauge symmetry:
\begin{equation}\label{gsym}
g(e^{i\theta}z)=e^{i\theta}g(z),\quad g(\bar z)=\overline{g(z)}.
\end{equation}
Note that $g$ is not necessarily twice differentiable at 0, e.g.
$g(z)=|z|^{\frac{5}{6}}z$.

We are going to show that the manifold of periodic solutions of
\eqref{u} (center manifold) is a global attractor for all small
initial data. More precisely, for $u_0\in H^1\cap
L^{\frac{\alpha_2+2}{\alpha_2+1}}$ with sufficiently small norm the
solution of \eqref{u}-\eqref{ic} can be decomposed as follows:
 $$u(t,x)=e^{i\theta(t)}\psi_{E(t)}(x)+r(t,x)$$
where for each fixed time $t_1\in\R ,\ E=E(t_1),$ we have that $u_E(t,x)=e^{-iEt}\psi_E(x)$
is a periodic solution of \eqref{u} and, as $|t|\rightarrow\infty,$  $r(t)\in H^1(\R^N)$
converges strongly to zero in $L^p(\R^N),\ 2<p<2N/(N-2),$ spaces and
weakly in $H^1(\R^N),$ see Section \ref{se:main}
for more details. We can also show that the full dynamics converges to a
certain periodic solution, i.e. $\psi_{E(t)}\rightarrow
\psi_{E_{\pm\infty}},$ for $t\rightarrow\pm\infty,$ provided we restrict the
range of nonlinearities to the supercritical regime $\alpha_1 >4/N,$
see Corollary \ref{cor:as}. In this case we generalize the results in
\cite{sw:mc2,pw:cm}.

Moreover, for the remaining range $\alpha_0(N)<\alpha_1\le 4/N,$ we
could show a type of asymptotic stability for periodic solutions
(bound states) of \eqref{u}. In other words, if $e^{-iEt}\psi_E(x),\
\psi_E\not\equiv 0$ is a (small) periodic solution of \eqref{u} then
there exists $\varepsilon>0$ depending on $\psi_E$ such that for all
initial data $u_0\in H^1\cap L^{\frac{\alpha_2+2}{\alpha_2+1}}$
satisfying
$$\inf_{0\leq\theta<2\pi}\|u_0-e^{i\theta }\psi_E\|_{H^1\cap
L^{\frac{\alpha_2+2}{\alpha_2+1}}}<\varepsilon(\psi_E)$$ the
solution of \eqref{u}-\eqref{ic} converges to a periodic solution
(close to $e^{-iEt}\psi_E(x)$) strongly in $L^p, 2<p<2N/(N-2),$
spaces and weakly in $H^1,$ as $t\rightarrow\pm\infty.$ The proof of this result is left for
another paper \cite{kz:asLs} because it involves a different
decomposition of the dynamics and a more delicate way to obtain the
linear estimates similar to the ones in Section \ref{se:lin}. It has
the advantage that it can be generalized to large periodic solutions
and the disadvantage that it only describes the evolution of initial
data in a conic like neighborhood (since $\varepsilon$ depends on
$\psi_E$) of the set of periodic solutions (center manifold) with
the zero solution removed. In fact the choice of $\varepsilon$ is
such that the solution stays away from zero, the point where the
center manifold fails to be $C^2$ smooth. Since the dynamics only
sees the $C^2$ smooth part of the center manifold a better
decomposition of the dynamics can be employed, see
\cite{km:asd3,kz:as2d2}, and convergence to a periodic solution
follows.

The main contribution of our paper is to describe the long time
evolution of all small initial data for rather general
nonlinearities including for the first time the subcritical ones
$\alpha_1<4/N,$ see Section \ref{se:main}. We accomplish this by
using a time-dependent projection of the solution of
\eqref{u}-\eqref{ic} onto the center manifold of periodic solutions described in
Section \ref{se:cm} and \eqref{ardef}. We first prove dispersive
estimates for the propagator of the linearized equation at the
time-dependent projection, see Section \ref{se:lin}. We then
estimate the error between the actual solution and its projection onto the center manifold
via a Duhamel principle with respect to the linearized operator at
the projection and a fixed point argument for the resulting integral
equation, see \eqref{sor}. Since the operator on its right hand side
contains no linear terms in the error, we are able to show that it is
contractive in appropriately chosen Banach spaces for a large
spectrum of nonlinearities $g.$ This is in contrast with the
approach in \cite{pw:cm,sw:mc2} where linear and nonlinear terms had
to be estimated at the same time, the linear ones requiring the use
of $L^2$ weighted (localized) estimates which applied to the
non-localized, nonlinear terms forced the assumption $\alpha_1>4/N.$
In the current approach, as in our previous 3D and 2D results, see
\cite{km:asd3,kz:as2d,kz:as2d2}, we completely separate estimates
for nonlinear terms from the ones for linear terms and use methods
tailored for each of them.

The most difficult part is to obtain dispersive estimates for the
propagator of the time-dependent linearized operator at the projections
onto the center manifold, see Section \ref{se:lin}. While estimates
for the Schr\" odinger group of operators:
\begin{equation}\label{introde}\|e^{-i(-\Delta+V(x))t}P_c\|_{L^{p'}\mapsto L^p}\le C_p |t|^{-N\left(\frac{1}{2}-\frac{1}{p}\right)},\qquad 1/p'+1/p=1,\ 2\le p\le\infty\end{equation}
are well known, see \cite{kn:jss} and references therein, they are
almost non-existent when the potential $V$ depends on time
$V=V(t,x).$ This is to be expected since the time-dependence of $V$
is the quantum mechanical analog of the parametric forcing in
ordinary differential equations and, in principle, can lead to very
different behavior compared to the time-independent case, see
\cite{sw:nh,kw:mbs,kw:diff,kirr:thesis,kn:CoSo,kn:CCLR2}. However,
in the absence of resonant phenomena one might expect similarities
between the two dynamics. Indeed, this is the case in \cite{rs:tdde}
which cannot be generalized to our situation mostly because of the
complex-valued potential, see \eqref{z} and \eqref{eq:dg}. To
overcome this issues we use smallness and localization of the time
dependent terms \eqref{F1loc}-\eqref{F1loc1} to first obtain
dispersive estimates in weighted (localized) norms, see Theorem
\ref{th:lw}. Then in Theorem \ref{th:lp} we rely on the
integrability in time of the group of operators generated by the
nearby time independent operator $-i(-\Delta+V(x)),$ see
\eqref{introde} with $p>2N/(N-2)$ and $t\ge 1,$ to remove the
weights and obtain dispersive estimates in non-localized $L^p$
norms. But the integrability in time for $t\ge 1$ comes at the cost
of a non-integrable singularity at $t=0$ which we remove by using
cancelations in highly oscillatory integrals. The method is similar
to the one we employed in \cite{km:asd3}, see also
\cite{kz:as2d,kz:as2d2} for an alternate way of dealing with the
singularity.

In a nutshell the results in this paper rely on shadowing the actual
solution of \eqref{u}-\eqref{ic} via a curve on the central manifold
of periodic solutions for \eqref{u}. Essential in showing that the
distance between the solution and its shadow goes to zero are the
new, apriori, dispersive estimates for the propagator of the
linearized equation along the shadowing curve. In this regard the
paper is an extension to higher dimensions, hence different linear
dispersive behavior, and to rougher nonlinearities of our previous
results in two respectively three space dimensions
\cite{kz:as2d,kz:as2d2,km:asd3}.

\bigskip

\noindent{\bf Notations:} $H=-\Delta+V;$

$L^p=\{f:\mathbb{R}^N\mapsto \mathbb{C}\  |\ f\ {\rm measurable\
and}\ \int_{\mathbb{R}^N}|f(x)|^pdx<\infty\},$
$\|f\|_p=\left(\int_{\mathbb{R}^N}|f(x)|^pdx\right)^{1/p}$ denotes
the standard norm in these spaces;

$<x>=(1+|x|^2)^{1/2},$ and for $\sigma\in\mathbb{R},$ $L^2_\sigma$
denotes the $L^2$ space with weight $<x>^{2\sigma},$ i.e. the space
of functions $f(x)$ such that $<x>^{\sigma}f(x)$ are square
integrable endowed with the norm
$\|f(x)\|_{L^2_\sigma}=\|<x>^{\sigma}f(x)\|_2;$

$\langle f,g\rangle =\int_{\mathbb{R}^N}\overline f(x)g(x)dx$ is the
scalar product in $L^2$ where $\overline z=$ the complex conjugate
of the complex number $z;$

$P_c$ is the projection on the continuous spectrum of $H$ in $L^2;$

$\hat u$ denotes the Fourier transform of the temperate distribution $u;$

$H^s,\ s\in\R$ denote the Sobolev spaces of temperate distributions
$u$ such that $(1+|\xi|^2)^{s/2}\hat u(\xi)\in L^2(\R^N)$ with norm
$\|u\|_{H^s}=\|(1+|\xi|^2)^{s/2}\hat u(\xi)\|_{L^2}.$ Note that for
$s=n, $ $n$ a natural number, this spaces coincide with the space of
measurable functions having all distributional partial derivatives
up to order $n$ in $L^2.$

\section{The Center Manifold}\label{se:cm}
The center manifold is formed by the collection of periodic solutions for (\ref{u}):
\begin{equation}\label{eq:per}
  u_E(t,x)=e^{-iEt}\psi_E(x)
\end{equation}
where $E\in\mathbb{R}$ and $0\not\equiv\psi_E\in H^2(\mathbb{R}^N)$
satisfy the time independent equation:
\begin{equation}\label{eq:ev}
[-\Delta+V]\psi_E+g(\psi_E)=E\psi_E
\end{equation}
Clearly the function constantly equal to zero is a solution of (\ref{eq:ev})
but
(iii) in the following hypotheses on the potential $V$ allows for a
bifurcation
with a nontrivial, one parameter family of solutions:

\bigskip

\noindent{\bf (H1)} Assume that
\begin{itemize}
  \item[(i)] $V(x)$ sutisfies the following properties:
   \begin{enumerate}
   \item $<x>^{\rho}V(x):H^\eta\rightarrow H^\eta$, for some $\rho>N+4$ and $\eta>0;$
   \item $\nabla V\in L^p(\mathbb{R}^N)$ for some $2\le p\le\infty$
   and $|\nabla V(x)|\rightarrow 0$ as $|x|\rightarrow\infty ;$
   \item the Fourier transform of $V$ is in $L^1.$
   \end{enumerate}
  \item[(ii)] $0$ is a regular point\footnote{see
 \cite[Definition 7]{ws:de2} or $M_\mu=\{0\}$ in relation (3.1) in \cite{mm:ae}}
 of the
spectrum of
  the linear operator $H=-\Delta+V$ acting on $L^2.$
  \item [(iii)]$H$ acting on $L^2$ has exactly one negative eigenvalue $E_0<0$ with corresponding normalized eigenvector $\psi_0.$ It is well known that $\psi_0(x)$ is  exponentially decaying as $|x|\rightarrow\infty,$ and can be chosen strictly positive.
\end{itemize}

\par\noindent Conditions (i) and (ii) guarantee the applicability of dispersive estimates in
\cite{mm:ae} and \cite{kn:jss} to the Schr\" odinger group
$e^{-iHt}P_c.$ Condition (i)2. implies certain regularity of the
nonlinear bound states while (i)3. allows us to use commutator type
inequalities, see \eqref{est:fm} and \cite[Theorem 5.2]{km:asd3}.
All these are needed to obtain estimates for the semigroup of
operators generated by our time dependent linearization,
 see Theorems \ref{th:lw} and \ref{th:lp} in
 section \ref{se:lin}. In particular (i)1. implies the local well
posedness in $H^1$ of the initial value problem
(\ref{u})-(\ref{ic}), see section \ref{se:main}.

By the standard bifurcation argument in Banach spaces \cite{ln:fa}
for (\ref{eq:ev}) at $E=E_0,$ condition (iii) guarantees existence
of nontrivial solutions. Moreover, these solutions can be organized
as a $C^1$ manifold (center manifold):

\begin{proposition}\label{pr:cm} There exist $\delta>0,$
the $C^1$ function
$$h:\{a\in\mathbb{C}\ :\ |a|\leq\delta\}\mapsto  H^2\cap L^2_\sigma
,\quad \sigma\in\R$$ and the $C^1$ function
$E:[-\delta,\delta]\mapsto\mathbb{R}$ such that for
$|E-E_0|\leq\delta$ and $|\langle\psi_0,\psi_E\rangle |\leq\delta,$
$\|\psi_E-\langle\psi_0,\psi_E\rangle\psi_0\|_{H^2}\leq\delta $ the
eigenvalue problem (\ref{eq:ev}) has a unique non-zero solution up
to multiplication with $e^{i\theta},\ \theta\in [0,2\pi),$ which can
be represented as a center manifold:
 \begin{equation}\label{eq:cm}
 \psi_E=a\psi_0+h(a),\ E=E(|a|), \quad \langle\psi_0,h(a)\rangle =0,\quad
 h(e^{i\theta}a)=e^{i\theta}h(a),\ {\rm for}\
|a|\leq\delta .\end{equation}
\end{proposition}
See \cite[section 2]{kz:as2d} for details.

\begin{remark}\label{rmk:reg} By regularity methods, see for example \cite[Theorem 8.1.1]{caz:bk}, one can show $\psi_E\in H^3\cap W^{2,p},\ 2\le p<\infty.$
Hence by Sobolev imbeddings both $\psi_E$ and $\nabla\psi_E$ are continuous and converge to zero as $|x|\rightarrow\infty.$ Moreover comparison techniques for elliptic equations, see \cite[Section 5.2]{km:asd3}, imply exponential decay, i.e. for each $0<A<-E$ there exists a constant $C_A>0$ such that:
$$\|e^{\sqrt{A}|x|}\psi_E(x)\|_{L^\infty}\leq C_A\|\psi_E(x)\|_{L^\infty} <\infty,\ {\rm and}\quad  \|e^{\sqrt{A}|x|}\nabla\psi_E(x)\|_{L^\infty}\leq C_A\|\nabla\psi_E(x)\|_{L^\infty}<\infty .$$
\end{remark}

\begin{remark}\label{rmk:pos} By variational methods, see for example \cite{rw:bs}, one can show that the real valued solutions of
\eqref{eq:ev} do not change sign. Then Harnack inequality for
$H^2\bigcap C(\mathbb{R}^N)$ solutions of \eqref{eq:ev} implies that
these real solution cannot take the zero value. Hence $\psi_E$ given
by \eqref{eq:cm} for $a\in\mathbb{R}$ is either strictly positive or
strictly negative.
\end{remark}

In section \ref{se:lin} we also need some smoothness for the
effective (linear) potential induced by the nonlinearity which
for the real valued bound states is:
\begin{equation}\label{eq:dg}Dg|_{\pe}[u+iv]=g'(\pe)u+i\frac{g(\pe)}{\pe}v,\qquad \psi_E > 0\end{equation}
while for an arbitrary bound state $\tilde\psi_{E}=e^{i\theta}\psi_E,\ \psi_E>0$ we have via the rotational symmetry of $g,$ see \eqref{gsym},
\begin{equation}\label{eq:dg1}Dg|_{\tilde\psi_E}[\beta]=e^{i\theta}\Dg [e^{-i\theta}\beta].\end{equation}

\noindent{\bf (H2)} Assume that for the positive solution of
\eqref{eq:ev} we have $\widehat{g'(\pe)},\
\widehat{\frac{g(\pe)}{\pe}}\in L^1(\mathbb{R}^N)$ where $\hat{f}$
stands for the Fourier transform of the function $f.$
\bigskip

\noindent In concrete cases the hypothesis may be checked directly
using the regularity of $\psi_E,$ the solution of an uniform
elliptic e-value problem. In general we can prove the following
result:

\begin{proposition}\label{pr:dg} If the following holds:

\bigskip

\noindent{\bf (H2')} $g$ restricted to reals
 has four derivatives except at zero and
 $$|g^{(m)}(s)|<C\left(\frac{1}{s^{m-1-\alpha_1}}+\frac{1}{s^{m-1-\alpha_2}}\right),\ {\rm for}\ m=3,4,\
 {\rm and}\ s>0.$$

\noindent Then for the positive solution of \eqref{eq:ev}, $\pe,$ we
have $\widehat{g'(\pe)}\in L^1$ and $\widehat{\frac{g(\pe)}{\pe}}\in
L^1$.
\end{proposition}
\textbf{Proof:}
$$\|\widehat{g'(\pe)}\|_{L^1}=\Big\|\frac{1}{(1+|\xi|^2)^{\frac{3}{2}}}(1+|\xi|^2)^{\frac{3}{2}}\widehat{g'(\pe)}\Big\|_{L^1}\leq\Big\|\underbrace{\frac{1}{(1+|\xi|^2)^{\frac{3}{2}}}}_{\in L^2(\R^N)}\Big\|_{L^2}\|\underbrace{(1+|\xi|^2)^{\frac{3}{2}}\widehat{g'(\pe)}}_{\in L^2\Leftrightarrow g'(\pe)\in H^3}\|_{L^2}$$
So it suffices to show that  $g'(\pe)\in H^3$ and
$\frac{g(\pe)}{\pe}\in H^3$. We have: %We will choose $A_1$ and $A_2$ such that $(3-\p)\sqrt{A_2}<\sqrt{A_1}$.
\begin{align}
\nabla g'(\pe)&=g''(\pe)\nabla\pe \nonumber \\
\Delta g'(\pe)&=g'''(\pe)|\nabla\pe|^2+g''(\pe)\Delta\pe \nonumber \\
              &=g'''(\pe)|\nabla\pe|^2+g''(\pe)[(V-E)\pe+g(\pe)] \nonumber \\
\nabla^3 g'(\pe)&=g^{(4)}(\pe)|\nabla\pe|^2\nabla\pe+3g'''(\pe)\nabla\pe\Delta\pe+g''(\pe)\nabla^3\pe \nonumber \\
                &=g^{(4)}(\pe)|\nabla\pe|^2\nabla\pe+3g'''(\pe)\nabla\pe[(V-E)\pe+g(\pe)]+g''(\pe)\nabla[(V-E)\pe+g(\pe)] \nonumber \end{align}
and, using \eqref{gest},
\begin{align}
|g'(\pe)|&\leq C(|\pe|^{\p_1}+|\pe|^{\p_2})\nonumber \\
|\nabla g'(\pe)|&\leq C|\nabla\pe|(\frac{1}{|\pe|^{1-\p_1}}+\frac{1}{|\pe|^{1-\p_2}}) \nonumber \\
|\Delta g'(\pe)|&\leq C|\nabla\pe|^2(\frac{1}{|\pe|^{2-\p_1}}+\frac{1}{|\pe|^{2-\p_2}})+C(|\pe|^{\p_1}+|\pe|^{\p_2})(|V|+|E|)+C(|\pe|^{2\p_1}+|\pe|^{2\p_2}) \nonumber \\
|\nabla^3 g'(\pe)|&\leq C|\nabla\pe|^3(\frac{1}{|\pe|^{3-\p_1}}+\frac{1}{|\pe|^{3-\p_2}})+2C|\nabla\pe|(|V|+|E|)(\frac{1}{|\pe|^{1-\p_1}}+\frac{1}{|\pe|^{1-\p_2}})\nonumber \\
&\quad+C|\nabla V|(|\pe|^{\p_1}+|\pe|^{\p_2})+2C|\nabla\pe|(|\pe|^{2\p_1-1}+|\pe|^{2\p_2-1}) \nonumber \end{align}

Similarly we get the estimates for $\frac{g(\pe)}{\pe}$ as follows:

\begin{align}
|\frac{g(\pe)}{\pe}|&\leq C(|\pe|^{\p_1}+|\pe|^{\p_2})\nonumber \\
|\nabla \frac{g(\pe)}{\pe}|&\leq C|\nabla\pe|(\frac{1}{|\pe|^{1-\p_1}}+\frac{1}{|\pe|^{1-\p_2}}) \nonumber \\
|\Delta \frac{g(\pe)}{\pe}|&\leq C|\nabla\pe|^2(\frac{1}{|\pe|^{2-\p_1}}+\frac{1}{|\pe|^{2-\p_2}})+C(|\pe|^{\p_1}+|\pe|^{\p_2})(|V|+|E|)+C(|\pe|^{2\p_1}+|\pe|^{2\p_2}) \nonumber \\
|\nabla^3 \frac{g(\pe)}{\pe}|&\leq C|\nabla\pe|^3(\frac{1}{|\pe|^{3-\p_1}}+\frac{1}{|\pe|^{3-\p_2}})+2C|\nabla\pe|(|V|+|E|)(\frac{1}{|\pe|^{1-\p_1}}+\frac{1}{|\pe|^{1-\p_2}})\nonumber \\
&\quad+C|\nabla V|(|\pe|^{\p_1}+|\pe|^{\p_2})+2C|\nabla\pe|(|\pe|^{2\p_1-1}+|\pe|^{2\p_2-1}) \nonumber \end{align}

Now, we will use the following bounds for $\pe$ and $\nabla\pe$, see
\cite[Section 5.2]{km:asd3}. For any $0<A<-E<A_2$ and any $0<A_1<-E$ there exist the constants $C_A,\ C_{A_1},\ C_{A_2}>0$ such that:
$$C_{A_2}e^{-\sqrt{A_2}|x|}\leq\pe(x)\leq C_A e^{-\sqrt{A}|x|},\quad
{\rm for\ all}\ x\in\R^N,$$ $$|\nabla\pe(x)|\leq
C_{A_1}e^{-\sqrt{A_1}|x|},\quad {\rm for\ all}\ x\in\R^N.$$ Choosing
$A_1$ and $A_2$ such that $(3-\p_1)\sqrt{A_2}<3\sqrt{A_1}$, we
obtain $g'(\pe)\in H^3$ and $\frac{g(\pe)}{\pe}\in H^3$. The
proposition is now completely proven. $\Box$

\section{Main Results}\label{se:main}
\begin{theorem}\label{mt}
Assume that the nonlinear term in \eqref{u} satisfies \eqref{gest},
\eqref{alpha12} and \eqref{gsym}. In addition assume that hypotheses
(H1) and either (H2) or (H2') hold. Then there exists an
$\varepsilon_0$ such that for all initial conditions $u_0(x)$
satisfying
$$\max\{\|u_0\|_{L^{p_2'}},\|u_0\|_{H^1}\}\leq\varepsilon_0,\qquad p_2=2+\p_2,\qquad \frac{1}{p_2'}+\frac{1}{p_2}=1$$
the initial value problem (\ref{u})-(\ref{ic}) is globally
well-posed in $H^1$ and the solution converges to the center
manifold.

More precisely, there exist a $C^1$ function
$a:\mathbb{R}\mapsto\mathbb{C}$ such that, for all $t\in\R$ we have:
\begin{equation}
u(t,x)=\underbrace{a(t)\psi_0(x)+h(a(t))}_{\psi_{E(t)}}+r(t,x)
\label{dc} \end{equation} where $\psi_{E(t)}$ is on the central
manifold (i.e it is a ground state) , see Proposition \ref{pr:cm}.
Moreover, for all $t\in\R $ $r(t,x)$ satisfies the following decay
estimates:
 \begin{eqnarray}
 \|r(t)\|_{L^2}&\leq &C_0(\alpha_1,\alpha_2)\varepsilon_0\nonumber\\
 \|r(t)\|_{L^{p_1}}&\leq &
 C_1(\alpha_1,\alpha_2)\frac{\varepsilon_0}{(1+|t|)^{N(\frac{1}{2}-\frac{1}{p_1})}},\quad p_1=2+\alpha_1\nonumber
 \end{eqnarray}
and, for $p_2=2+\alpha_2:$
\begin{enumerate}
\item [(i)] if $\alpha_1\ge\frac{4}{N}$ or $p_2<\frac{2N}{2+N-N\p_1}$ then $\|r(t)\|_{L^{p_2}}\leq
 C_2(\alpha_1,\alpha_2)\frac{\varepsilon_0}{(1+|t|)^{N(\frac{1}{2}-\frac{1}{p_2})}}
 $
\item [(ii)] if $\alpha_1<\frac{4}{N}$ and $p_2=\frac{2N}{2+N-N\p_1}$ then $\|r(t)\|_{L^{p_2}}\leq
 C_2(\alpha_1,\alpha_2)\varepsilon_0\frac{\log(2+|t|)}{(1+|t|)^{N(\frac{1}{2}-\frac{1}{p_2})}}
 $
\item [(iii)] if $\alpha_1<\frac{4}{N}$ and $p_2>\frac{2N}{2+N-N\p_1}$ then
 $\|r(t)\|_{L^{p_2}}\leq
 C_2(\alpha_1,\alpha_2)\frac{\varepsilon_0}{(1+|t|)^{\frac{N\alpha_1}{2}-1}}
 $
\end{enumerate}
where the constants $C_0,\ C_1$ and $C_2$ are independent of
$\varepsilon_0$.
\end{theorem}

Before proving the theorem let us note that (\ref{dc}) decomposes
the evolution of the solution of (\ref{u})-(\ref{ic}) into an
evolution on the center manifold $\psi_{E(t)}$ and the ``distance"
from the center manifold $r(t)$. The estimates on the latter show
collapse of solution onto the center manifold. A more
precise decay of the `radiative" part, $r(t),$ in different $L^p$
spaces is given in the following Corollary. It shows same decay as
the solution of the free Schr\" odinger equation up to the threshold
$p=\frac{2N}{2+N-N\alpha_1}$ where it saturates:
\begin{corollary}\label{cor} Consider $2\leq p<\frac{2N}{N-2},$ and $1/p'+1/p=1.$
Under the hypotheses of Theorem \ref{mt}, assuming also $u_0\in L^{p'},$  we have the following  decay estimates:

if $\alpha_1\geq 4/N$ then
$$\|r(t)\|_{L^p}\leq\frac{C(p)\max\{\|u_0\|_{L^{p'}}, \varepsilon_0\}}{(1+|t|)^{N\ip}},\quad \textrm{for all}\ 2\leq p< \frac{2N}{N-2},$$

otherwise
\begin{align}
\|r(t)\|_{L^p}&\leq\frac{C(p) \max\{\|u_0\|_{L^{p'}}, \varepsilon_0\}}{(1+|t|)^{N\ip}},\quad  \textrm{if $p_1<p<\frac{2N}{2+N-N\alpha_1}$}\nonumber \\
\|r(t)\|_{L^p}&\leq\frac{C(p)\log(2+|t|)\max\{\|u_0\|_{L^{p'}}, \varepsilon_0\}}{(1+|t|)^{N\ip}},\quad  \textrm{if $p=\frac{2N}{2+N-N\alpha_1}$} \nonumber \\
\|r(t)\|_{L^p}&\leq\frac{C(p)\max\{\|u_0\|_{L^{p'}},
\varepsilon_0\}}{(1+|t|)^{\frac{N\p_1}{2}-1}},\quad \textrm{if
$\frac{2N}{2+N-N\alpha_1}<p<\frac{2N}{N-2}.$}\nonumber \end{align}
\end{corollary}

The dynamics on the
center manifold is determined by equation (\ref{a}) below. In supercritical regimes $\alpha_1>4/N$ we can actually show that it converges to a certain periodic orbit:
\begin{corollary}\label{cor:as} Under the hypotheses of Theorem \ref{mt}, assuming also $\alpha_2\ge\alpha_1>4/N$ and $u_0\in L^{p'}(\R^N)$ for some $1\le p'<2N/(N+2),$ we have in addition to the conclusion of Theorem \ref{mt} that there exists $a_{\pm\infty}\in\R$ such that $\lim_{t\rightarrow\pm\infty}|a(t)|=a_{\pm\infty}.$ Moreover, if we denote by $e^{-itE_{\pm\infty}}\psi_{E_{\pm\infty}}(x)$ the periodic solutions of \eqref{u} corresponding to the parameters $a_{\pm\infty}$ on the center manifold:
$$\psi_{E_{\pm\infty}}=a_{\pm\infty}\psi_0+h(a_{\pm\infty}),\qquad E_{\pm\infty}=E(a_{\pm\infty}),$$
see Proposition \ref{pr:cm}, then there exists a $C^1$ function $\theta:\R\mapsto\R$ such that:
$$\lim_{t\rightarrow\pm\infty}\|\psi_{E(t)}-e^{-it(E_{\pm\infty}+\theta(t))}\psi_{E_{\pm\infty}}\|_{H^2\cap L^2_\sigma}=0,\quad \lim_{|t|\rightarrow\infty}\theta(t)=0,$$
where $\psi_{E(t)}$ is the component on the central manifold of the actual solution of \eqref{u}-\eqref{ic}, see \eqref{dc}.
\end{corollary}
The corollary extends the results in \cite{sw:mc2,pw:cm} to nonlinearities in more general form than pure power and to initial data that are not necessarily localized, i.e. in $L^2_\sigma(\R^N),\ \sigma > N.$ As we shall see in Remark \ref{rmk:whysc} the supercriticallity restriction $\alpha_1 >4/N$ comes from the fact that the dispersive part $r(t)$ appears linearly in the equation \eqref{a} for the central manifold parameter $a.$ Note that in \cite{kz:as2d2,km:asd3} we use an improved decomposition of the type \eqref{dc} in which the equation on the central manifold corresponding to \eqref{a} contains only quadratic and higher order terms in $r(t).$ While this decomposition allows us to show convergence to a periodic solution even for subcritical regimes $\alpha_1 \le 4/N$ in dimensions $N=2,3,$ it requires the central manifold to be $C^2,$ i.e. $\alpha_1\ge 1.$ However our central manifold is $C^2$ except at zero, so if the initial data is chosen in an appropriate manner such that the dynamics stays away from zero then Corollary \ref{cor:as} can be proven even in subcritical regimes, see \cite{kz:asLs}.
\begin{remark}\label{rmk:as} In conclusion the approach in \cite{kz:as2d2,km:asd3} would allow us to obtain Corollary \ref{cor:as} for the critical regime $\alpha_1=1$ in dimension $N=4,$ but it would require the stronger hypothesis $\alpha_1\ge 1$ in dimension $N=5.$ However for initial data in a conic like neighborhood of the manifold of non-zero periodic solutions of \eqref{u} the conclusion of Corollary \ref{cor:as} is valid for all $\alpha_0(N)<\alpha_1\le \alpha_2 <4/(N-2),$ see \cite{kz:asLs}.
\end{remark}

We now proceed with the proofs:

 \textbf{Proof of Theorem \ref{mt}} It is well known
that under hypothesis (H1)(i) the initial value problem (1)-(2) is
locally well posed in the energy space $H^1$ and its $L^2$ norm is
conserved, see for example \cite[Cor. 4.3.3 at p. 92]{caz:bk}.
Global well posedness follows via energy estimates from
$\|u_0\|_{H^1}$ small, see \cite[Remark 6.1.5 at p. 165]{caz:bk}.

In particular we can define
$$a(t)=\<\psi_0,u(t)\>, \quad \text{for} {\ } t\in\R$$ Cauchy-Schwarz inequality implies
 \begin{equation}\label{abound}
 |a(t)|\leq\|u(t)\|_{L^2}\|\psi_0\|_{L^2}=\|u_0\|_{L^2}\leq\varepsilon_0, \quad \text{for all}{\ } t\in\R .
 \end{equation}
Hence, if we choose $\varepsilon_0\leq\delta$ we can define $h(a(t))$, $t\in\R$, see Proposition \ref{pr:cm}. We then obtain (\ref{dc}) where
\begin{equation}\label{ardef}r(t)=u(t)-a(t)\psi_0-h(a(t)), \quad a(t)=\<\psi_0,u(t)\>, \quad \<\psi_0, r(t)\>\equiv0.\end{equation}

The solution is now described by the $C^1$ scalar function
$a(t)\in\C$ and $r(t)\in C(\R,H^1)\cap C^1(\R,H^{-1})$. To obtain
their equations we plug in (\ref{dc}) into (\ref{u}) and project onto $\psi_0$ and its orthogonal complement:
\begin{align}
i\frac{da}{dt}&=E(|a(t)|)a(t)+\<\psi_0, g(\psi_{E(t)}+r(t))-g(\psi_{E(t)})\> \label{a} \\
i\frac{dr}{dt}&=Hr(t)+P_c[g(\psi_{E(t)}+r(t))-g(\psi_{E(t)})]+iDh|_{a(t)}i\<\psi_0,
g(\psi_{E(t)}+r(t))-g(\psi_{E(t)})\> \label{r}
\end{align}
where $H=-\Delta+V$.
\begin{equation}
g(\psi_{E(t)}+r(t))-g(\psi_{E(t)})=\underbrace{\frac{d}{d\varepsilon}g(\psi_{E(t)}+\varepsilon r(t))|_{\varepsilon=0}}_{Dg|_{\psi_{E(t)}}[r(t)]}+F_2(\psi_{E(t)},r(t)) \label{g} \end{equation}
where $F_2$ contains only higher order corrections in $r$ to the linear term $\Dg [r].$
The linear part of (\ref{r}) is:
\begin{align}
i\frac{dz}{dt}&=Hz(t)+P_cDg|_{\psi_{E(t)}} [z(t)]+iDh|_{a(t)}i\<\psi_0,Dg|_{\psi_{E(t)}}[z(t)]\> \label{z} \\
z(s)&=v \nonumber
\end{align}

Define $\Omega(t,s)v=z(t)$. Then using Duhamel's principle (\ref{r}) becomes
\begin{equation}
r(t)=\om(t,0)r(0)-\int_0^t
\om(t,s)[P_ciF_2(\psi_{E(s)},r(s))-Dh|_{a(s)}\<\po,iF_2(\psi_{E(s)},r(s))\>]ds.
\label{sor}
\end{equation}
In order to apply the linear estimates of Section \ref{se:lin} to
$\om(t,s),$ we fix $\sigma>N/2$ and
$\frac{2N}{N-2}<q_2<\frac{2N}{N-4},$ then we consider the
$\varepsilon_1(q_2)>0$ given by Theorem \ref{th:lw} and choose
$\varepsilon_0>0$ in the hypotheses such that
 \begin{equation}\label{eq:lee0}
 \|\<x\>^{4\sigma}\psi_{E(t)}(x)\|_{L^\infty(\R^N)}\leq\varepsilon_1,\qquad {\rm for\ all}\ t\in\R .
 \end{equation}
The latter is possible because $E(t)=E(|a(t)|)$ and $E$ is a $C^1$ function from the compact interval $[-\delta,\delta]$ to the real numbers with $E(0)=E_0<0.$ So, there exists $\varepsilon_0$ such that $E(|a|)\le E_0/2 <0$ for all $|a|\leq\varepsilon_0.$ \eqref{eq:lee0} now follows from the exponential decay estimates in Remark \ref{rmk:reg} and the observation $\|\psi_{E(t)}\|_{L^\infty}\leq C\varepsilon_0,$ for some constant $C>0.$ This is a consequence of Sobolev imbeddings and $\|\psi_{E(t)}\|_{H^3}\leq C\varepsilon_0$ which follows from $\psi_{E(t)}=a(t)\psi_0+h(a(t)),$ the norm $\|\psi_0\|_{H^2}$ is fixed, $|a(t)|\leq \varepsilon_0$ for all $t\in\R,$ see \eqref{abound}, and $h$ is a $C^1$ function on the compact ball of radius $\delta$ in complex plane to $H^2.$ Hence, there exists a constant $C>0$ such that $\|\psi_{E(t)}\|_{H^2}\leq C\varepsilon_0$ for all $t\in\R ,$ and, by regularity arguments, see for example \cite[Theorem 8.1.1]{caz:bk} we get the same estimate for the $H^3$ norm with a possible larger constant $C.$

Consider now the nonlinear operator in (\ref{sor}):
$$N(u)(t)=\int_0^t \om(t,s)[P_ciF_2(\psi_{E(s)},u(s))-Dh|_{a(s)}\<\po,iF_2(\psi_{E(s)},u(s))\>]ds$$
In order to apply a contraction mapping argument for (\ref{sor}) we
use the following Banach spaces. Let $p_1=2+\p_1,$  $p_2=2+\p_2$,
then
 \begin{eqnarray}
 Y_i&=&\{u\in C\left(\R,\ L^2(\R^N)\cap L^{p_1}(\R^N)\cap L^{p_2}(\R^N)\right):\nonumber\\
 &&\sup_{t\in\R}
(1+|t|)^{N(\frac{1}{2}-\frac{1}{p_1})}\|u(t)\|_{L^{p_1}}<\infty,\ \sup_{t\in\R}
\frac{(1+|t|)^{n_i}}{[\log(2+|t|)]^{m_i}}\|u(t)\|_{L^{p_2}}<\infty,\
\sup_{t\in\R} \|u(t)\|_{L^2}<\infty \}\nonumber\end{eqnarray} endowed with the norm $$
\|u\|_{Y_i}=\max\{\sup_{t\in\R}
(1+|t|)^{N(\frac{1}{2}-\frac{1}{p_1})}\|u\|_{L^{p_1}},\ \sup_{t\in\R}
\frac{(1+|t|)^{n_i}}{[\log(2+|t|)]^{m_i}}\|u\|_{L^{p_2}},\ \sup_{t\in\R}
\|u\|_{L^2} \}$$ for $i=1,2,3$, where
$n_{1,2}=N(\frac{1}{2}-\frac{1}{p_2})$, $n_3=\frac{N\p_1}{2}-1$,
$m_1=m_3=0$ and $m_2=1$.

\begin{lemma}\label{lm}
Consider the cases:
$$ 1.\ N\left(\frac{\p_1-1}{2}+\frac{1}{2+\p_2}\right)>1 ; \quad 2.\ N\left(\frac{\p_1-1}{2}+\frac{1}{2+\p_2}\right)=1; \quad 3.\ N\left(\frac{\p_1-1}{2}+\frac{1}{2+\p_2}\right)<1;$$
and assume that \eqref{eq:lee0} holds for some $\sigma>N/2.$ Then,
for each case number i, $N : Y_i\rightarrow Y_i$ is well defined,
and locally Lipschitz, i.e. there exists $\tilde{C_i}>0$, such that
$$\|Nu_1
-Nu_2\|_{Y_i}\leq\tilde{C_i}(\|u_1\|_{Y_i}+\|u_2\|_{Y_i}+\|u_1\|_{Y_i}^{\p_1}+\|u_2\|_{Y_i}^{\p_1}+\|u_1\|_{Y_i}^{\p_2}+\|u_2\|_{Y_i}^{\p_2})\|u_1
-u_2\|_{Y_i}. $$
\end{lemma}
Note that the Lemma finishes the proof of Theorem \ref{mt}. Indeed,
$\alpha_1$ and $\alpha_2$ determine in which case $i=1,2,3$ we are.
But in all of them, if we denote:$$v=\Omega(t,0)r(0),$$ then:
$$\|v\|_{Y_i}\leq C_0\|r(0)\|_{L^{p_2'}\cap L^2}\leq C_0(1+\|\psi_0\|_{L^{p'_2}})\underbrace{\|u_0\|_{L^{p_2'}\cap L^2}}_{\leq\varepsilon_0},$$ where
$C_0=\max\{C_2,C_{p_2}\}$, see Theorem \ref{th:lp}. We choose
$\varepsilon_0$ in the hypotheses of theorem \ref{mt}, such that
$R=2\|v\|_{Y_i}$ satisfies
$$Lip=2\tilde{C_i}(R+R^{\p_1}+R^{\p_2})<1,$$ where $\tilde{C_i}$ is given by the above Lemma. In this case the integral
operator given by the right hand side of the \eqref{sor}:
$$K(r)=v-N(r)$$
leaves the closed ball $B(0,R)=\{z\in Y_i :\|z\|_{Y_i}\leq R\}$
invariant and it is a contraction on it with Lipschitz constant
$Lip$. Consequently the equation \eqref{sor} has a unique solution
in $B(0,R)$. In particular, $r(t)$ satisfies the $L^p$ estimates
claimed in the Theorem \ref{mt}. We now have two solutions of
\eqref{sor}, one in $C(\R,H^1)$ from classical well posedness theory
and one in $C(\R,L^2\cap L^{p_1}\cap L^{p_2})$, $p_1=2+\p_1$,
$p_2=2+\p_2$ from the above argument. Using uniqueness and the
continuous embedding of $H^1$ in $L^2\cap L^{p_1}\cap L^{p_2}$, we
infer that the solutions must coincide. Therefore, the time decaying
estimates in the space $Y_{i}$ hold also for the $H^1$ solution.

\textbf{Proof of Lemma \ref{lm}} Let us first consider the difference
\begin{align}
F_2(\pe,u_1)-&F_2(\pe,u_2)=g(\pe+u_1)-g(\pe+u_2)-\Dg[u_1]+\Dg[u_2] \nonumber \\
&=\int_0^1\Big(Dg(\pe+u_2+\tau(u_1-u_2))[u_1-u_2]-Dg(\pe)[u_1-u_2]\Big)d\tau \nonumber \\
&=\int_0^1 \int_0^1 D^2g(\pe+s(u_2+\tau(u_1-u_2)))[u_2+\tau(u_1-u_2)][u_1-u_2]ds d\tau \label{f2dif}
\end{align}
For a fixed $(t,x)\in\R\times\R^N$, the inside integral is a line integral connecting
$z_1=\pe(x)$ and $z_2=\pe(x)+u_2(t,x)+\tau(u_1(t,x)-u_2(t,x)).$ Note that, from
\eqref{gest} and \eqref{alpha12}, $D^2g(\cdot)$ is integrable along
any segment in the complex plane. Using the equivariance under
rotations \eqref{gsym} we will reduce the line integral to a
horizontal segment. First let us observe the behavior of $Dg$ and
$D^2g$ under rotation:
$$Dg(z)[e^{i\theta}\beta]=\lim_{\varepsilon\rightarrow0}\frac{g(z+\varepsilon e^{i\theta}\beta)-g(z)}{\varepsilon}=e^{i\theta}\lim_{\varepsilon\rightarrow0}\frac{g(e^{-i\theta}z+\varepsilon\beta)-g(e^{-i\theta}z)}{\varepsilon}=e^{i\theta}Dg(e^{-i\theta}z)[\beta]$$
and
\begin{align}D^2g(e^{i\theta}z)[\alpha][\beta]&=\lim_{\varepsilon\rightarrow0}\frac{Dg(e^{i\theta}z+\varepsilon\alpha)[\beta]-Dg(e^{i\theta}z)[\beta]}{\varepsilon}\nonumber \\
&=e^{i\theta}\lim_{\varepsilon\rightarrow0}\frac{Dg(z+\varepsilon e^{-i\theta}\alpha)[e^{i\theta}\beta]-Dg(z)[e^{i\theta}\beta]}{\varepsilon}\nonumber \\
&=e^{i\theta}D^2
g(z)[e^{-i\theta}\alpha][e^{-i\theta}\beta]\nonumber \end{align}

From \eqref{gest} we have $|D^2g(z)|\leq
C(|z|^{\alpha_1-1}+|z|^{\alpha_2-1})$. For $z_2\not= z_1\in\C $ let
$z_2-z_1=|z_2-z_1|e^{i\theta}.$ Then:
\begin{align}|(Dg(z_2)-Dg(z_1))[\beta]|&\leq \left|\int_0^1 D^2 g(z_1+s(z_2-z_1))[z_2-z_1][\beta]ds\right|\nonumber \\
&\leq \left|\int_0^1 e^{i\theta}D^2 g(e^{-i\theta}z_1+se^{-i\theta}(z_2-z_1))[e^{-i\theta}(z_2-z_1)][e^{-i\theta}\beta ] ds\right|\nonumber \\
&\leq\int_0^1  C\left(|e^{-i\theta}z_1+s|z_2-z_1|\ |^{\p_1-1}+|e^{-i\theta}z_1+s|z_2-z_1|\ |^{\p_2-1}\right)|z_2-z_1||\beta|ds \nonumber
\end{align}
Now for $0<\p<1$ we have
 $$
 \int_0^1  |e^{-i\theta}z_1+s|z_2-z_1|\ |^{\p-1}|z_2-z_1|ds\leq  \int_0^1  |\Re
 z_1+s|z_2-z_1|\
 |^{\p-1}|z_2-z_1|ds
\leq |z_2-z_1|^\p ,$$
 while for $\alpha\geq 1$,
 $$
 \int_0^1  |e^{-i\theta}z_1+s|z_2-z_1|\ |^{\p-1}ds\leq \max_{t\in[0,1]}|e^{-i\theta}z_1+t|z_2-z_1|\ |^{\p-1}\int_0^1 1ds
\leq \max\{|z_1|^{\p -1},\ |z_2|^{\p-1}\}. $$ Hence, depending on
the powers $\p_1$ and $\p_2$ we have
\begin{align}
|F_2(\pe,u_1)-F_2(\pe,u_2)|\leq C\big[&\underbrace{(|u_1|^{\p_1}+|u_2|^{\p_1})|u_1-u_2|}_{A_1}+\underbrace{(|u_1|^{\p_2}+|u_2|^{\p_2})|u_1-u_2|}_{A_2} \label{def:a13} \\
 &+\underbrace{|\pe|^{\p_2-1}(|u_1|+|u_2|)|u_1-u_2|}_{A_3}+\underbrace{|\pe|^{\p_1-1}(|u_1|+|u_2|)|u_1-u_2|}_{A_4}\big] \nonumber \end{align}
where the $A_3$ term appears only when $\p_2> 1$ and the $A_4$ term
appears only when $\p_1> 1$.

Now let us consider the difference $Nu_1-Nu_2$

\begin{align}
(Nu_1-Nu_2)(t)=-i\int_0^t &\Omega(t,s)P_c \big[F_2(\pe,u_1)-F_2(\pe,u_2)\big]ds \nonumber \\
&+\int_0^t \Omega(t,s)Dh|_{a(s)}i
\<\po,F_2(\pe,u_1)-F_2(\pe,u_2)\>ds \label{eq:nu12} \end{align}

\textbf{$\bullet\ L^{p_2}$ Estimate :}
 \begin{align}
    \|Nu_1-Nu_2\|_{L^{p_2}}&\leq\int_0^{|t|} \|\om(t,s)\|_{L^{p'_2}\rightarrow L^{p_2}}(\|A_1\|_{L^{p'_2}}+\|A_2\|_{L^{p'_2}}+\|A_{3,4}\|_{L^{p_2'}})ds \nonumber \\
    &+\int_0^{|t|} \|\om(t,s)\|_{\Ls\rightarrow L^{p_2}}\|Dh|_{a(s)}\|_{\Ls}(\underbrace{|\<\po,A_1\>|}_{B_1}+\underbrace{|\<\po,A_2\>|}_{B_2}+\underbrace{|\<\po,A_{3,4}\>|}_{B_{3,4}})ds \nonumber \end{align}
The term $A_1$ satisfies via H\" older inequality:
 \begin{eqnarray} \|(|u_1|^{\p_1}+|u_2|^{\p_1})|u_1-u_2|\|_{L^{p'_2}}&\leq &\left(\|u_1\|_{L^{p_1}}^{\theta\p_1}\|u_1\|_{L^{2}}^{(1-\theta)\p_1}+\|u_2\|_{L^{p_1}}^{\theta\p_1}\|u_2\|_{L^{2}}^{(1-\theta)\p_1}\right) \nonumber\\
 &&\times\|u_1-u_2\|^\theta_{L^{p_1}}\|u_1-u_2\|^{1-\theta}_{L^{2}}\nonumber
 \end{eqnarray}
 where $\frac{1}{p'_2}=(1+\p_1)(\frac{1-\theta}{2}+\frac{\theta}{p_1}),\ 0\le\theta\le 1. $ Using Theorem \ref{th:lp} we get
\begin{align}
    \int_0^{|t|} & \|\om(t,s)\|_{L^{p'_2}\rightarrow L^{p_2}}\|A_1\|_{L^{p'_2}}ds \nonumber \\
    &\leq\int_0^{|t|} \frac{C(p_2)}{|t-s|^{N\ipt}}\cdot\frac{(\|u_1\|_{Y_i}^{\p_1}+\|u_2\|_{Y_i}^{\p_1})\|u_1-u_2\|_{Y_i}}{(1+|s|)^{N(\frac{\p_1-1}{2}+\frac{1}{p_2})}} ds \nonumber \\
&\leq\frac{C(p_2)C_3[\log(2+|t|)]^{m_i}}{(1+|t|)^{n_i}}(\|u_1\|_{Y_i}^{\p_1}+\|u_2\|_{Y_i}^{\p_1})\|u_1-u_2\|_{Y_i}
\nonumber \end{align} where the different decay rates $n_i$ depend
on the case number in the hypotheses of this Lemma:

 1. corresponds to
$N(\frac{\p_1-1}{2}+\frac{1}{p_2})>1$, and in this case
 $$C_3=\sup_{t\in\R}
(1+|t|)^{N\ipt}\int_0^{|t|}
\frac{1}{|t-s|^{N\ipt}}\frac{1}{(1+|s|)^{N(\frac{\p_1-1}{2}+\frac{1}{p_2})}}ds<\infty;$$

 2. corresponds to $N(\frac{\p_1-1}{2}+\frac{1}{p_2})=1$, and in this case
 $$C_3=\sup_{t\in\R}
\frac{(1+|t|)^{N\ipt}}{\log(2+|t|)}\int_0^{|t|}
\frac{1}{|t-s|^{N\ipt}}\frac{1}{(1+|s|)}ds<\infty;$$

 3. corresponds to $N(\frac{\p_1-1}{2}+\frac{1}{p_2})<1$, and in this case
 $$C_3=\sup_{t\in\R}
(1+|t|)^{\frac{N\p_1}{2}-1}\int_0^{|t|}
\frac{1}{|t-s|^{N\ipt}}\frac{1}{(1+|s|)^{N(\frac{\p_1-1}{2}+\frac{1}{p_2})}}ds<\infty.$$

To estimate the term containing $A_2$, observe that via H\" older inequality
 $$\|(|u_1|^{\p_2}+|u_2|^{\p_2})|u_1-u_2|\|_{L^{p'_2}}\leq\big(\|u_1\|_{L^{p_2}}^{\p_2}+\|u_2\|_{L^{p_2}}^{\p_2}\big)\|u_1-u_2\|_{L^{p_2}}$$
since $\frac{1}{p'_2}=\frac{1+\p_2}{p_2}$. Again, using Theorem
\ref{th:lp}, we have
    \begin{align}
    \int_0^{|t|} & \|\om(t,s)\|_{L^{p'_2}\rightarrow L^{p_2}}\|A_2\|_{L^{p'_2}}ds \label{a2lp2} \\
    &\leq\int_0^{|t|} \frac{C(p_2)}{|t-s|^{N\ipt}}\cdot\frac{[\log(2+|s|)]^{(1+\p_2)m_i}}{(1+|s|)^{(1+\p_2)n_i}}(\|u_1\|_{Y_i}^{\p_2}+\|u_2\|_{Y_i}^{\p_2})\|u_1-u_2\|_{Y_i} ds \nonumber \\
&\leq\frac{C(p_2)C_4C_5[\log(2+|t|)]^{m_i}}{(1+|t|)^{n_i}}(\|u_1\|_{Y_i}^{\p_2}+\|u_2\|_{Y_i}^{\p_2})\|u_1-u_2\|_{Y_i} \nonumber \end{align}
    where $C_5=\sup_t \frac{(1+|t|)^{n_i}}{[\log(2+|t|)]^{m_i}}\int_0^{|t|} \frac{[\log(2+|s|)]^{(1+\p_2)m_i}ds}{|t-s|^{N\ipt}(1+|s|)^{(1+\p_2)n_i}}<\infty$ since $(1+\p_2)n_i >1$.

As for the $A_{3,4}$ terms note that they only appear when
$\alpha_2> 1$ respectively $\alpha_1> 1.$ To estimate them observe
that
$$\|(|\pe|^{\p_1-1}+|\pe|^{\p_2-1})(|u_1|+|u_2|)|u_1-u_2|\|_{L^{p'_2}}\leq\||\pe|^{\p_1-1}+|\pe|^{\p_2-1}\|_{L^{\beta}}(\|u_1\|_{L^{p_2}}+\|u_2\|_{L^{p_2}})\|u_1-u_2\|_{L^{p_2}}$$
where $\frac{1}{\beta}+\frac{2}{p_2}=\frac{1}{p'_2}.$ By Theorem
\ref{th:lp} we have for each case number $i$ and $u_1,\ u_2\in Y_i:$
\begin{align}
\int_0^{|t|} & \|\om(t,s)\|_{L^{p'_2}\rightarrow L^{p_2}}\|A_{3,4}\|_{L^{p'_2}}ds \label{a34lp2} \\
    &\leq\int_0^{|t|} \frac{C(p_2)}{|t-s|^{N\ipt}}\||\pe|^{\p_1-1}+|\pe|^{\p_2-1}\|_{L^{\beta}}\frac{[\log(2+|s|)]^{2m_i}}{(1+|s|)^{2n_i}}(\|u_1\|_{Y_i}+\|u_2\|_{Y_i})\|u_1-u_2\|_{Y_i} ds \nonumber \\
    &\leq\frac{C(p_2)C_1 C_2}{(1+|t|)^{N\ipt}}(\|u_1\|_{Y_i}+\|u_2\|_{Y_i})\|u_1-u_2\|_{Y_i} \nonumber \end{align}
where $C_1=\sup_t \||\pe|^{\p_1-1}+|\pe|^{\p_2-1}\|_{L^{\beta}}$ and
$C_2=\sup_t \frac{(1+|t|)^{n_i}}{[\log(2+|t|)]^{m_i}}\int_0^{|t|}
\frac{[\log(2+|s|)]^{2m_i}ds}{|t-s|^{N\ipt}(1+|s|)^{2n_i}}<\infty$
since $2n_i >1$ (for $p_2=\alpha_2+2> 3,$ and $\alpha_1$ satisfying
\eqref{alpha12}). The uniform bounds in $t\in\R$ for
$\|\pe\|_{L^{(\alpha_j-1)\beta}}^{\alpha_j-1},\ j=1,2$ follow from
\eqref{eq:lee0}.

For $B$ terms we have:
$$|B_1|\leq\|\po\|_{L^{p_1}}\|A_1\|_{L^{p'_1}},\quad
|B_2|\leq\|\po\|_{L^{p_2}}\|A_2\|_{L^{p'_2}},\quad \textrm{and}\quad
|B_{3,4}|\leq\|\po\|_{L^{p_2}}\|A_{3,4}\|_{L^{p'_2}}.$$ Note that
\begin{equation}\label{a1p1}\|(|u_1|^{\p_1}+|u_2|^{\p_1})|u_1-u_2|\|_{L^{p'_1}}\leq\big(\|u_1\|_{L^{p_1}}^{\p_1}+\|u_2\|_{L^{p_1}}^{\p_1}\big)\|u_1-u_2\|_{L^{p_1}}\end{equation}
since $\frac{1}{p'_1}=\frac{1+\p_1}{p_1}.$ Using Theorem \ref{th:lw}
we have \begin{align}
    \int_0^{|t|} & \|\om(t,s)\|_{\Ls\rightarrow L^{p_2}}\|Dh|_{a(s)}\|_{\Ls}\|\po\|_{L^{p_1}}\|A_1\|_{L^{p'_1}}ds \nonumber \\
    &\leq\int_0^{|t|} \frac{C(p_2) C_2 \|\psi_0\|_{L^{p_1}}}{|t-s|^{N\ipt}}\cdot\frac{(\|u_1\|_{Y}^{\p_1}+\|u_2\|_{Y}^{\p_1})\|u_1-u_2\|_{Y}}{(1+|s|)^{N\ipo(1+\p_1)}} ds \nonumber \\
&\leq\frac{C(p_2)C_1
C_2}{(1+|t|)^{N\ipt}}(\|u_1\|_{Y}^{\p_1}+\|u_2\|_{Y}^{\p_1})\|u_1-u_2\|_{Y},
\nonumber \end{align} where we used
$$N\left(\frac{1}{2}-\frac{1}{p_1}\right)(\alpha_1+1)>1$$
because $p_1=\alpha_1+2$ and $\alpha_1 $ satisfies \eqref{alpha12}, and the uniform estimates
 \begin{equation}\label{est:udh} \|Dh|_{a(s)}\|_{L^2_\sigma}\leq C_2,\qquad \textrm{for all }s\in\R ,\end{equation} which follow from $h$ being $C^1$ on $a\in\C,\ |a|\leq\delta ,$ with values in $L^2_\sigma ,$ see Proposition \ref{pr:cm} and $|a(s)|\leq\varepsilon_0\leq\delta$ for all $s\in\R ,$ see \eqref{abound}.

Now $$\int_0^{|t|} \|\om(t,s)\|_{\Ls\rightarrow
L^{p_2}}\|Dh\|_{\Ls}\|\po\|_{L^{p_2}}\|A_2\|_{L^{p'_2}}ds$$ is
estimated as \eqref{a2lp2}, and $$\int_0^{|t|}
\|\om(t,s)\|_{\Ls\rightarrow
L^{p_2}}\|Dh\|_{\Ls}\|\po\|_{L^{p_2}}\|A_{3,4}\|_{L^{p'_2}}ds$$ is
estimated as \eqref{a34lp2}.
\bigskip

\textbf{$\bullet\ L^{p_1}$ Estimate :}

 \begin{eqnarray}
    \|Nu_1-Nu_2\|_{L^{p_1}}(t)&\le &\|\int_0^t \Omega(t,s)
[F_2(\pe(s),u_1(s))-F_2(\pe(s),u_2(s))]ds\|_{L^{p_1}}\nonumber\\
 &+&C\int_0^{|t|}\|\Omega(t,s)\|_{\Ls\mapsto
L^{p_1}}\|Dh|_{a(s)}\|_{\Ls}(\underbrace{|\<\po,A_1\>|}_{B_1}+\underbrace{|\<\po,A_2\>|}_{B_2}+\underbrace{|\<\po,A_{3,4}\>|}_{B_{3,4}})ds \nonumber
  \end{eqnarray}
The second integral is estimated as in the previous $L^{p_2}$ estimates on $B_1$, $B_2$ and $B_{3,4}$ to obtain the required
bound. For the first integral moving the norm inside the integration
and applying $L^{p'_1}\mapsto L^{p_1}$ estimates for $\Omega(t,s)$
and \eqref{def:a13} for the nonlinear term would require the control
of $A_2$ and $A_{3,4}$ in $L^{p'_1}.$ The latter, unfortunately, can no longer be
interpolated between $L^2$ and $L^{p_2}.$ To avoid this difficulty
we separate and treat differently the part of the nonlinearity
having an $A_2$ and $A_{3,4}$ like behavior by decomposing $\R^N$ in two disjoints
measurable sets related to the inequality \eqref{def:a13}:
 \begin{eqnarray}
 V_1(s)&=&\{x\in\R^N\ |\
 |F_2(\pe(s,x),u_2(s,x))-F_2(\pe(s,x),u_1(s,x))|\le C
 (A_2(s,x)+A_{3,4}(s,x))\},\nonumber\\
 V_2(s)&=&\R^N\setminus V_1(s)\nonumber\end{eqnarray}
On $V_2(s),$ using polar representation of complex numbers, we
further split the nonlinear term into:
 \begin{eqnarray}
 \lefteqn{F_2(\pe(s,x),u_1(s,x))-F_2(\pe(s,x),u_2(s,x))=e^{i\theta(s,x)}C(A_2(s,x)+A_{3,4}(s,x))}\nonumber\\
 &&+\underbrace{e^{i\theta(s,x)}[|F_2(\pe(s,x),u_1(s,x))-F_2(\pe(s,x),u_2(s,x))|-C(A_2(s,x)+A_{3,4}(s,x))]}_{G(s,x)}\nonumber\end{eqnarray}
where, due to inequality \eqref{def:a13}, $|G(s,x)|\le
CA_1(s,x)$ on $V_2(s).$ Then we have:
 \begin{eqnarray}\lefteqn{\int_0^t \Omega(t,s)
[F_2(\pe(s),u_1(s))-F_2(\pe(s),u_2(s))]ds=\int_0^t \Omega(t,s)
 (1-\chi (s))G(s)ds}\nonumber\\
 &&+\underbrace{\int_0^t \Omega(t,s)
[\chi(s)(F_2(\pe(s),u_1(s))-F_2(\pe(s),u_2(s)))+(1-\chi(s))e^{i\theta(s)}C(A_2(s)+A_{3,4}(s))]ds}_{I(t)},\nonumber\end{eqnarray}
where $\chi(s)$ is the characteristic function of $V_1(s).$ Now
 $$\|\int_0^t \Omega(t,s)
 (1-\chi (s))G(s)ds\|_{L^{p_1}}\le \int_0^{|t|}\| \Omega(t,s)\|_{L^{p'_1}\mapsto
 L^{p_1}}C\|A_1(s)\|_{L^{p'_1}}ds$$
 and estimates as in the previous step for $A_1$ give the required
 decay, see \eqref{a1p1} and the inequalities following it. For $I(t)$ we use interpolation:
        $$\|I(t)\|_{L^{p_1}}\leq\|I(t)\|_{L^2}^{1-\theta}\|I(t)\|_{L^{p_2}}^{\theta}\leq\|I(t)\|_{L^2}^{1-\theta}\left(\int_0^{|t|} \|\om(t,s)\|_{L^{p'_2}\mapsto L^{p_2}}\|A_2+A_{3,4}\|_{L^{p'_2}} ds\right)^{\theta}$$
where $\frac{1}{p_1}=\frac{1-\theta}{2}+\frac{\theta}{p_2}$. We know
from previous step that the above integral decays as
$(1+|t|)^{-N\ipt}$ and below we will show its $L^2$ norm will be
bounded. Therefore, since
 $\theta N\left(\frac{1}{2}-\frac{1}{p_2}\right)=N\left(\frac{1}{2}-\frac{1}{p_1}\right)$ we have: $$\sup_t
(1+|t|)^{N\ipo}\|I(t)\|_{L^{p_1}}<\infty$$ and the $L^{p_1}$
estimates are complete.
\bigskip

\textbf{$\bullet\ L^2$ Estimate :}
\bigskip

We will use $L^2\rightarrow L^2$ bounds for $\om(t,s),$ see Theorem
\ref{th:lp}, to control the $B_{1-4}$ terms. For the $A_{1-4}$ terms
we avoid $L^2\rightarrow L^2$ bounds because that would require us
to control the $L^{2(\alpha_2+1)}$ norm of functions in $Y_i$ which
is impossible since it can no longer be interpolated between the
norms in $L^2$ and $L^{p_2},\ p_2=\alpha_2+2.$ Instead we use the
decomposition:
    $$\om(t,s)=\ehs P_c+\tilde{T}(t,s)+(T(t,s)-\tilde{T}(t,s))$$
where for $t\geq s$
 \begin{align}
 \tilde{T}(t,s)&=\int_s^{\rm{min} \{t,s+1 \}} e^{-iH(t-\tau)}P_c g_u(\tau)e^{-iH(\tau-s)}P_c d\tau \nonumber \\
 &=e^{-iH(t-s)}P_c\int_s^{\rm{min}\{t,s+1\}}e^{-iH(t-s)}P_c e^{iH(\tau-s)}P_c
 g_u(\tau)e^{-iH(\tau-s)}P_c d\tau \nonumber \end{align}
while for $t<s$
 \begin{align}
 \tilde{T}(t,s)&=\int_s^{\max\{t,s-1 \}} e^{-iH(t-\tau)}P_c g_u(\tau)e^{-iH(\tau-s)}P_c d\tau \nonumber \\
 &=e^{-iH(t-s)}P_c\int_s^{\max\{t,s-1\}} e^{iH(\tau-s)}P_c
 g_u(\tau)e^{-iH(\tau-s)}P_c d\tau \nonumber \end{align}

 For $\ehs
P_c$ and $\tilde{T}(t,s)$ we will use  Stricharz estimates, while
for $T(t,s)-\tilde{T}(t,s)$ we will use $L^{p'}\rightarrow L^2$
estimates, see Theorem \ref{th:lp}. All in all we have:
    \begin{align}
    \|Nu_1-Nu_2\|_{L^2}&\leq \int_0^{|t|} \|\om(t,s)\|_{L^2\rightarrow L^2}\|Dh\|_{L^2}(|B_1|+|B_2|+|B_{3,4}|)ds
    \nonumber\\
    &+\int_0^{|t|} \|T(t,s)-\tilde{T}(t,s)\|_{L^{p'_1}\rightarrow L^2}C\|A_1\|_{L^{p'_1}}ds\nonumber \\
    & +\int_0^{|t|} \|T(t,s)-\tilde{T}(t,s)\|_{L^{p'_2}\rightarrow L^2}C(\|A_2\|_{L^{p'_2}}+\|A_{3,4}\|_{L^{p'_2}})ds \nonumber \\
    &+\|\int_0^t  \ehs P_c A_1(s) ds\|_{L^2}+\|\int_0^t \ehs P_c (A_2(s)+A_{3,4}(s)) ds\|_{L^2}\nonumber \\
    &+\|\int_0^t \tilde{T}(t,s) A_1 ds\|_{L^{2}}+\|\int_0^t \tilde{T}(t,s)(A_2+A_{3,4})ds\|_{L^{2}} \nonumber
    \end{align}
First three integrals are estimated similar to the previous cases. We deduce that this integrals are
uniformly bounded by:
$$\tilde{C_i}(\|u_1\|_{Y_i}^{\p_1}+\|u_2\|_{Y_i}^{\p_1}+\|u_1\|_{Y_i}^{\p_2}+\|u_2\|_{Y_i}^{\p_2})\|u_1 -u_2\|_{Y_i} $$

For the fourth integral we use Stricharz estimate:
   $$\sup_{t\in\R}\|\int_0^t \ehs P_c A_1 ds\|_{L^2}\leq C_s\Big(\int_\R \|A_1(s)\|_{L^{p_1'}}^{\gamma_1'}ds\Big)^{\frac{1}{\gamma_1'}}$$
where $\frac{1}{\gamma_1'}+\frac{1}{\gamma_1}=1,$ and
$\frac{2}{\gamma_1}=N\left(\frac{1}{2}-\frac{1}{p_1}\right).$
Furthermore we have
\begin{align}
\| A_1\|_{L^{\gamma_1'}L^{{p'_1}}}&\leq C_{13}\Big[\int_\R \frac{ds}{(1+|s|)^{N(1+\p_1)\gamma_1'\ipo}}\Big]^{\frac{1}{\gamma_1'}}(\|u_1\|_{Y_i}^{\p_1}+\|u_2\|_{Y_i}^{\p_1})\|u_1-u_2\|_{Y_i} \nonumber \\
&\leq C_{13}C_{10}
(\|u_1\|_{Y_i}^{\p_1}+\|u_2\|_{Y_i}^{\p_1})\|u_1-u_2\|_{Y_i}
\label{A1L2} \end{align} where $C_{10}=\int_\R
\frac{ds}{(1+|s|)^{N(1+\p_1)\gamma_1'\ipo}}ds<\infty$ since
$N(1+\p_1)\gamma_1'\ipo>1$ .

Similarly, for the fifth integral:
$$\sup_{t\in\R}\|\int_0^t \ehs P_c (A_2+A_{3,4})ds\|_{L^2}\leq C_s\left[\left(\int_\R \|A_2(s)\|_{L^{p'_2}}^{\gamma_2'}ds\right)^{\frac{1}{\gamma_2'}}+\left(\int_\R \|A_{3,4}(s)\|_{L^{p'_2}}^{\gamma_2'}ds\right)^{\frac{1}{\gamma_2'}}\right]$$
where $\frac{1}{\gamma_2'}+\frac{1}{\gamma_2}=1,$ and
$\frac{2}{\gamma_2}=N\left(\frac{1}{2}-\frac{1}{p_2}\right).$ Using
again the estimates we obtained before for $A_2$ and $A_{3,4}$ we
get:
\begin{align}
\| A_{3,4}\|_{L^{\gamma_2'}_s L^{{p'_2}}}&\leq C_{11}\Big[\int_\R \frac{(\log(2+|s|))^{2m_i\gamma_2'}ds}{(1+|s|)^{2n_i\gamma_2'}}\Big]^{\frac{1}{\gamma_2'}}(\|u_1\|_{Y_i}+\|u_2\|_{Y_i})\|u_1-u_2\|_{Y_i} \nonumber \\
&\leq  C_{11} C_8 (\|u_1\|_{Y_i}+\|u_2\|_{Y_i})\|u_1-u_2\|_{Y_i}
\label{A34L2} \end{align} where $C_8= \int_\R
\frac{(\log(2+|s|))^{2m_i\gamma_2'}ds}{(1+|s|)^{2n_i\gamma_2'}}<\infty$
since $2n_i\gamma_2'>1$ and:
\begin{align}
\| A_2\|_{L^{\gamma_2'}_s L^{{p'_2}}}&\leq C_{12}\Big[\int_\R \frac{(\log(2+|s|))^{(1+\p_2)m_i\gamma_2'}ds}{(1+|s|)^{(1+\p_2)n_i\gamma_2'}}\Big]^{\frac{1}{\gamma_2'}}(\|u_1\|_{Y_i}^{\p_2}+\|u_2\|_{Y_i}^{\p_2})\|u_1-u_2\|_{Y_i} \nonumber \\
&\leq  C_9
(\|u_1\|_{Y_i}^{\p_2}+\|u_2\|_{Y_i}^{\p_2})\|u_1-u_2\|_{Y_i}
\label{A2L2} \end{align} where $C_9=\int_\R
\frac{(\log(2+|s|))^{(1+\p_2)m_i\gamma_2'}ds}{(1+|s|)^{(1+\p_2)n_i\gamma_2'}}<\infty$
since $(1+\p_2)n_1\gamma_2'>1$.

Now for the last two integrals consider
$$\tilde{A}_i(t)=\|\int_0^t \tilde{T}(t,s) A_i(s) ds \|_{L^2},\qquad i=\overline{1,4}$$
Fix $t\geq 0$ and $i\in\overline{1,4}.$ The case $t<0$ is treated analogously. We have
\begin{align}
\tilde{A}_i(t)&\leq\sup_{\|\tilde{v}\|_{L^2}=1} \Big|\<\tilde{v}, \int_0^t \tilde{T}(t,s) A_i(s) ds\>\Big| \nonumber \\
&\leq\sup_{\|\tilde{v}\|_{L^2}=1}\int_0^t \Big|\<e^{iH(t-s)}P_c \tilde{v}, \int_s^{\rm{min}\{t,s+1\}}e^{iH(\tau-s)}P_c g_u(\tau)e^{-iH(\tau-s)}P_c A_i(s)d\tau\>\Big|ds \nonumber \\
&\leq\sup_{\|\tilde{v}\|_{L^2}=1}\int_0^t \|e^{iH(t-s)}P_c \tilde{v}\|_{L^p} \int_s^{\rm{min}\{t,s+1\}}\|e^{iH(\tau-s)}P_c g_u(\tau)e^{-iH(\tau-s)}P_c\|_{L^{p'}\rightarrow L^{p'}} d\tau \|A_i\|_{L^{p'}}ds\nonumber \\
&\leq\sup_{\|\tilde{v}\|_{L^2}=1}\int_0^t \|e^{iH(t-s)}P_c
\tilde{v}\|_{L^p}\ C\
\sup_{\tau\in[s,s+1]}\|\widehat{g_u(\tau)}\|_{L^1}
\|A_i\|_{L^{p'}}ds \nonumber \end{align} where we used the Fourier multiplier
type estimate $\|e^{+iH(\tau-s)}F e^{-iH(\tau-s)}\|_{L^p\rightarrow
L^p}\leq C\|\hat{F}\|_{L^1}$ for all $1\leq p\leq\infty$ and $|\tau-s|\leq 1;$ see
Appendix in \cite[Theorem 5.2]{km:asd3}. Note
that by Stricharz estimates there exist a fixed constant $C>0$ such that for all $\tilde{v}\in L^2:$
$$\|e^{iH(t-s)}\tilde{v}\|_{L^{\gamma_i} L^{p_i}_x}\leq C
\|\tilde{v}\|_{L^2}$$  and using \eqref{A1L2},
\eqref{A34L2} and \eqref{A2L2} for $\|A_i\|_{L^{\gamma'_i}
L^{p'_i}}$ we get that $\tilde{A}_i(t)$ are bounded uniformly for $t\in\R .$

%
%    $$\sup_{t\in\R}\|\int_0^t \ehs P_c A_2 ds\|_{L^2}\leq C_s\Big(\int_\R \|A_2\|_{L^{p_2'}}^{\gamma_2'}ds\Big)^{\frac{1}{\gamma_2'}}$$
%where $\frac{1}{\gamma_2'}+\frac{1}{\gamma_2}=1,$ and
%$\frac{2}{\gamma_2}=N\left(\frac{1}{2}-\frac{1}{p_2}\right).$
%Furthermore we have
%\begin{align}
%\| A_2\|_{L^{\gamma_2'}_sL^{{p'_2}}}&\leq C_{13}\Big[\int_\R \frac{ds}{(1+|s|)^{(1+\p_2)b_i\gamma_2'}}\Big]^{\frac{1}{\gamma_2'}}(\|u_1\|_{Y_i}^{\p_2}+\|u_2\|_{Y_i}^{\p_2})\|u_1-u_2\|_{Y_i} \nonumber \\
%&\leq C_{13}C_{10}
%(\|u_1\|_{Y_i}^{\p_2}+\|u_2\|_{Y_i}^{\p_2})\|u_1-u_2\|_{Y_i}
%\nonumber \end{align} where $C_{10}=\int_\R
%\frac{ds}{(1+|s|)^{(1+\p_2)b_i\gamma_2'}}ds<\infty$ since
%$(1+\p_1)b_i\gamma_2'>1$.

The $L^2$ estimates are now complete and the proofs of Lemma
\ref{lm} and Theorem \ref{mt} are finished. $\Box$

To obtain Corollary \ref{cor} in the case (i), i.e.
 $$\alpha_1\ge\frac{4}{N},\quad {\rm or}\quad p_2=\alpha_2+2<\frac{2N}{2+N-N\p_1}$$
we use Riesz-Thorin interpolation for $2<p<p_2$ while for $p>p_2$ we return to the integral form
of the equation for $r(t):$
$$
r(t)=\om(t,0)r(0)-\int_0^t
\om(t,s)P_c[iF_2(\pe,r)-Dh|_{a(s)}\<\po,iF_2(\pe,r)\>]ds
=\om(t,0)r(0)-N(r),
$$
see \eqref{sor}, and use the same arguments as in $L^{p_2}$ estimates for Lemma \ref{lm}
with $u_1=r,$ $u_2=0$ and $p_2$ replaced by an arbitrary $p,$
$p_2<p<2N/(N-2).$ Same works for $p_2<p<2N/(N-2)$ in the cases (ii) and (iii), i.e.
 $$\alpha_1<\frac{4}{N}\quad {\rm and}\quad p_2=\alpha_2+2=\frac{2N}{2+N-N\p_1},$$
respectively
 $$\alpha_1<\frac{4}{N}\quad {\rm and}\quad p_2=\alpha_2+2>\frac{2N}{2+N-N\p_1}.$$
In these last two cases, for $\alpha_1+2=p_1<p<p_2$ we use the arguments for $L^{p_1}$ estimates in Lemma {\ref{lm} with $p_1 $ replaced by $p.$ Riesz-Thorin interpolation gives the required estimates for $2<p<p_1$ and finishes the proof of Corollary \ref{cor}.

For the proof of Corollary \ref{cor:as} we return to the dynamics on the center manifold given by equation \eqref{a}:
 $$i\frac{da}{dt}=E(|a(t)|)a(t)+\<\psi_0, g(\psi_{E(t)}+r(t))-g(\psi_{E(t)})\>.$$
Changing to $\tilde a(t)=e^{i\int_0^tE(|a(s)|)ds}a(t)$ which eliminates the fast phase oscillations of the complex valued function $a(t),$ and using the symmetries of $g(z)$ with respect to rotations of the complex plane \eqref{gsym} we get:
 \begin{equation}\label{eq:tildea}i\frac{d\tilde a}{dt}=\<\psi_0, g(\tilde\psi_{E(t)}+\tilde r(t))-g(\tilde\psi_{E(t)})\>=\<\psi_0, Dg|_{\tilde\psi_{E(t)}}[\tilde r(t)]+F_2(\tilde\psi_{E(t)},\tilde r(t))\>,\end{equation}
where $\tilde\psi_{E(t)}=e^{i\int_0^tE(|a(s)|)ds}\psi_{E(t)}=\tilde a(t)\psi_0+h(\tilde a(t))$ and $\tilde r(t)=e^{i\int_0^tE(|a(s)|)ds}r(t).$ Now the right hand side of \eqref{eq:tildea} is integrable in time on $t\in [-\infty,\infty].$ Indeed, for the nonlinear term we use \eqref{def:a13} with $u_1=\tilde r$ and $u_2\equiv 0$ to get:
 \begin{eqnarray}|\<\psi_0,F_2(\tilde\psi_{E(t)},\tilde r(t))\>|&\leq&C\left[\|\psi_0\|_{L^{p'_1}}\|\tilde r\|_{L^{p_1}}^{\alpha_1+1}+\|\psi_0\|_{L^{p'_2}}\left(\|\tilde r\|_{L^{p_2}}^{\alpha_2+1}+\|\ |\tilde\psi_E|^{\alpha_1-1}+|\tilde\psi_E|^{\alpha_2-1}\|_{L^q}\|\tilde r\|^2_{L^{p_2}}\right)\right]\nonumber\\
 &\leq&\frac{\tilde C}{(1+|t|)^{\beta}} \label{est:anlin}\end{eqnarray}
where $$\beta=\min\left\{(\alpha_1+1) N\left(\frac{1}{2}-\frac{1}{p_1}\right),(\alpha_2+1) N\left(\frac{1}{2}-\frac{1}{p_2}\right),2N\left(\frac{1}{2}-\frac{1}{p_2}\right)\right\}>1,$$ see the estimates for $|r(t)|=|\tilde r(t)|$ in Theorem \ref{mt}.

For the linear term we are forced to use  $L^2_{-\sigma}$ estimates.
By \eqref{eq:dg1}, \eqref{eq:dg}, and \eqref{gest} combined with Cauchy-Schwarz and H\" older inequalities:
 \begin{eqnarray}\lefteqn{|\<\psi_0, Dg|_{\tilde\psi_{E(t)}}[\tilde r(t)]\>| \leq  }\nonumber\\ &&\|\psi_0\|_{L^2}\left(\|<x>^\sigma |\tilde\psi_{E(t)}|^{\alpha_1}\|_{L^\infty } +\|<x>^\sigma |\tilde\psi_{E(t)}|^{\alpha_2}\|_{L^\infty}\right)\|<x>^{-\sigma}\tilde r(t)\|_{L^2}\leq C \|r(t)\|_{L^2_{-\sigma}},\nonumber\end{eqnarray}
where the uniform in time bounds for the norms involving
$|\tilde\psi_{E(t)}|=|\psi_{E(t)}|$ follow from \eqref{eq:lee0}. For
the $L^2_{-\sigma}$ estimate of  $r(t)$ we turn to \eqref{sor} which
combined with \eqref{def:a13} where $u_1=r$ and $u_2\equiv 0$ gives:
 \begin{eqnarray} \|r(t)\|_{\Lsn}&\leq & \|\Omega(t,0)\|_{L^{q_0}\mapsto\Lsn}\|r(0)\|_{L^{q_0}}\label{est:rw}\\
 &&+\int_0^{|t|}\|\Omega(t,s)\|_{L^{q_1}\mapsto\Lsn}\|r(s)\|^{\alpha_1+1}_{L^{(\alpha_1+1)q_1}}ds
 +\int_0^{|t|}\|\Omega(t,s)\|_{L^{q_2}\mapsto\Lsn}\|r(s)\|^{\alpha_2+1}_{L^{(\alpha_2+1)q_2}}ds\nonumber\\
 &&+\int_0^{|t|}\|\Omega(t,s)\|_{\Ls\mapsto\Lsn}\left(\|<x>^{\sigma} |\psi_{E(s)}|^{\alpha_1-1}+<x>^{\sigma}|\psi_{E(s)}|^{\alpha_2-1}\|_{L^q}\right)\|r(s)\|^2_{L^{p_2}}ds\nonumber\\
 &&+\int_0^{|t|}\|\Omega(t,s)\|_{\Ls\mapsto\Lsn}\|Dh|_{a(s)}\|_{\Ls} |\<\psi_0,F_2(\psi_{E(s)},r(s))\>|ds\nonumber\end{eqnarray}
Now, from the weighted estimates for $\Omega(t,s)$ in Theorem
\ref{th:lw}, uniform bounds for $\psi_{E(s)}$ and $Dh|_{a(s)}$ see
\eqref{eq:lee0} respectively \eqref{est:udh}, and from
\eqref{est:anlin} we get that the last two integrals decay like
$(1+|t|)^{-\tilde\beta},$ where $ \tilde\beta=\min\{N/2,\beta\}>1.$
In the remaining terms, to avoid singularities at $t=0$ we use
$q_0=2,\ q_1=p'_1,\ q_2=p'_2$ for $|t|\leq 1,$ while for $|t|>1$ we
fix
 \begin{equation}\label{def:p3}\frac{2N}{N-2}<p_3<\frac{2N}{N+2-N\alpha_1}\end{equation}
such that $p_3\leq p$ with $p$ given in the hypotheses of Corollary
\ref{cor:as}. Then we use $q_0=p'_3$ and split the first two
integrals on the right hand side of \eqref{est:rw} into integrals
from zero to $|t|-1$ where we choose $q_1=q_2=p'_3$ and integrals
from $|t|-1$ to  $|t|$ where we choose $q_1=p'_1,\ q_2=p'_2.$ This
way the first three terms on the right hand side of \eqref{est:rw}
decay like $(1+|t|)^{-\tilde\beta}$ where
$$\tilde\beta=\min\left\{N\left(\frac{1}{2}-\frac{1}{p_3}\right),(\alpha_1+1)N\left(\frac{1}{2}-\frac{1}{(\alpha_1+1)p'_3}\right)\right\}>1.$$

\begin{remark}\label{rmk:whysc} The restriction $\alpha_1>4/N$ is necessary for the existence of a $p_3$ with the properties
\eqref{def:p3} which in turn insures the integrability in time of
the first two integrals in \eqref{est:rw}. Such a restriction was
not needed for the integrability of \eqref{est:anlin}, the nonlinear
terms on the right hand side of \eqref{eq:tildea}. Consequently a
decomposition that removes the linear term on the right hand side of
\eqref{eq:tildea} will also remove this restriction to asymptotic
stability.\end{remark}

All in all we now have
 $$\frac{d\tilde a}{dt}=-i\<\psi_0, Dg|_{\tilde\psi_{E(t)}}[\tilde r(t)]+F_2(\tilde\psi_{E(t)},\tilde r(t))\>,\quad |\<\psi_0, Dg|_{\tilde\psi_{E(t)}}[\tilde r(t)]+F_2(\tilde\psi_{E(t)},\tilde r(t))\>|\leq \frac{C}{(1+|t|)^{\beta}},\ \beta>1.$$
Consequently there exist $\tilde a_{\pm\infty}\in\C$ such that
 $$\lim_{t\rightarrow\pm\infty}\tilde a(t)=\tilde a_{\pm\infty},\quad |\tilde a(t)-\tilde a_{\pm\infty}|\leq \frac{C_1}{(1+|t|)^{\beta-1}}.$$
Moreover, because $E(|a|),\ |a|\leq\varepsilon_0$ is a $C^1$ function hence Lipschitz, we get $E(t)=E(|a(t)|)=E(|\tilde a(t)|)\rightarrow E(|a_{\pm\infty}|)$ as $t\rightarrow\pm\infty$ and the function
 $$\theta(t)=\left\{\begin{array}{ll}  \frac{1}{t}\int_0^t E(|a(s)|)-E(|a_{+\infty}|)ds & {\rm if}\ t\geq 1\\
                                       \frac{1}{t}\int_0^t E(|a(s)|)-E(|a_{-\infty}|)ds & {\rm if}\ t\leq -1
                      \end{array}\right.$$
converges to zero as $t\rightarrow\pm\infty.$ Finally, from $\psi_{E(t)}=a(t)\psi_0+h(a(t)),$ by the continuity of $h(a)$ and its equivariance with respect to rotations see Proposition \ref{pr:cm} we get:
 $$\lim_{t\rightarrow\pm\infty}\|\psi_{E(t)}-e^{-it(E(|a_{\pm\infty}|)+\theta(t))}\psi_{E(|a_{\pm\infty}|)}\|_{H^2\cap L^2_\sigma}=0.$$
and the proof of Corollary \ref{cor:as} is finished.

In the next section we obtain the estimates for the propagator $\Omega(t,s),\ t,s\in\R$ of \eqref{z}. Note that they were essential in proving Theorem \ref{mt} and Corollaries \ref{cor} and \ref{cor:as}.

\section{Linear Estimates}\label{se:lin}
Consider the linear Schr\" odinger equation with a potential in four
and five space dimensions:
\begin{align}
i\frac{\partial u}{\partial t}&=\underbrace{(-\Delta+V(x))}_{H}u \nonumber \\
u(0)&=u_0. \nonumber \end{align} If V satisfies hypothesis (H1) (i)
1. and (ii) it is known, see \cite[Example 7.8]{mm:ae}, that for
$N=4,5,$ and $\sigma>N/2,$ there exists a constant $C_N>0$ such that
\begin{equation}
    \|e^{-iHt}P_c u_0\|_{\Lsn}\leq \frac{C_N}{(1+|t|)^{\frac{N}{2}}}\|u_0\|_{\Ls}
    \label{eq:ls}
\end{equation}
where $P_c$ is the projection onto the continuous spectrum of
$H=-\Delta+V.$

In addition, if V satisfies (H1) (i) 1., 3. and (ii) then for each
$2\leq p\leq\infty,\ 1/p'+1/p=1$ there exists a constant $C_p>0$
such that:
\begin{equation}
    \|e^{-iHt}P_c u_0\|_{L^p}\leq \frac{C_p}{|t|^{N(\frac{1}{2}-\frac{1}{p})}}\|u_0\|_{L^{p'}},\quad N=4, 5
    \label{eq:lp}
\end{equation}
see for example \cite{kn:jss}.

We would like to extend these estimates to the linearized dynamics
around the time dependent motion on center manifold. We consider the
linear equation with initial data at time $s$ in the range of $P_c:$
\begin{align}
i\frac{dz}{dt}&=Hz(t)+P_cDg|_{\psi_{E(t)}} [z(t)]+iDh|_{a(t)}i\<\psi_0,Dg|_{\psi_{E(t)}} [z(t)]\> \nonumber \\
z(s)&=v\in{\rm range}P_c \nonumber
\end{align}
where $\Dg [z]=\frac{d}{d\varepsilon}g(\pe+\varepsilon z)|_{\varepsilon=0}=\frac{\partial}{\partial u}g(u)|_{u=\pe}z+\frac{\partial}{\partial \bar{u}}g(u)|_{u=\pe}\overline{z}$.

By Duhamel's principle we have: \begin{equation} z(t)=\ehs
P_cv-\int_s^t \eht P_c\left(iDg|_{\psi_{E(\tau)}}
[z(\tau)]-Dh|_{a(\tau)}\<\psi_0,iDg|_{\psi_{E(\tau)}}
[z(\tau)]\>\right) d\tau \label{soz} \end{equation}

In the next theorems we will extend estimates of type (\ref{eq:ls})-(\ref{eq:lp}) to the operators
 $$\om(t,s)v=z(t),\quad {\rm and}\quad T(t,s)=\om(t,s)-\ehs P_c,$$
relying on the fact that $\psi_{E(t)}$ is small and localized in
space, see \eqref{eq:lee0}. The arguments can be extended for large
$\psi_{E(t)}$ provided for a certain fixed solution $\psi_{E}$ of
\eqref{eq:ev} we have
$\inf_{\theta\in\R}\|\psi_{E(t)}-e^{i\theta}\psi_{E}\|_{H^1}$ is
small uniformly in $t\in\R ,$ see \cite{kz:asLs}. We start with
weighted estimates. While the approach is similar to the one in
\cite{kz:as2d}, see also \cite{km:asd3}, we include the proofs for
completeness.
\begin{theorem}\label{th:lw} Fix $\sigma > N/2,$ and $\frac{2N}{N-2}<q_2<\frac{2N}{N-4}$. There exists
$\varepsilon_1(q_2)>0$ such that if
$\|<x>^{4\sigma}\psi_{E(t)}\|_{L^\infty}<\varepsilon_1$ for all $t\in\R ,$ then there are
constants $C_\sigma,$ $C_p,$ $C$ and $C(q_2)>0$ with the property that for any $t,s\in\R$
the following estimates hold:
\begin{align}
(i)\quad &\|\om(t,s)\|_{\Ls\mapsto\Lsn}\leq\frac{C_\sigma}{(1+|t-s|)^{\frac{N}{2}}}; \nonumber \\
(ii)\quad &\|\om(t,s)\|_{\Ls\mapsto L^p}\leq\frac{C_p}{|t-s|^{N\ip}},\quad\textrm{for all}\ 2\leq p< \frac{2N}{N-2};\nonumber \\
(iii)\quad &T(t,s)\in L^2(\R,L^2\mapsto\Lsn)\cap
L^\infty(\R,L^2\mapsto\Lsn); \nonumber \\
(iv)\quad &\|\om(t,s)\|_{L^{p'}\mapsto\Lsn}\leq\frac{C_p}{|t-s|^{N(\frac{1}{2}-\frac{1}{p})}}\quad\textrm{for all $2\leq p\leq q_2$} ;\nonumber \\
&\|T(t,s)\|_{L^{p'}\mapsto\Lsn}\leq\left\{ \begin{array}{ll}
C & \textrm{for $|t-s|\leq 1$ and $2\leq p\leq \frac{2N}{N-2}$}, \\
\frac{C_p}{|t-s|^{N(\frac{1}{2}-\frac{1}{p})-1}} & \textrm{for $|t-s|\leq 1$ and $\frac{2N}{N-2}<p\leq q_2$}, \\
\frac{C_p}{(1+|t-s|)^{N(\frac{1}{2}-\frac{1}{p})}} & \textrm{for
$|t-s|>1$ and $2\leq p\leq q_2$}. \end{array} \right. \nonumber
 \end{align}

\end{theorem}

\textbf{Proof of Theorem \ref{th:lw}} Fix $s\in\R$ and let $q_1=\frac{2N}{N-2}.$

$(i)$ By definition, we have $\om(t,s)v=z(t)$ where $z(t)$ satisfies
equation \eqref{soz}. We are going to prove the estimate by showing
that the nonlinear equation \eqref{soz} can be solved via
contraction principle argument in an appropriate functional space.
To this extent let us consider the functional space
$$X_1:=\{u\in C(\R,L_{-\sigma}^2(\R^N))\ :\ \sup_{t\in\R}(1+|t-s|)^{\frac{N}{2}}\|u(t)\|_{L_{-\sigma}^2}<\infty\}$$
endowed with the norm
$$\|u\|_{X_1}:=\sup_{t\in\R}\{(1+|t-s|)^{\frac{N}{2}}\|u(t)\|_{L_{-\sigma}^2}\}<\infty$$
Note that the inhomogeneous term in (\ref{soz}) $z_0=\ehs P_c v$
satisfies $z_0\in X_1$ and
\begin{equation}
    \|z_0\|_{X_1}\leq C_N\|v\|_{\Ls}
    \label{eq:z0}
\end{equation}
because of (\ref{eq:ls}). We collect the $z$ dependent part of the right hand side of (\ref{soz}) in a linear operator $L(s):X_1\rightarrow X_1$
\begin{equation}
    [L(s)z](t)=-\int_s^t \eht P_c\left(iDg|_{\psi_{E(\tau)}}[z(\tau)]-Dh|_{a(\tau)}i\<\psi_0,Dg|_{\psi_{E(\tau)}}[z(\tau)]\>\right) d\tau  \label{L} \end{equation}
We will show that $L$ is a well defined bounded operator from $X_1$
to $X_1$ whose operator norm can be made less or equal to 1/2 by
choosing $\varepsilon_1$ sufficiently small. Consequently $Id-L$ is
invertible and the solution of the equation (\ref{soz}) can be
written as $z=(Id-L)^{-1}z_0$. In particular
$$\|z\|_{X_1}\leq(1-\|L\|)^{-1}\|z_0\|_{X_1}\leq 2\|z_0\|_{X_1}$$
which in combination with the definition of $\om$, the definition of
the norm $X_1$ and the estimate (\ref{eq:z0}), finishes the proof of
$(i)$.

By computing the $\Lsn $ norm of both the left hand side and right
hand side of \eqref{L}, for $t>s$ we have:
\begin{align}
\|[L(s)z](t)\|_{L_{-\sigma}^2}\leq & \int_s^t \|\eht P_c\|_{L_{\sigma}^2\rightarrow L_{-\sigma}^2}\Big[\|\Dg[z]\|_{\Ls}+\|Dh|_{a(\tau)}\|_{\C\rightarrow\Ls}|\langle\psi_0, \Dg[z]\rangle|\Big]d\tau \nonumber \\
&\leq \int_s^t \|\eht P_c\|_{L_{\sigma}^2\rightarrow L_{-\sigma}^2}[\|\Dg[z]\|_{\Ls} \nonumber \\
&\quad\quad\quad\quad\quad+\|Dh|_{a(\tau)}\|_{\C\rightarrow\Ls}\|\psi_0\|_{L^2}\|\Dg[z]\|_{L^2}]d\tau
\nonumber \end{align} On the other hand, from \eqref{eq:dg},
\eqref{eq:dg1}, and \eqref{gest} we obtain:
\begin{eqnarray}
\|\Dg[z]\|_{\Ls}&\leq
&\|\xst(|\pe|^{\p_1}+|\pe|^{\p_2})\|_{L^{\infty}}\|z\|_{\Lsn}\leq\varepsilon_1^{\p_1}\|z\|_{\Lsn}\label{F1loc}\\
\|\Dg[z]\|_{L^2}&\leq &\|\xs
(|\pe|^{\p_1}+|\pe|^{\p_2})\|_{L^{\infty}}\|z\|_{\Lsn}\leq\varepsilon_1^{\p_1}\|z\|_{\Lsn}\label{F1loc1}
\end{eqnarray}
Also
$$\|Dh|_{a(\tau)}\|_{\C\rightarrow\Ls}\leq C_2,
\quad \textrm{for all}\ \tau\in\R$$ which follow from $h$ being
$C^1$ on $a\in\C,\ |a|\leq\delta ,$ with values in $L^2_\sigma ,$
see Proposition \ref{pr:cm} and
$|a(\tau)|\leq\varepsilon_0\leq\delta$ for all $\tau\in\R ,$ see
\eqref{abound}.

Using the last three relations, as well as the estimate
(\ref{eq:ls}) and the fact that $z\in X_1$ we obtain that
\begin{align}
\sup_{t>s}(1+|t-s|)^{\frac{N}{2}} \|[L(s)z](t)\|_{\Lsn}&\leq\varepsilon_1^{\p_1} \sup_{t>s}(1+|t-s|)^{\frac{N}{2}}\int_s^t \frac{C_N}{(1+|t-\tau|)^{\frac{N}{2}}} \cdot \frac{\|z\|_{X_1}}{(1+|\tau-s|)^{\frac{N}{2}}}{\ } d\tau \nonumber \\
&\leq\varepsilon_1^{\p_1}
\sup_{t>s}(1+|t-s|)^{\frac{N}{2}}\frac{\tilde
C\|z\|_{X_1}}{(1+|\frac{t-s}{2}|)^{\frac{N}{2}}}\leq
C\varepsilon_1^{\p_1}\|z\|_{X_1} \nonumber \end{align} Similar
arguments lead to $\sup_{t<s}(1+|t-s|)^{\frac{N}{2}}
\|[L(s)z](t)\|_{\Lsn}\leq C\varepsilon_1^{\p_1}\|z\|_{X_1},$ hence
$\|L(s)\|_{X_1\mapsto X_1}\leq C\varepsilon_1^{\p_1}.$ Now choosing
$\varepsilon_1$ small enough we get that $L(s)$ is a contraction
operator on the Banach space $X_1,$ therefore:
$$\|\Omega(t,s)\|_{\Ls\mapsto\Lsn}\leq\frac{C}{(1+|t-s|)^{\frac{N}{2}}}$$

$(ii)$ From part $(i)$ we already know that \eqref{soz} has a unique
solution in $\Lsn$ provided that $v\in\Ls$. We are going to show
that the right had  side of \eqref{soz} in $L^p$. Indeed, using
\eqref{eq:lp} and $L^2_\sigma\hookrightarrow L^{p'}$ we have:
\begin{equation} \|e^{-iH(t-s)}P_c
v\|_{L^p}\leq\frac{C_p}{|t-s|^{N\ip}}\|v\|_{\Ls}. \label{lpls}
\end{equation} The remaining term satisfies for $t>s:$
\begin{align}
\|[L(s)z](t)\|_{L^p}&\leq\int_s^t \|e^{-iH(t-\tau)}P_c\|_{L^{p'}\rightarrow L^p}\Big[\|\Dg[z]\|_{L^{p'}}+\|Dh\|_{\C\rightarrow L^{p'}}|\<\po, \Dg[z]\>|\Big]d\tau \nonumber \\
&\leq\int_s^t\frac{C}{|t-\tau|^{N\ip}}\Big[\|\xs\pe\|_{L^\beta}+\|Dh\|\|\po\|_{L^2}\|\xs\pe\|_{L^\infty}\Big]\|z(\tau)\|_{\Lsn}d\tau \nonumber\\
&\leq\int_s^t \frac{C}{|t-\tau|^{N\ip}}\frac{\|v\|_{\Ls}}{(1+|\tau-s|)^{\frac{N}{2}}}d\tau \nonumber\\
&\leq\frac{C\|v\|_{\Ls}}{|t-s|^{N\ip}},\quad\textrm{for all $2\leq
p<\frac{2N}{N-2}$} \label{lpls2} \end{align} and same estimate can
be obtained for $t<s.$ Pluginig \eqref{lpls} and \eqref{lpls2} into
\eqref{soz} we get:
$$\|z(t)\|_{L^p}\leq\frac{C}{|t-s|^{N\ip}}\|v\|_{\Ls},\qquad \textrm{for all }t\in\R$$ which by the
definition $\om(t,s)v=z(t)$ finishes the proof of part $(ii)$.

$(iii)$ Denote: \begin{equation}
    T(t,s)v=W(t)
    \label{w}
\end{equation}
then, by plugging in (\ref{soz}), $W(t)$ satisfies the following equation:
\begin{align}
W(t)=&-i\int_s^t \eht P_c\Dg[\ehts P_c v]d\tau \nonumber \\
&+\int_s^t \eht P_c Dh|_{a(\tau)}i\<\po,\Dg[\ehts P_c v]\>d\tau \nonumber \\
&+[L(s)W](t) \label{sow}
\end{align}
Consequently \begin{align} \xsn W(t)=&-\int_s^t \xsn\eht P_c\left( i\Dg[\ehts P_c v]-Dh|_{a(\tau)}i\<\po,\Dg[\ehts P_c v]\>\right)d\tau \nonumber \\
& -\int_s^t \xsn\eht P_c\left(i\Dg[W(\tau)]-Dh|_{a(\tau)}i\<\po,
\Dg[W(\tau)]\>\right)d\tau .\nonumber \end{align} Then, by
\eqref{eq:ls}:
\begin{align} \|\xsn& W(t)\|_{L_t^2 L_x^2}\leq\Big\|\int_s^t \frac{C}{(1+|t-\tau|)^{\frac{N}{2}}}(1+\|Dh\|_{\Ls}\|\po\|_{\Ls}) \|\xs \Dg[\ehts P_c v]\|_{L^2}d\tau\Big\|_{L^2_t} \nonumber \\
&+\Big\|\int_s^t \frac{C}{(1+|t-\tau|)^{\frac{N}{2}}}(1+\|Dh\|_{\Ls}\|\po\|_{\Ls})(\|\xst g_u\|_{L^\infty}+\|\xst g_{\bar{u}}\|_{L^\infty})\|\xsn W(\tau)\|_{L^2_x}d\tau \Big\|_{L^2_t} \nonumber \\
&\leq \varepsilon_1^{\alpha_1} C\|K\|_{L^1}\|v\|_{L^2}+\varepsilon_1^{\alpha_1} C\|K\|_{L^1}\|\xsn
W\|_{L^2_t L^2_x} \nonumber\end{align} where we used Young inequality $\|K*f\|_{L^2}\leq\|K\|_{L^1}\|f\|_{L^2},$
with $K(t)=(1+|t|)^{-\frac{N}{2}}\in L^1_t,$ while for the term $$\xs \Dg[e^{-iHt}
P_c v]=\xs g_u\xs\xsn e^{-iHt} P_c v+\xs g_{\bar{u}}\xs\xsn e^{iHt}
P_c v$$ we used $\|\xst g_u\|_{L^\infty}$ and $\|\xst
g_{\bar{u}}\|_{L^\infty}$ is uniformly bounded in $t$ by $\varepsilon_1^{\alpha_1}$ since
$|g_u|=|g_{\bar{u}}|\leq C(|\pe|^{\p_1}+|\pe|^{\p_2})$ and the Kato
smoothing estimate $\|\xsn e^{-iHt} P_c v\|_{L^2(\R,L^2_x)}\leq
C\|v\|_{L^2_x}$. Choosing $\varepsilon_1<1/(C\|K\|_{L^1})$ we get
$\|\xsn W\|_{L^2_t L^2_x}<\infty$. In other words $T(t,s)\in
L^2(\R, L^2\rightarrow\Lsn)$.

Similarly, using now \eqref{eq:lp} with $p=2$ and $u_0=v,$ we obtain
for $t>s:$
\begin{align} \|\xsn W(t)\|_{L_x^2}&\leq\int_s^t \frac{C}{(1+|t-\tau|)^{\frac{N}{2}}}(1+\|Dh\|_{\Ls}\|\po\|_{\Ls}) \|\xs \Dg[\ehts P_c v]\|_{L^2_x}d\tau \nonumber \\
&+\int_s^t \frac{C}{(1+|t-\tau|)^{\frac{N}{2}}}(1+\|Dh\|_{\Ls}\|\po\|_{\Ls})(\|\xst g_u\|_{L^\infty}+\|\xst g_{\bar{u}}\|_{L^\infty})\|\xsn W(\tau)\|_{L^2_x}d\tau \nonumber \\
&\leq \varepsilon_1^{\alpha_1} C\|v\|_{L^2}+\varepsilon_1^{\alpha_1}
C\sup_{\tau\in\R}\|\xsn W(\tau)\|_{L^2_x} \nonumber\end{align} Same
argument works for $t<s.$ Then passing to supremum over $t\in\R$ on
the left hand side we get the required estimate provided
$\varepsilon_1$ is small enough.

$(iv)$ By definition of $T(t,s)$ (\ref{w}) and the similarity
between $t>s$ and $t<s$ estimates it is sufficient to prove that the
solution of (\ref{sow}) satisfies
$$\|W(t)\|_{\Lsn}\leq\left\{
\begin{array}{ll}C\|v\|_{L^{q_1'}} & \textrm{for $s\leq t\leq s+1$} \\
\frac{C\|v\|_{L^{q_2'}}}{|t-s|^{N(\frac{1}{2}-\frac{1}{q_2})-1}} & \textrm{for $s\leq t\leq s+1$} \\
\frac{C\|v\|_{L^{q_2'}}}{(1+|t-s|)^{N(\frac{1}{2}-\frac{1}{q_2})}} &
\textrm{for $t>s+1$} \end{array} \right.$$ The estimates for $2\leq
p\leq q_2$ are then obtained by Riesz-Thorin interpolation.

Let us also observe that it suffices to obtain estimates only for
the forcing terms in \eqref{sow}:
 \begin{eqnarray}
 f(t)&=&-i\int_s^t \eht P_c\Dg[\ehts P_c v]d\tau\label{def:f}\\
 \tilde f(t)&=&\int_s^t \eht P_c Dh|_{a(\tau)}i\<\po,\Dg[\ehts P_c v]\>d\tau\label{def:tildef}
 \end{eqnarray} because then we will be able to do the
contraction principle in the functional space in which $f(t),\tilde
f(t)$ are, and thus obtain the same decay for $W$ as for $f(t)$ and
$\tilde f(t).$ This time we will consider the functional spaces (recall that $s\in\R$ is a fixed number)
 \begin{eqnarray}
 X_1&=&\Big\{u\in C([s-1,s+1],\Lsn(\R^N))\ :\
 \sup_{|t-s|\leq1} \|u(t)\|_{\Lsn}<\infty
 \Big\}\nonumber\\
 X_2&=&\Big\{u\in C(\R,\Lsn(\R^N))\ :\  \sup_{|t-s|>1}
(1+|t-s|)^{N(\frac{1}{2}-\frac{1}{q_2})}\|u(t)\|_{\Lsn}<\infty
\nonumber\\
 &&\qquad \sup_{|t-s|\leq1}
|t-s|^{N(\frac{1}{2}-\frac{1}{q_2})-1}\|u(t)\|_{\Lsn}<\infty
 \Big\}\nonumber\end{eqnarray} endowed with the norms
 \begin{eqnarray}
 \|u\|_{X_1}&=&\sup_{|t-s|\leq1} \|u(t)\|_{\Lsn}\nonumber\\
  \|u\|_{X_2}&=&\max\left\{
\sup_{|t-s|\leq 1}
|t-s|^{N(\frac{1}{2}-\frac{1}{q_2})-1}\|u(t)\|_{\Lsn},\
\sup_{|t-s|>1}
(1+|t-s|)^{N(\frac{1}{2}-\frac{1}{q_2})}\|u(t)\|_{\Lsn}
\right\}\nonumber\end{eqnarray}

First we will investigate the short time behavior of the forcing
terms. If $s< t\leq s+1:$
\begin{align} \|f(t)\|_{\Lsn}&\leq\|\xsn\int_s^t \eht P_c \Dg [\ehts P_c v] d\tau\|_{L^2} \nonumber \\
&\leq\|\xsn\|_{L^\beta}\int_s^t \|\ehs e^{iH(\tau-s)} P_c g_u \ehts P_c v\|_{L^{q_2}} d\tau \nonumber\\
&\quad+\|\xsn\|_{L^\beta}\int_s^{s+\frac{t-s}{4}} \|e^{-iH(t+s-2\tau)} \ehts P_c g_{\bar{u}} e^{iH(\tau-s)} P_c v\|_{L^{q_2}} d\tau \nonumber\\
&\quad+\int_{s+\frac{t-s}{4}}^t \|\eht P_c\|_{\Ls\rightarrow\Lsn}\|\xs g_{\bar{u}} \ehts P_c v\|_{L^2} d\tau \nonumber \\
&\leq \int_s^t\frac{C}{|t-s|^{N(\frac{1}{2}-\frac{1}{q_2})}}\sup_{\tau\in[s,t]} \|\widehat{g_u}(\tau)\|_{L^1}\|v\|_{L^{q_2'}}d\tau \nonumber \\
&\quad+\int_s^{s+\frac{t-s}{4}}\frac{C}{|t+s-2\tau|^{N(\frac{1}{2}-\frac{1}{q_2})}}\sup_{\tau\in[s,t]} \|\widehat{g_{\bar{u}}(\tau)}\|_{L^1}\|v\|_{L^{q_2'}}d\tau \nonumber \\
&\quad+\int_{s+\frac{t-s}{4}}^t \frac{C}{(1+|t-\tau)^{\frac{N}{2}}}\|\xs g_{\bar{u}}\|_{L^\beta} \|\ehts P_c v\|_{L^{q_2}} d\tau \nonumber \\
&\leq \frac{C\sup_{\tau\in[s,t]}
(\|\widehat{g_u}(\tau)\|_{L^1}+\|\widehat{g_{\bar{u}}(\tau)}\|_{L^1}+\|\xs
g_{\bar{u}}(\tau)\|_{L^\beta})\|v\|_{L^{q_2'}}}{|t-s|^{N(\frac{1}{2}-\frac{1}{q_2})-1}}\nonumber
\end{align}
where we used the Fourier multiplier type estimates:
\begin{equation}\label{est:fm}\|e^{-iH(\tau-s)}F(x) e^{iH(\tau-s)}\|_{L^p(\R^N)\rightarrow L^p(\R^N)}\leq
C\|\hat{F}\|_{L^1(\R^N)},\quad \textrm{for all }|\tau-s|\leq
1,\textrm{ and } 1\leq p\leq\infty,\end{equation} where $\hat F$ is
the Fourier transform of $F$ with respect to the variable $x\in\R^N$
and $C>0$ is a fixed constant see \cite[Theorem 5.2]{km:asd3}.
Similarly we obtain $\|f\|_{\Lsn}\leq C\|v\|_{L^{q_1'}}$ for
$q_1=\frac{2N}{N-2}$ and $|t-s|\leq1.$ For the second term we have:
\begin{align}
\|\tilde{f}(t)\|_{\Lsn}&\leq\int_s^t \|\eht P_c v\|_{\Ls\rightarrow\Lsn}\|Dh|_a(\tau)\|_{\Ls}\< \po,\Dg[\pe,\ehts P_c v] \>| \nonumber \\
&\leq\int_s^t \frac{C}{(1+|t-\tau|)^{\frac{N}{2}}}\|Dh|_a(\tau)\|_{\Ls}[|\< \po,g_u \ehts P_c v \>|+|\< \po,g_{\bar{u}} e^{iH(\tau-s)} P_c \bar{v} \>|] \nonumber \\
&\leq\int_s^t \frac{C\|Dh\|}{(1+|t-\tau|)^{\frac{N}{2}}}[|\< e^{iH(\tau-s)}\po,e^{iH(\tau-s)}g_u \ehts P_c v \>| \nonumber \\
&\quad\quad\quad\quad\quad\quad\quad\quad\quad\quad\quad\quad\quad\quad+|\<\ehts \po,\ehts g_{\bar{u}} e^{iH(\tau-s)} P_c \bar{v} \>|] \nonumber \\
&\leq\int_s^t \frac{C\|Dh\|}{(1+|t-\tau|)^{\frac{N}{2}}}\|e^{iE_0(\tau-s)}\po\|_{L^{q_2}}\sup_{\tau\in[s,t]}(\|\widehat{g_u(\tau)}\|_{L^1}+\|\widehat{g_{\bar{u}}(\tau)}\|_{L^1})\|v\|_{L^{q_2'}} \nonumber \\
&\leq
\frac{C\|Dh\|\sup_{\tau\in[s,t]}(\|\widehat{g_u(\tau)}\|_{L^1}+\|\widehat{g_{\bar{u}}(\tau)}\|_{L^1})\|\po\|_{L^{q_2}}\|v\|_{L^{q_2'}}}{(1+|t-s|)^{\frac{N}{2}-1}}\leq
\frac{C\|v\|_{L^{q_2'}}}{(1+|t-s|)^{\frac{N}{2}-1}} \nonumber
\end{align} where we also used $\psi_0\in H^2\hookrightarrow L^{q_2}.$ Similarly $\|\tilde{f}\|_{\Lsn}\leq\frac{C\|v\|_{L^{q_1'}}}{(1+|t-s|)^{\frac{N}{2}-1}}\leq C\|v\|_{L^{q_1'}}$ for $q_1=\frac{2N}{N-2}$.

For $t>s+1$ we will split these two integral in two parts to be
estimated differently:
$$f(t)=\underbrace{\int_s^{s+\frac{1}{4}}}_{I_{1}}+\underbrace{\int_{s+\frac{1}{4}}^{t}-i\eht P_c\Dg[\ehts P_c v]d\tau}_{I_{2}}$$ and
$$\tilde{f}(t)=\underbrace{\int_s^{s+\frac{1}{4}}}_{II_{1}}+\underbrace{\int_{s+\frac{1}{4}}^{t}\eht P_c Dh|_{a(\tau)}i\<\po,\Dg[\ehts P_c v]\>d\tau}_{II_{2}}.$$
Then we have:
\begin{align} \|I_1\|_{\Lsn}&\leq\|\xsn\int_s^{s+\frac{1}{4}}-i \ehs P_c \Dg[\ehts P_c v] d\tau\|_{L^2} \nonumber \\
&\leq\|\xsn\|_{L^\beta}\int_s^{s+\frac{1}{4}} \|\ehs e^{iH(\tau-s)} P_c g_u\ehts P_c v\|_{L^{q_2}} d\tau \nonumber \\
&\quad+\|\xsn\|_{L^\beta}\int_s^{s+\frac{1}{4}} \|e^{-iH(t+s-2\tau)} \ehts P_c g_{\bar{u}}e^{iH(\tau-s)} P_c v\|_{L^{q_2}} d\tau \nonumber \\
&\leq\|\xsn\|_{L^\beta}\int_s^{s+\frac{1}{4}} \frac{C_{q_2}}{|t-s|^{N(\frac{1}{2}-\frac{1}{q_2})}}\|e^{iH(\tau-s)} P_c g_u\ehts P_c v\|_{L^{q_2'}}d\tau \nonumber \\
&\quad+\|\xsn\|_{L^\beta}\int_s^{s+\frac{1}{4}} \frac{C_{q_2}}{|t+s-2\tau|^{N(\frac{1}{2}-\frac{1}{q_2})}}\|\ehts P_c g_{\bar{u}}e^{iH(\tau-s)} P_c v\|_{L^{q_2'}} d\tau \nonumber \\
&\leq C(q_2)\|\xsn\|_{L^\beta}\Big(\frac{1}{|t-s|^{N(\frac{1}{2}-\frac{1}{q_2})}}+\frac{1}{|t-s-\frac{1}{2}|^{N(\frac{1}{2}-\frac{1}{q_2})}}\Big)\int_s^{s+\frac{1}{4}} (\|\widehat{g_u}(\tau)\|_{L^1}+\|\widehat{g_{\bar{u}}}(\tau)\|_{L^1})\|v\|_{L^{q_2'}}d\tau \nonumber \\
&\leq C(q_2)\|\xsn\|_{L^\beta}\sup_{\tau\in\R}(\|\widehat{g_u}(\tau)\|_{L^1}+\|\widehat{g_{\bar{u}}}(\tau)\|_{L^1})\Big(\frac{1}{|t-s|^{N(\frac{1}{2}-\frac{1}{q_2})}}+\frac{1}{|t-s-\frac{1}{2}|^{N(\frac{1}{2}-\frac{1}{q_2})}}\Big)\|v\|_{L^{q_2'}} \nonumber \\
&\leq
C(q_2)\frac{\|v\|_{L^{q_2'}}}{(1+|t-s|)^{N(\frac{1}{2}-\frac{1}{q_2})}}
\nonumber
\end{align} For the second integral we have \begin{align}
\|I_2\|_{\Lsn}&\leq\int_{s+\frac{1}{4}}^{t} \|\eht P_c\|_{\Ls\rightarrow\Lsn}\|\Dg[\ehts P_c v]\|_{\Ls} d\tau \nonumber \\
&\leq\int_{s+\frac{1}{4}}^{t} \frac{C_N}{(1+|t-\tau|)^{\frac{N}{2}}}\|\xs(|\pe|^{\p_1}+|\pe|^{\p_2})\|_{L^\beta}\|\ehts P_c v\|_{L^{q_2}}d\tau \nonumber \\
&\leq
C(q_2)\|v\|_{L^{q_2'}}\int_{s+\frac{1}{4}}^{t}\frac{d\tau}{(1+|t-\tau|)^{\frac{N}{2}}
|\tau-s|^{N(\frac{1}{2}-\frac{1}{q_2})}}\leq\frac{C(q_2)\|v\|_{L^{q_2'}}}{(1+|t-s|)^{N(\frac{1}{2}-\frac{1}{q_2})}}
\nonumber \end{align} For the second forcing term $\tilde{f}(t)$, we
use again $\psi_0\in H^2\hookrightarrow L^{q_2}:$
\begin{align}
\|II_1\|_{\Lsn}&\leq\int_s^{s+\frac{1}{4}} \frac{C}{(1+|t-\tau|)^{\frac{N}{2}}}\|Dh|_a(\tau)\|_{\Ls}|\< \po, \Dg[\ehts P_c v] \>|d\tau \nonumber \\
&\leq\frac{C}{(1+|t-s-\frac{1}{4}|)^{\frac{N}{2}}}\int_s^{s+\frac{1}{4}}\Big[ |\<e^{iH(\tau-s)}\po,e^{iH(\tau-s)}g_u\ehts P_c v \>| \nonumber \\
&\quad\quad\quad\quad\quad\quad\quad\quad\quad\quad\quad\quad\quad\quad\quad\quad\quad\quad+|\< \ehts\po,\ehts g_{\bar{u}}e^{iH(\tau-s)} P_c \bar{v}\>|\Big] d\tau \nonumber \\
&\leq\frac{C\sup_{\tau\in\R}(\|\widehat{g_u}(\tau)\|_{L^1}+\|\widehat{g_{\bar{u}}}(\tau)\|_{L^1})\|v\|_{L^{q'_2}}}{(1+|t-s|)^{\frac{N}{2}}}\int_s^{s+\frac{1}{4}} \|e^{\pm iE_0(\tau-s)}\po\|_{L^{q_2}} d\tau \nonumber \\
&\leq\frac{C\|v\|_{L^{q'_2}}}{(1+|t-s|)^{\frac{N}{2}}} \nonumber
\end{align}

$II_2$ is estimated exactly the same way as $I_2.$ Let us observe
that the above estimates are for the case $t>s+1$. Because of that
we can replace the $C/|t-s|$ term by $C/(1+|t-s|)$ in the  $I_1$,
$I_2$ and $II_2$ integrals. The estimates for $s-1\leq t\leq s$ respectively $t<s-1$ are obtained in the same way as the ones for $s\leq t\leq s+1$ respectively $t>s+1.$

Theorem \ref{th:lw} is now completely proven. $\Box$

The next step is to obtain estimates for $\om(t,s)$ and $T(t,s)$ in unweighted $L^p$ spaces.
\begin{theorem}{\label{th:lp}} Fix $\sigma > N/2$ and
$\frac{2N}{N-2}<q_2<\frac{2N}{N-4}.$ Assume that
$\|<x>^{4\sigma}\psi_{E(t)}\|_{L^\infty} < \varepsilon_1(q_2),$ for
all $t\in\R$ (where $\varepsilon_1(q_2)$ is the one used in Theorem
\ref{th:lw}). Then there exist constants $C_2$, $C_2'$ and $C_p$
such that for all $t$, $s$ $\in$ $\R$ the following estimates hold:
\begin{align}
(i)\quad &\|\om(t,s)\|_{L^2\rightarrow L^2}\leq C_2,\quad \|T(t,s)\|_{L^2\rightarrow L^2}\leq C_2 \nonumber \\
(ii)\quad &\|T(t,s)\|_{L^{p'}\rightarrow L^p}\leq\left\{ \begin{array}{ll}
\frac{C_p}{|t-s|^{N(1-2/p)-1}} & \textrm{for $|t-s|\leq 1$} \\
\frac{C_p}{|t-s|^{N\ip}} & \textrm{for $|t-s|>1$} \end{array} \right.,\quad\textrm{for}\quad 2\leq p\leq q_2<\frac{2N}{N-4} \nonumber \\
& \|\om(t,s)\|_{L^{p'}\rightarrow L^p}\leq
\frac{C_p}{|t-s|^{N\ip}},\quad\textrm{for}\quad 2\leq p\leq q_2<\frac{2N}{N-4} \nonumber \\
 (iii)\quad &\|T(t,s)-\tilde{T}(t,s)\|_{L^{p'}\rightarrow L^2}\leq
C_2', \quad\textrm{for }\quad 2\leq p\leq\frac{2N}{N-2} \nonumber
\end{align} where
$$\tilde{T}(t,s)=\left\{\begin{array}{ll}-i\int_s^{\min\{t,
s+1\}}e^{-iH(t-\tau
)}P_c g_u e^{iH(\tau-s)}P_c d\tau, & \textrm{if}\quad t\geq s,\\
 -i\int_s^{\max\{t, s-1\}}e^{-iH(t-\tau
)}P_c g_u e^{iH(\tau-s)}P_c d\tau , & \textrm{if}\quad t<s.\end{array}\right.$$
\end{theorem}

%\begin{remark}\label{rmk:lp}
% By Riesz-Thorin interpolation from $(i)$ and $(ii)$, and from $(i)$ and $(iii)$ we get
%\begin{eqnarray}
%\|T(t,s)\|_{L^{p'}\rightarrow L^p}&\leq&\frac{C_p}{|t-s|^{N\ip}},\quad\textrm{for all}\quad 2\leq p\leq\frac{2N}{N-2} \nonumber \\
%\|\om(t,s)\|_{L^{p'}\rightarrow L^p}&\leq&\frac{C_p}{|t-s|^{N\ip}},\quad\textrm{for all}\quad 2\leq p<\frac{2N}{N-2} \nonumber \\
%\|T(t,s)\|_{L^{p'}\rightarrow L^p}&\leq&\left\{ \begin{array}{ll}
%\frac{C_p}{|t-s|^3} & \textrm{for $|t-s|\leq 1$} \\
%\frac{C}{|t-s|^{N\ip}} & \textrm{for $|t-s|>1$} \end{array} \right.,\quad\textrm{for all}\quad \frac{2N}{N-2}<p<\frac{2N}{N-4} \nonumber \\
%\|T(t,s)-\tilde{T}(t,s)\|_{L^{p'}\rightarrow L^2}&\leq &C_p,
%\quad\textrm{for all}\quad 2\leq p\leq\frac{2N}{N-2}
%\nonumber \end{eqnarray}
%\end{remark}

\noindent\textbf{Proof of Theorem \ref{th:lp}} Fix $s\in\R .$
Because of the estimate (\ref{eq:lp}) and relation
$\om(t,s)=T(t,s)+\ehs P_c$, it suffices to prove the theorem for
$T(t,s)$. Throughout this proof we will repeatedly use the equations
defining $T(t,s)$ \eqref{w}-\eqref{sow} where the linear operator
$L(s)$ is given by \eqref{L}. Note that we have already denoted the
remaining forcing terms in \eqref{sow} by $f$ respectively $\tilde
f$ see \eqref{def:f} and \eqref{def:tildef}.

$(i)$ To estimate the $L^2$ norm we will use the following duality
argument: \begin{align}
&\|f(t)\|_{L^2}^2=\< f(t),f(t)\> \nonumber \\
&=\int_s^t\int_s^t \< \eht P_c \Dg[\ehts P_c v],e^{-iH(t-\tau')} P_c \Dg[e^{-iH(\tau' -))} P_c v]\> d\tau' d\tau  \nonumber \\
&=\int_s^t \int_s^t \< \Dg[\ehts P_c v],e^{-iH(\tau-\tau')} P_c \Dg[e^{-iH(\tau'-s)} P_c v]\> d\tau' d\tau  \nonumber \\
&=\int_s^t \int_s^t \< \xs \Dg[\ehts P_c v], \xsn e^{-iH(\tau-\tau')} P_c \Dg[e^{-iH(\tau'-s)} P_c v] \> d\tau' d\tau  \nonumber \\
&\leq\int_s^t \int_s^t \|\Dg[\ehts P_c v]\|_{\Ls}\|e^{-iH(\tau-\tau')} P_c \Dg[e^{-iH(\tau'-s)} P_c v]\|_{\Lsn} d\tau' d\tau \nonumber \\
&\leq\int_s^t \|\xs \Dg[\ehts P_c v]\|_{L^2}\int_s^t\frac{C}{(1+|\tau-\tau'|)^{\frac{N}{2}}}\|\xs \Dg[e^{-iH(\tau'-s)} P_c v]\|_{L^2} d\tau' d\tau\nonumber \\
&\leq C\|\xs \Dg[\ehts P_c v]\|_{L_{\tau}^2 L_x^2}
\Big\|\int_s^t\frac{C}{(1+|\tau-\tau'|)^{\frac{N}{2}}} \|\xs \Dg[e^{-iH(\tau'-s)} P_c v]\|_{L_x^2} d\tau'\Big\|_{L^2_{\tau}}\nonumber \\
&\leq C\|K\|_{L^1}\|\xs \Dg[e^{-iHt} P_c v]\|_{L_t^2 L^2_x}^2\leq
C\|v\|^2_{L^2}<\infty \nonumber \end{align} At the last line we used
Young inequality $\|K*f\|_{L^2}\leq\|K\|_{L^1}\|f\|_{L^2}$ with
$K(t)=(1+|t|)^{-\frac{N}{2}}\in L^1$ and  for the term $\xs
\Dg[e^{-iHt} P_c v]=\xs (g_u e^{-iHt} P_c v+g_{\bar{u}} e^{iHt} P_c
v)$ we used the Kato smoothing estimate $\|\xsn e^{-iHt} P_c
v\|_{L^2_t(\R,L^2_x)}\leq C\|v\|_{L^2_x},$ together with the uniform
bounds in time $\|\xst g_u(\tau')\|_{L^\infty},\ \|\xst
g_{\bar{u}}(\tau')\|_{L^\infty}\leq C \varepsilon_1^{\alpha_1}$
since $|g_u(\tau')|,\ |g_{\bar{u}}(\tau')|\leq
C(|\psi_{E(\tau')}|^{\p_1}+|\psi_{E(\tau')}|^{\p_2}).$ Similarly we
have,
\begin{align}
\|\tilde{f}\|_{L^2}^2=&\int_s^t \int_s^t \< \eht P_c Dh|_{a(\tau)} \<\po, \Dg[\ehts P_c v]\>, \nonumber\\
&\quad\quad\quad\quad e^{-iH(t-\tau')} P_c Dh|_{a(\tau)} \<\po, \Dg[e^{-iH(\tau' -s)} P_c v]\> \> d\tau' d\tau  \nonumber \\
&\leq\int_s^t \int_s^t \|Dh\|_{\Ls}\|\po\|_{L^2}\|\Dg[\ehts P_c v]\|_{L^2} \nonumber \\
&\quad\quad\times\frac{C}{(1+|\tau-\tau'|)^{N/2}}\|Dh\|_{\Ls}\|\po\|_{L^2}\|\Dg[e^{-iH(\tau'-s)} P_c v]\|_{L^2} d\tau' d\tau  \nonumber \\
&\leq C\|v\|^2_{L^2}<\infty \nonumber \end{align}

We will estimate $L^2$ norm of $L(s)$ see \eqref{L} in the same way as for $f:$
\begin{align}
&\|L(s)W\|_{L^2}^2=\<L(s)W,L(s)W\> \nonumber \\
%&=\<\int_s^t \eht P_c [\Dg[W(\tau))-Dh\<\po,\Dg[W(\tau))\>]d\tau,\int_s^t e^{-iH(t-\tau')} P_c [\Dg[W(\tau'))-Dh\<\po,\Dg[W(\tau'))\>d\tau'\> %]\nonumber \\
&=\int_s^t \int_s^t \left\< \Dg[W(\tau)]-Dh\<\po,\Dg[W(\tau)]\>,e^{-iH(\tau-\tau')} P_c \left(\Dg[W(\tau')]-Dh\<\po,\Dg[W(\tau')]\>\right)\right\> d\tau' d\tau  \nonumber \\
&\leq\int_s^t (\|\xs g_u\|_{L^\infty_x}+\|\xs g_{\bar{u}}\|_{L^\infty_x})(1+\|Dh\|_{\Ls}\|\po\|_{L^2})\|\xsn W\|_{L^2_x} \nonumber \\
&\quad\quad\quad\quad\quad\quad\times\int_s^t C K(\tau-\tau')(\|\xs g_u\|_{L^\infty_x}+\|\xs g_{\bar{u}}\|_{L^\infty_x})(1+\|Dh\|_{\Ls}\|\po\|_{L^2})\|\xsn W\|_{L^2_x} d\tau' d\tau\nonumber \\
%&\leq C\|\xsn W\|_{L_x^2 L_\tau^2}
%\Big\|\int_s^t CK(\tau-\tau') \|\xsn W\|_{L_x^2} d\tau'\Big\|_{L^2_{\tau}}\nonumber \\
&\leq C\|K\|_{L^1}\|\xsn W\|_{L_{\tau}^2 L^2_x}^2<\infty \nonumber \end{align}
 By Theorem \ref{th:lw} $(iii)$, $\|\xsn W\|_{L_\tau^2 L^2_x}<\infty$.

\noindent Therefore we conclude $\|T(t,s)\|_{L^2\rightarrow L^2}\leq C$ and $\|\om(t,s)\|_{L^2\rightarrow L^2}\leq C$

$(ii)$ It suffices to prove the estimates for $p=2$ and $p=q_2.$ The
estimates for $2<p<q_2$ will follow from Riesz-Thorin interpolation.
We will also assume $t>s$ since the case $t<s$ can be treated
similarly.

We start with short time estimates, $s<t\leq s+1$ where the difficult part is
to remove the non-integrable singularities of  $e^{-iHt}P_c
$ at $t=0,$ see \eqref{eq:lp} for $p>2N/(N-2),$ which appears in the
convolution integrals in \eqref{sow}. For this purpose we will use
the Fourier multiplier estimates \eqref{est:fm}. Let us first investigate
the short time behavior of the terms $f(t)$ and $\tilde{f}(t)$. In what follows
$p=2$ and $p=q_2>2N/(N-2):$
\begin{align}
\|f(t)\|_{L^p}&=\|\int_s^t \eht P_c \Dg[\ehts P_c v]d\tau\|_{L^p} \nonumber \\
&\leq\int_s^t \|e^{-iH(t-s)}\|_{L^{p'}\rightarrow L^p}\|e^{iH(\tau-s)}P_c g_u \ehts P_c v\|_{L^{p'}}d\tau \nonumber \\
&+\int_s^{s+\frac{t-s}{4}} \|e^{-iH(t+s-2\tau)}\|_{L^{p'}\rightarrow L^p}\|\ehts P_c g_{\bar{u}} e^{iH(\tau-s)}  P_c \bar{v}\|_{L^{p'}}d\tau \nonumber \\
&+\int_{s+\frac{t-s}{4}}^{t-\frac{t-s}{4}} \|e^{-iH(t-\tau)}P_c\|_{L^{p'}\rightarrow L^p}\|g_{\bar{u}} e^{iH(\tau-s)}  P_c \bar{v}\|_{L^{p'}}d\tau \nonumber \\
&+\int_{t-\frac{t-s}{4}}^t \|e^{-iH(t+s-2\tau)}\|_{L^{p'}\rightarrow L^p}\|\ehts P_c g_{\bar{u}} e^{iH(\tau-s)}  P_c \bar{v}\|_{L^{p'}}d\tau \nonumber \\
&\leq\int_s^t \frac{C_p}{|t-s|^{N\ip}}C\|\widehat{g_u}(\tau)\|_{L^1}\|v\|_{L^{p'}}d\tau+\int_s^{s+\frac{t-s}{4}} \frac{C_p}{|t+s-2\tau|^{N\ip}}C\|\widehat{g_{\bar{u}}}(\tau)\|_{L^1}\|v\|_{L^{p'}}d\tau \nonumber \\
&+\int_{s+\frac{t-s}{4}}^{t-\frac{t-s}{4}}
\frac{C_p}{|t-\tau|^{N\ip}}\|g_{\bar{u}}(\tau)\|_{L^{\frac{p}{p-2}}}\frac{C_p\|v\|_{L^{p'}}}{|\tau-s|^{N\ip}}d\tau+\int_{t-\frac{t-s}{4}}^t
\frac{C_p}{|t+s-2\tau|^{N\ip}}C
\|\widehat{g_{\bar{u}}}(\tau)\|_{L^1}\|v\|_{L^{p'}}d\tau \nonumber \\
&\leq \frac{\tilde
C\max\{\sup_{\tau\in\R}\|g_{\bar{u}}(\tau)\|_{L^{\frac{p}{p-2}}},\sup_{\tau\in\R}\|\widehat{g_u}(\tau)\|_{L^1},\sup_{\tau\in\R}\|\widehat{g_{\bar{u}}}(\tau)\|_{L^1}\}\|v\|_{L^{p'}}}{|t-s|^{N(1-\frac{2}{p})-1}}
\nonumber \end{align} For $\tilde f$ we use the Sobolev imbeddings
$\|F\|_{L^p}\leq C\|F\|_{H^2}$ for $2\leq p\leq\frac{2N}{N-4}$
together with $\|e^{-iHt} F\|_{H^2}\leq C\|F\|_{H^2}:$
\begin{align}
\|\tilde{f}(t)\|_{L^p}&\leq\int_s^t C\|Dh|_{a(\tau)}\|_{\C\rightarrow H^2}|\<\po,\Dg[\ehts P_c v]\>|d\tau \nonumber \\
&\leq\int_s^t \|Dh|_{a(\tau)}\|_{\C\rightarrow H^2}[\
|\<e^{iH(\tau-s)}\po,e^{iH(\tau-s)}g_u\ehts P_c v\>|\nonumber\\
&+|\<\ehts\po,\ehts g_{\bar{u}}e^{iH(\tau-s)} P_c \bar{v}\>|\ ]d\tau \nonumber \\
&\leq\int_s^t \|Dh|_{a(\tau)}\|_{{\C\rightarrow H^2}}\|e^{iH(\tau-s)}\po\|_{L^p}\|e^{iH(\tau-s)}g_u\ehts P_c v\|_{L^{p'}}d\tau \nonumber \\
&+\int_s^t \|Dh|_{a(\tau)}\|_{{\C\rightarrow H^2}} \|\ehts\po\|_{L^p}\|\ehts g_{\bar{u}}e^{iH(\tau-s)} P_c \bar{v}\|_{L^{p'}} d\tau \nonumber \\
&\leq C\sup_{\tau\in\R}\|Dh|_{a(\tau)}\|_{\C\rightarrow
H^2}(\sup_{\tau\in\R}\|\widehat{g_u}(\tau)\|_{L^1}+\sup_{\tau\in\R}\|\widehat{g_{\bar{u}}}(\tau)\|_{L^1})\|\po\|_{H^2}|t-s|\
\|v\|_{L^{p'}}\leq \frac{\tilde C
\|v\|_{L^{p'}}}{|t-s|^{N(1-\frac{2}{p})-1}} \nonumber \end{align}

We now move to the short time estimate: $s<t\leq s+1$ of
$\|L(s)W\|_{L^p},\ p=2,\ q_2,$ see \eqref{L} for the definition of
the integral operator $L(s).$ The main difference compared to the
$f$ and $\tilde f$ terms is the fact that the singularity at
$\tau=s$ is integrable due to Theorem \ref{th:lw} part (iv):
 $$[L(s)W](t)=-i\int_s^t \eht P_c \Dg[W(\tau)] d\tau+\int_s^t \eht
 P_cDh|_{a(\tau)}i\<\po,\Dg[W(\tau)]\>
 d\tau$$ where
 \begin{eqnarray}
 \lefteqn{\|\int_s^t \eht P_cDh|_{a(\tau)}i\<\po,\Dg[W(\tau)]\>
 d\tau\|_{L^p}}\nonumber\\
 &\leq &\int_s^t C\|\eht P_c\|_{H^2\mapsto H^2}\|Dh|_{a(\tau)}\|_{\C\mapsto H^2}|\<\po,\Dg[W(\tau)]\>|
 d\tau\nonumber\\
 &\leq & \int_s^t C\|Dh|_{a(\tau)}\|_{\C\mapsto H^2}\|\po\|_{L^2}(\|\xs g_u\|_{L^\infty}+\|\xs g_{\bar
 u}\|_{L^\infty})\|W(\tau)\|_{\Lsn}d\tau\nonumber\\
 &\leq &\|\po\|_{L^2}\sup_{\tau\in\R}\{\|Dh|_{a(\tau)}\|_{\C\mapsto H^2}(\|\xs g_u\|_{L^\infty}+\|\xs g_{\bar
 u}\|_{L^\infty})\}\int_s^t\frac{C_p\|v\|_{L^{p'}}}{|\tau-s|^{N\ip -1}}d\tau\leq \frac{C\|v\|_{L^{p'}}}{|t-s|^{N\ip
-1}}\nonumber
 \end{eqnarray} However to remove the non-integrable singularity at $\tau=t$ in the
remaining integral we need to plug in (\ref{sow}) in it:
\begin{eqnarray}
\lefteqn{\int_s^t \eht P_c \Dg[W(\tau)] d\tau}\nonumber \\
&=&\int_s^t \eht P_c g_u\Big(\int_s^\tau e^{-iH(\tau-\tau')} P_c(-i\Dg[e^{-iH(\tau'-s)}P_c v]+Dhi\<\po,\Dg[e^{-iH(\tau'-s)}P_c v]\>)d\tau'\nonumber \\
&&\quad\quad\quad\quad\quad\quad\quad\quad\quad\quad\quad\quad+\int_s^\tau e^{-iH(\tau-\tau')}P_c (-i \Dg[W(\tau')]+Dhi\<\po,\Dg[W(\tau')]\>)d\tau'\Big)d\tau \nonumber \\
&&+\int_s^t \eht P_c g_{\bar{u}}\Big(\int_s^\tau e^{iH(\tau-\tau')}P_c(\overline{-i \Dg[e^{-iH(\tau'-s)}P_c v]}+\overline{Dhi\<\po,\Dg[e^{-iH(\tau'-s)}P_c v]\>})d\tau' \nonumber \\
&&\quad\quad\quad\quad\quad\quad\quad\quad\quad\quad\quad\quad+\int_s^\tau
e^{iH(\tau-\tau')} P_c(\overline{-i
\Dg[W(\tau')]}+\overline{Dhi\<\po,\Dg[W(\tau')]\>})d\tau'\Big)d\tau
\nonumber \end{eqnarray} All the terms will be either of the
following forms
$$L_1=\int_s^t \eht P_c g_u\int_{s}^\tau e^{-iH(\tau-\tau')} P_c X(\tau')d\tau'd\tau$$
$$L_2=\int_s^t \eht P_c g_{\bar{u}}\int_{s}^\tau e^{iH(\tau-\tau')} P_c \overline{X(\tau')}d\tau'd\tau$$
$$\tilde{L}_1=\int_s^t \eht P_c g_u\int_{s}^\tau e^{-iH(\tau-\tau')} P_c Dh\<\po,X(\tau')\>d\tau'd\tau$$
$$\tilde{L}_2=\int_s^t \eht P_c g_{\bar{u}}\int_{s}^\tau e^{iH(\tau-\tau')} P_c \overline{Dh\<\po,X(\tau')\>}d\tau'd\tau$$
where $X(\tau')$ can be either of $-ig_u e^{-iH(\tau'-s)}P_c v,
-ig_{\bar{u}} e^{iH(\tau'-s)}P_c \bar{v}, -ig_u W(\tau'),$ or $
-ig_{\bar{u}}\overline{ W(\tau')}.$ We can now remove the
singularity of $\|\eht P_c\|_{L^{p'}\mapsto L^p}$ at $\tau=t$ via
\eqref{est:fm}:
\begin{equation}
L_1=\int_s^t \eht P_c g_u e^{iH(t-\tau)}\int_{s}^\tau e^{-iH(t-\tau')} P_c X(\tau')d\tau'd\tau \label{l1} \end{equation}
\begin{equation}L_2=\int_s^t \eht P_c g_{\bar{u}} e^{iH(t-\tau)}\int_{s}^\tau e^{-iH(t-2\tau+\tau')} P_c \overline{X(\tau')}d\tau'd\tau \label{l2} \end{equation}
\begin{equation} \tilde{L}_1=\int_s^t \eht P_c g_u e^{iH(t-\tau)}\int_{s}^\tau e^{-iH(t-\tau')} P_c Dh\<\po,X(\tau')\>d\tau'd\tau \label{lt1} \end{equation}
\begin{equation} \tilde{L}_2=\int_s^t \eht P_c g_{\bar{u}} e^{iH(t-\tau)}\int_{s}^\tau e^{-iH(t-2\tau+\tau')} P_c \overline{Dh\<\po,X(\tau')\>}d\tau'd\tau \label{lt2} \end{equation}
\begin{itemize}
\item For $X(\tau')=-ig_u e^{-iH(\tau'-s)}P_c v$ we have:
\begin{align}
\|L_1\|_{L^p}&\leq\int_s^t \|\eht P_c g_u e^{iH(t-\tau)}\|_{L^p\rightarrow L^p}\int_s^{\tau}\|e^{-iH(t-s)}P_c\|_{L^{p'}\rightarrow L^p} \|e^{iH(\tau'-s)}P_c g_u e^{-iH(\tau'-s)}P_c v\|_{L^{p'}}d\tau'd\tau \nonumber \\
&\leq\int_s^t \|\widehat{g_u}\|_{L^1}\int_s^\tau \frac{C}{|t-s|^{N\ip}}\|\widehat{g_u}\|_{L^1}\|v\|_{L^{p'}}d\tau'd\tau\leq \frac{C\|v\|_{L^{p'}}}{|t-s|^{N\ip-2}} \nonumber \end{align}
\begin{align}
\|L_2\|_{L^p}&\leq\int_s^{t-\frac{t-s}{4}} \|\eht P_c\|_{L^{p'}\rightarrow L^p}\|g_{\bar{u}}\int_s^\tau e^{-iH(\tau-s)}e^{-iH(\tau'-s)}P_c g_{\bar u} e^{iH(\tau'-s)}P_c \bar{v}\|_{L^{p'}}d\tau' d\tau \nonumber \\
&+\int_{t-\frac{t-s}{4}}^t \|\eht P_c g_{\bar{u}} e^{iH(t-\tau)}\|_{L^p\rightarrow L^p} \nonumber \\
&\quad\quad\quad\quad\times \int_s^{\tau}\|e^{-iH(t+s-2\tau)}P_c\|_{L^{p'}\rightarrow L^p}\|e^{-iH(\tau'-s)}P_c g_{\bar{u}} e^{iH(\tau'-s)}P_c \bar{v}\|_{L^{p'}}d\tau'd\tau \nonumber \\
&\leq\int_s^{t-\frac{t-s}{4}} \frac{C}{|t-\tau|^{N\ip}}\|g_{\bar{u}}\|_{L^{\frac{p}{p-2}}}\int_s^\tau\frac{C}{|\tau-s|^{N\ip}}\|\widehat{g_{\bar{u}}}\|_{L^1}\|v\|_{L^{p'}}d\tau' d\tau \nonumber \\
&+\int_{t-\frac{t-s}{4}}^t \|\widehat{g_{\bar{u}}}\|_{L^1}\int_s^\tau \frac{C}{|t+s-2\tau|^{N\ip}}\|\widehat{g_{\bar{u}}}\|_{L^1}\|v\|_{L^{p'}}d\tau'd\tau\leq \frac{C\|v\|_{L^{p'}}}{|t-s|^{2N\ip-2}} \nonumber \end{align}
\begin{align}
\|\tilde{L}_{1,2}\|_{L^p}&\leq\int_s^t \|e^{-iH(t-\tau)} P_c g_u e^{iH(t-\tau)}\|_{L^p\rightarrow L^p}\int_s^{\tau} C\|Dh\|_{\C\rightarrow H^2}|\<\po, g_u e^{-iH(\tau'-s)}P_c v\>|d\tau'd\tau \nonumber \\
&\leq\int_s^t \|\widehat{g_u}(\tau)\|_{L^1}\int_s^\tau C|\<e^{iH(\tau'-s)}\po,e^{iH(\tau'-s)}g_u e^{-iH(\tau'-s)}P_c v\>|d\tau'd\tau \nonumber \\
&\leq\int_s^t \|\widehat{g_u}(\tau)\|_{L^1}\int_s^\tau C\|e^{iH(\tau'-s)}\po\|_{L^p}\|e^{iH(\tau'-s)}g_u e^{-iH(\tau'-s)}P_c v\|_{L^{p'}}d\tau'd\tau \nonumber \\
&\leq C|t-s|^2 \|v\|_{L^{p'}} \nonumber \end{align}
\item For $X(\tau')=-ig_{\bar{u}}e^{iH(\tau'-s)}P_c \bar{v}$ in $L_1$ we first change the order of integration then split and use \eqref{est:fm}:
\begin{align}
\|L_1\|_{L^p}&\leq\int_s^{s+\frac{t-s}{4}}\int_{\tau'}^{t} \|\eht P_c g_u e^{iH(t-\tau)} e^{-iH(t+s-2\tau')} P_c e^{-iH(\tau'-s)}P_c g_{\bar{u}} e^{iH(\tau'-s)}P_c \bar{v}\|_{L^p}d\tau d\tau'\nonumber \\
&+\int_{s+\frac{t-s}{4}}^{t}\int_{\tau'}^t\|\eht P_c g_u e^{iH(t-\tau)}\|_{L^p\mapsto L^p}\|e^{-iH(t-\tau')}\|_{L^{p'}\rightarrow L^p}\|g_{\bar{u}} e^{iH(\tau'-s)}P_c \bar{v}\|_{L^{p'}}d\tau d\tau' \nonumber \\
&\leq\int_s^{s+\frac{t-s}{4}}\int_{\tau'}^{t} \|\widehat{g_u}\|_{L^1}\frac{C\|v\|_{L^1}}{|t+s-2\tau'|^{N\ip}}\|\widehat{g_{\bar{u}}}\|_{L^1}d\tau d\tau'\nonumber \\
&+\int_{s+\frac{t-s}{4}}^{t}\int_{\tau'}^t\|\widehat{g_u}\|_{L^1}\frac{C}{|t-\tau'|^{N\ip}}\|g_{\bar{u}}\|_{L^\beta}\frac{C\|v\|_{L^{p'}}}{|\tau'-s|^{N\ip}}d\tau d\tau '\nonumber \\
&\leq \frac{C\|v\|_{L^{p'}}}{|t-s|^{2N\ip-2}} \nonumber\end{align}
For $L_2$ we do not change the order of integration but we have to split both integrals to avoid singularities:
\begin{align}
\|L_2\|_{L^p}&\leq\int_s^{t-\frac{t-s}{16}} \|\eht P_c\|_{L^{p'}\rightarrow L^p}\nonumber \\
&\quad\quad\quad\times\Big[ \|g_{\bar{u}}\|_{L^\beta}\int_s^{s+\frac{\tau-s}{4}} \|e^{-iH(2\tau'-\tau-s)}\|_{L{p'}\rightarrow L^p}\|e^{iH(\tau'-s)} P_c g_u e^{-iH(\tau'-s)}P_c v\|_{L^{p'}}d\tau'\|_{L^p} \nonumber \\
&\quad\quad\quad+\|g_{\bar{u}}\|_{L^{2\beta}}\int_{s+\frac{t-s}{4}}^\tau \|e^{-iH(\tau'-\tau)}P_c\|_{L^2\rightarrow L^2}\|g_u\|_{L^{2\beta}}\|e^{-iH(\tau'-s)}P_c v\|_{L^{p'}}d\tau'\Big]d\tau \nonumber \\
&+\int_{t-\frac{t-s}{16}}^t \|\eht P_c g_{\bar{u}} e^{iH(t-\tau)}\|_{L^p\rightarrow L^p} \nonumber \\
&\quad\quad\quad\quad\times\Big[\int_{s}^{s+\frac{t-s}{4}}\| e^{-iH(t-2\tau+2\tau'-s)}e^{iH(\tau'-s)} P_c g_u R_a e^{-iH(\tau'-s)}P_c v\|_{L^p}d\tau' \nonumber \\
&\quad\quad\quad\quad\quad\quad\quad\quad+\int_{s+\frac{t-s}{4}}^{t-\frac{t-s}{4}} \|e^{-iH(t-2\tau+\tau')} P_c g_u e^{-iH(\tau'-s)}P_c v\|_{L^p}d\tau' \nonumber \\
&\quad\quad\quad\quad\quad\quad\quad\quad+\int_{t-\frac{t-s}{4}}^\tau\| e^{-iH(t-2\tau+2\tau'-s)}e^{iH(\tau'-s)} P_c g_u e^{-iH(\tau'-s)}P_c v\|_{L^p}d\tau' \Big]d\tau \nonumber \\
&\leq\int_s^t \|\widehat{g_{\bar{u}}}\|_{L^1}\Big[\int_s^{s+\frac{t-s}{4}} \frac{C}{|t-2\tau+2\tau'-s|^{N\ip}}\|\widehat{g_{\bar{u}}}\|_{L^1} \|v\|_{L^{p'}}d\tau'\nonumber \\ &\quad\quad\quad\quad\quad\quad\quad\quad+\int_{s+\frac{t-s}{4}}^\tau\frac{C}{|t-2\tau+\tau'|^{N\ip}}\|g_{\bar{u}}\|_{L^{\frac{p}{p-2}}}\frac{\|v\|_{L^{p'}}}{|\tau'-s|^{N\ip}}d\tau'\Big]d\tau \nonumber \\
&+\int_{t-\frac{t-s}{16}}^t \|\widehat{g_{\bar{u}}}\|_{L^1}\Big[\int_s^{s+\frac{t-s}{4}} \frac{C\|\widehat{g_{\bar{u}}}\|_{L^1} \|v\|_{L^{p'}}}{|t-2\tau+2\tau'-s|^{N\ip}}d\tau'+\int_{s+\frac{t-s}{4}}^\tau \frac{C\|g_{\bar{u}}\|_{L^1}}{|t-2\tau+\tau'|^{N\ip}} \frac{\|v\|_{L^{p'}}}{|\tau'-s|^{N\ip}}d\tau' \nonumber \\
&\quad\quad\quad\quad\quad\quad\quad\quad+\int_{t-\frac{t-s}{4}}^{\tau} \frac{C}{|t-2\tau+2\tau'-s|^{\frac{3}{2}}}\|\widehat{g_{\bar{u}}}\|_{L^1} \|v\|_{L^1}d\tau'\Big] d\tau \nonumber \\
&\leq \frac{C\|v\|_{L^{p'}}}{|t-s|^{2N\ip-2}} \nonumber \end{align}
$\tilde{L}_1$ and $\tilde{L}_2$ are estimated as in the previous case.

\item For $X(\tau')=-ig_u W(\tau')$ and $-ig_{\bar{u}} \overline{W(\tau')}$ we will change the order of the integration and use Theorem \ref{th:lw} part (iv):
\begin{align}
\|L_1\|_{L^p}&\leq\int_s^t\int_{\tau'}^t \|\eht P_c g_u e^{iH(t-\tau)}\|_{L^p\rightarrow L^p}\|e^{-iH(t-\tau')}P_c\|_{L^{p'}\rightarrow L^p} \|g_u W(\tau')\|_{L^{p'}}d\tau d\tau' \nonumber \\
&\leq\int_s^t\int_{\tau'}^t\|\widehat{g_u}(\tau)\|_{L^1}\frac{C}{|t-\tau'|^{N\ip}}\|\xs g_u(\tau')\|_{L^{\frac{2p}{p-2}}}\|W(\tau')\|_{\Lsn}d\tau d\tau' \nonumber \\
&\leq \sup_{\tau\in[s,t]}\|\widehat{g_u}(\tau)\|_{L^1}\sup_{\tau'\in[s,t]}\|\xs g_u(\tau')\|_{L^{\frac{2p}{p-2}}}\int_s^t \frac{C}{|t-\tau'|^{N\ip-1}}\frac{C\|v\|_{L^{p'}}}{|\tau'-s|^{N\ip-1}}d\tau' \nonumber \\
&\leq \frac{C\|v\|_{L^{p'}}}{|t-s|^{2N\ip-3}} \nonumber \end{align}
\begin{align}
\|L_2\|_{L^p}&\leq\int_s^t\int_{\tau'}^{t-\frac{t-\tau'}{4}} \|\eht P_c\|_{L^{p'}\rightarrow L^p} \|g_u e^{-iH(\tau'-\tau)} P_c g_u W(\tau')d\tau'\|_{L^{p'}}d\tau d\tau' \nonumber \\
&+\int_s^t\int_{t-\frac{t-\tau'}{4}}^t \|e^{-iH(t-\tau)}P_c g_u e^{iH(t-\tau)}\|_{L^p\rightarrow L^p}\|e^{-iH(t+\tau'-2\tau)} P_c g_u W(\tau')\|_{L^p}d\tau d\tau' \nonumber \\
&\leq \int_s^t\int_{\tau'}^{t-\frac{t-\tau'}{4}} \frac{C}{|t-\tau|^{N\ip}}\|g_u\|_{L^{\frac{2p}{p-2}}}\|e^{-iH(\tau-\tau')}\|_{L^2\rightarrow L^2} \|g_u W(\tau')\|_{L^2}d\tau d\tau' \nonumber \\
&+\int_s^t\int_{t-\frac{t-\tau'}{4}}^t \|\widehat{g_u}(\tau)\|_{L^1} \frac{C}{|t+\tau'-2\tau|^{N\ip}} \|g_u W(\tau')\|_{L^2}d\tau d\tau' \nonumber \\
&\leq \sup_{\tau\in [s,t]}\|g_u(\tau)\|_{L^{\frac{2p}{p-2}}}\int_s^t  \frac{C}{|t-\tau'|^{N\ip-1}} \|\xs g_u(\tau')\|_{L^{\infty}}\| W(\tau')\|_{\Lsn}d\tau' \nonumber \\
&+\sup_{\tau\in [s,t]}\|\widehat{g_u(\tau)}\|_{L^1} \int_s^t \frac{C}{|t-\tau'|^{N\ip-1}}\|\xs g_u(\tau')\|_{L^{\infty}}\| W(\tau')\|_{\Lsn} d\tau'  \nonumber \\
&\leq \frac{C\|v\|_{L^{p'}}}{|t-s|^{2N\ip-3}} \nonumber \end{align}
While $\tilde L_{1,2}$ are estimated similarly with the previous cases:
\begin{align}
\|\tilde{L}_{1,2}\|_{L^p}&\leq\int_s^t \|\eht P_c g_u e^{iH(t-\tau)}\|_{L^p\rightarrow L^p}\int_s^\tau C \|Dh\|_{\C\mapsto H^2} |\<\po,g_u W(\tau')\>|d\tau' d\tau \nonumber \\
&\leq C\int_s^t \|\widehat{g_u}(\tau)\|_{L^1}\int_s^\tau \|Dh\|_{\C\mapsto H^2}\|\po\|_{L^2}\|\xs g_u\|_{L^\infty}\|W(\tau')\|_{\Lsn} d\tau' d\tau \nonumber \\
&\leq C\int_s^t \|\widehat{g_u}(\tau)\|_{L^1}\int_s^\tau \|Dh\|_{\C\mapsto H^2}\|\po\|_{L^2}\|\xs g_u\|_{L^\infty}\frac{C}{|\tau'-s|^{N\ip-1}}d\tau' d\tau \nonumber \\
&\leq C|t-s|^{3-N\ip}\|v\|_{L^{p'}}\nonumber   \end{align}\end{itemize}

The short time estimates for $|t-s|\leq 1$ are now done. For the long time estimates we will assume $t>s+1$ since the case $t<s-1$ can be treated in the same manner. Note that part (i) already gives the required estimate in $L^p(\R^N),\  p=2.$ It remains to obtain the $L^p(\R^N),\ p=q_2$ estimate since the ones for $2<p<q_2$ will be a consequence of Riesz-Thorin interpolation. In what follows it is essential that the kernel $\|e^{-iHt}P_c\|_{L^{p'}\mapsto L^p}$ is integrable in time on $t>s+1,$ see \eqref{eq:lp} for $p=q_2>2N/(N-2).$ This will allow us to use the well known convolution estimate:
 \begin{equation}\label{est:conv}
 \int_{s+\delta}^{t-\delta}|t-\tau|^{-a}|\tau-s|^{-b}d\tau\leq C(\delta,a,b) |t-s|^{-\min\{a,b\}}\qquad \textrm{for }0<\delta\leq 1,\ a,b>0
 \end{equation}
provided $a>1$ or $b>1.$ Note that if $a,b<1$ one has to replace $\min\{a,b\}$ above with $a+b-1< \min\{a,b\}$ which is not sufficient for our linear estimates nor for closing the nonlinear estimates in Section \ref{se:main}.

As before we start with  the long time behavior of $f(t),\ \tilde{f}(t),$ see \eqref{def:f}-\eqref{def:tildef}, and separate them into three integrals:
 \begin{eqnarray}
 f(t)&=&\underbrace{\int_s^{s+1/4}}_{I_1}+\underbrace{\int_{s+1/4}^{t-1/4}}_{I_2}+\underbrace{\int_{t-1/4}^t -i\eht P_c\Dg[\ehts P_c v]d\tau}_{I_3}\nonumber\\
 \tilde{f}(t)&=&\underbrace{\int_s^{s+1/2}}_{II_1}+\underbrace{\int_{s+1/2}^{t-1/2}}_{II_2}+\underbrace{\int_{t-1/2}^t\eht P_c Dh|_{a(\tau)}i\<\po,\Dg[\ehts P_c v]\>d\tau}_{II_3}.\nonumber
 \end{eqnarray}
For the middle integrals, $I_2$ and $II_2,$ we simply use \eqref{eq:lp} combined with \eqref{est:conv}. For the remaining integrals we remove the singularities at $\tau=s$ respectively at $\tau=t$ as in the above short time estimates. More precisely, for $p=q_2>2N/(N-2),$ we have:
\begin{align}
\|I_2\|_{L^p}&\leq\int_{s+1/4}^{t-1/4} \|\eht P_c\|_{L^{p'}\mapsto L^p} (\|g_u\ehts P_c v\|_{L^{p'}}+\|g_{\bar{u}}e^{iH(\tau-s)} P_c \bar{v}\|_{L^{p'}})\nonumber \\
&\leq\int_{s+1/4}^{t-1/4} \frac{C_p}{|t-\tau|^{N\ip}} ( \|g_u\|_{L^{\frac{p}{p-2}}}+\|g_{\bar{u}}\|_{L^{\frac{p}{p-2}}})\|\ehts P_c v\|_{L^p}d\tau \nonumber \\
&\leq\sup_{\tau\in\R}\{ \|g_u(\tau)\|_{L^{\frac{p}{p-2}}}+\|g_{\bar{u}}(\tau)\|_{L^{\frac{p}{p-2}}}\}\int_{s+1/4}^{t-1/4} \frac{C_p}{|t-\tau|^{N\ip}}\frac{C_p\|v\|_{L^{p'}}}{|\tau-s|^{N\ip}}d\tau \nonumber \\
&\leq C\frac{\|v\|_{L^{p'}}}{|t-s|^{N\ip}} \nonumber  \end{align}
 \begin{align}
 \|II_2\|_{L^p}&\leq\int_{s+1/2}^{t-1/2}\|\eht P_c\|_{L^{p'}\mapsto L^p} \|Dh|_{a(\tau)}\|_{L^{p'}}\|\po\|_{L^{2}}(\|g_u\|_{L^{\frac{2p}{p-2}}}+\|g_{\bar{u}}\|_{L^{\frac{2p}{p-2}}})\|\ehts P_c v\|_{L^p}d\tau \nonumber \\
 &\leq\|\po\|_{L^2}\sup_{\tau\in\R}\{\|Dh|_{a(\tau)}\|_{L^{p'}}(\|g_u(\tau)\|_{L^{\frac{2p}{p-2}}}+\|g_{\bar{u}}(\tau)\|_{L^{\frac{2p}{p-2}}})\}
 \int_{s+1/2}^{t-1/2} \frac{C_p}{|t-\tau|^{N\ip}}\frac{C_p\|v\|_{L^{p'}}}{|\tau-s|^{N\ip}}d\tau \nonumber \\
 &\leq\frac{C\|v\|_{L^{p'}}}{|t-s|^{N\ip}}\nonumber  \end{align}

\begin{align}
\|I_1\|_{L^p}&\leq\int_s^{s+1/4} \|\ehs P_c\|_{L^{p'}\mapsto L^p}\| e^{iH(\tau-s)} P_c g_u\ehts P_c v\|_{L^{p'}}d\tau \nonumber \\
&\quad+\int_s^{s+1/4} \|e^{-iH(t+s-2\tau)}\|_{L^{p'}\mapsto L^p}\| \ehts P_c g_{\bar{u}} e^{iH(\tau-s)} P_c \bar{v}\|_{L^{p'}}d\tau \nonumber \\
&\leq\int_s^{s+1/4}\frac{C}{|t-s|^{N\ip}}C\|\widehat{g_u}\|_{L^1}\|v\|_{L^{p'}}d\tau +\int_s^{s+1/4}\frac{C_p}{|t+s-2\tau|^{N\ip}}C\|\widehat{g_{\bar{u}}}\|_{L^1}\|v\|_{L^{p'}}d\tau \nonumber \\
&\leq C\|v\|_{L^{p'}}\sup_{\tau\in\R}\{\|\widehat{g_u}(\tau)\|_{L^1}+\|\widehat{g_{\bar{u}}}(\tau)\|_{L^1}\}\int_s^{s+1/4}  \frac{C_p}{|t-s|^{N\ip}}+\frac{C_p}{|t-s-1/2|^{N\ip}}d\tau \nonumber \\
&\leq C\frac{\|v\|_{L^{p'}}}{|t-s|^{N\ip}} \nonumber \end{align}
 \begin{align}
\|I_3\|_{L^p}&\leq\int_{t-1/4}^t \| e^{-iH(t-\tau)} P_c g_ue^{iH(t-\tau)} \|_{L^{p}\mapsto L^p}\|\ehs P_c v\|_{L^p}d\tau \nonumber \\
&\quad+\int_{t-1/4}^t \| e^{-iH(t-\tau)} P_c g_{\bar{u}} e^{iH(t-\tau)} \|_{L^{p}\mapsto L^p} \|e^{-iH(t+s-2\tau)}P_c \bar{v}\|_{L^{p}}d\tau \nonumber \\
&\leq\int_{t-1/4}^tC\|\widehat{g_u}\|_{L^1}\frac{C_p\|v\|_{L^{p'}}}{|t-s|^{N\ip}}d\tau +\int_{t-1/4}^tC\|\widehat{g_{\bar{u}}}\|_{L^1}\frac{C_p\|v\|_{L^{p'}}}{|t+s-2\tau|^{N\ip}}d\tau \nonumber \\
&\leq C\|v\|_{L^{p'}}\sup_{\tau\in\R}\{\|\widehat{g_u}(\tau)\|_{L^1}+\|\widehat{g_{\bar{u}}}(\tau)\|_{L^1}\}\int_{t-1/4}^t  \frac{C_p}{|t-s|^{N\ip}}+\frac{C_p}{|t-s-1/2|^{N\ip}}d\tau \nonumber \\
&\leq C\frac{\|v\|_{L^{p'}}}{|t-s|^{N\ip}} \nonumber \end{align}
 \begin{align}
 \|II_1\|_{L^p}&\leq\int_s^{s+1/2} \|\eht P_c\|_{L^{p'}\mapsto L^p} \|Dh|_{a(\tau)}\|_{L^{p'}} |\<\po,\Dg\ehts P_c v\>|\nonumber \\
 &\leq\int_s^{s+1/2} \frac{C_p}{|t-\tau|^{N\ip}}\|Dh|_{a(\tau)}\|_{L^{p'}}\Big[\< e^{iH(\tau-s)}\po,e^{iH(\tau-s)}g_u\ehts P_c v\> \nonumber \\
 &\quad\quad\quad\quad\quad \quad+\<\ehts \po,\ehts g_{\bar{u}}e^{iH(\tau-s)} P_c \bar{v}\>\Big]d\tau \nonumber \\
 &\leq \frac{C_p}{|t-s-1/2|^{N\ip}}\int_s^{s+1/2}\|Dh|_{a(\tau)}\|_{L^{p'}}\Big[\|e^{iH(\tau-s)}\po\|_{L^p}\|e^{iH(\tau-s)}g_u\ehts P_c v\|_{L^{p'}}\nonumber \\
 &\quad\quad\quad\quad\quad\quad\quad\quad\qquad\qquad +\|\ehts\po\|_{L^p}\|\ehts g_{\bar{u}} e^{iH(\tau-s)} P_c \bar{v}\|_{L^{p'}}\Big]d\tau \nonumber \\
 &\leq\frac{C\|\po\|_{L^p}}{|t-s|^{N\ip}}\int_s^{s+1/2}\|Dh|_{a(\tau)}\|_{L^{p'}}(\|\widehat{g_u}\|_{L^1}+\|\widehat{g_{\bar{u}}}\|_{L^1}) \|v\|_{L^{p'}} d\tau \nonumber \\
 &\leq\frac{C\|\po\|_{H^2}\|v\|_{L^{p'}}}{|t-s|^{N\ip}}\sup_{\tau\in\R}\{\|Dh|_{a(\tau)}\|_{L^{p'}}(\|\widehat{g_u}(\tau)\|_{L^1}+\|\widehat{g_{\bar{u}}}(\tau)\|_{L^1})\}\leq\frac{C\|v\|_{L^{p'}}}{|t-s|^{N\ip}} \nonumber \end{align}
 \begin{align}
 \|II_3\|_{L^p}&\leq\int_{t-1/2}^t \|\eht P_c\|_{H^2\mapsto H^2}\|Dh|_{a(\tau)}\|_{\C\rightarrow H^2}|\<\po,\Dg\ehts P_c v]\>|d\tau \nonumber \\
 &\leq\int_{t-1/2}^t C\|Dh|_{a(\tau)}\|_{\C\rightarrow H^2}\|\po\|_{L^2}(\|g_u\|_{L^{\beta}}+\|g_{\bar{u}}\|_{L^{\beta}})\|\ehts P_c v\|_{L^{p}}d\tau \nonumber \\
 &\leq\|\po\|_{L^2}\sup_{\tau\in\R}\{\|Dh|_{a(\tau)}\|_{\C\rightarrow H^2}(\|g_u(\tau)\|_{L^{\beta}}+\|g_{\bar{u}}(\tau)\|_{L^{\beta}})\}\int_{t-1/2}^t \frac{C_p\|v\|_{L^{p'}}}{|\tau-s|^{N\ip}} d\tau \nonumber \\
 &\leq \frac{C\|v\|_{L^{p'}}}{|t-s|^{N\ip}} \nonumber \end{align}

Similarly we will investigate the long time behavior, $t>s+1,$ of
$L(s)W.$ We split it into three integrals with
$s_1=\min\{s+1,t-1/16\}$:
 $$[L(s)W](t)=\underbrace{\int_s^{s_1}}_{L_3}
 +\underbrace{\int_{s_1}^{t-\frac{1}{16}}}_{L_4}
 +\underbrace{\int_{t-\frac{1}{16}}^t\eht P_c \left( -iDg|_{\psi_{E(\tau)}}[W(\tau)]+Dh|_{a(\tau)}i\<\psi_0,Dg|_{\psi_{E(\tau)}}[W(\tau)]\>\right) d\tau}_{L_5}$$
Due to Theorem \ref{th:lw} part (iv) and $W(\tau)=T(\tau,s)v,$ the
integral $L_3$ has an integrable singularity at $\tau=s$ while $L_4$
has no singularities. A combination of \eqref{eq:lp}, estimates in
Theorem \ref{th:lw} part (iv), and \eqref{est:conv} gives the
required result for $L_3$ and $L_4.$ In $L_5$ we will first remove
the singularity at $\tau=t$ in a similar manner we did it for short
time estimates. More precisely, for $p=q_2,\ 2N/(N-2)<q_2<2N/(N-4),$
we have:
 \begin{eqnarray}
 \|L_3\|_{L^p}&\leq &\int_s^{s_1}\|\eht P_c\|_{L^{p'}\mapsto
 L^p}(\|\xs g_u\|_{L^\beta}+\|\xs g_{\bar
 u}\|_{L^\beta})\|W(\tau)\|_{\Lsn}d\tau\nonumber\\
 &&+\int_s^{s_1}\|\eht P_c\|_{L^{p'}\mapsto
 L^p}\|Dh|_{a(\tau)}\|_{\C\mapsto L^{p'}}|\<\psi_0,Dg|_{\psi_{E(\tau)}}[W(\tau)]\>|
 d\tau\nonumber\\
 &\leq & \sup_{\tau\in\R}\{(1+\|Dh|_{a(\tau)}\|_{\C\mapsto
 L^{p'}}\|\po\|_{L^2})(\|\xs g_u(\tau)\|_{L^\beta}+\|\xs g_{\bar
 u}(\tau)\|_{L^\beta})\}\int_s^{s_1}\frac{C_p\|W(\tau)\|_{\Lsn}d\tau}{|t-\tau|^{N\ip}}\nonumber\\
 &\leq &
 \frac{C}{|t-s-1/16|^{N\ip}}\int_s^{s_1}\frac{C_p\|v\|_{L^{p'}}}{|\tau-s|^{N\ip
 -1}}d\tau\leq \frac{C\|v\|_{L^{p'}}}{|t-s|^{N\ip}}\nonumber
 \end{eqnarray}
We estimate $L_4$ exactly as $L_3$ but now $\|W(\tau)\|_{\Lsn}\leq
C_p(1+|\tau-s|)^{-N\ip},$ see Theorem \ref{th:lw} part (iv), and the
convolution estimate \eqref{est:conv} is now employed to yield
exactly the same result. For $L_5,$ one of the integrands has no
singularities:
\begin{eqnarray}
 \lefteqn{\|\int_{t-\frac{1}{16}}^t \eht
 P_cDh|_{a(\tau)}i\<\po,\Dg[W(\tau)]\>
 d\tau\|_{L^p}}\nonumber\\
 &\leq &\int_{t-\frac{1}{16}}^t C\|\eht P_c\|_{H^2\mapsto H^2}\|Dh|_{a(\tau)}\|_{\C\mapsto H^2}|\<\po,\Dg[W(\tau)]\>|
 d\tau\nonumber\\
 &\leq & \int_{t-\frac{1}{16}}^t C\|Dh|_{a(\tau)}\|_{\C\mapsto H^2}\|\po\|_{L^2}(\|\xs g_u\|_{L^\infty}+\|\xs g_{\bar
 u}\|_{L^\infty})\|W(\tau)\|_{\Lsn}d\tau\nonumber\\
 &\leq &\|\po\|_{L^2}\sup_{\tau\in\R}\{\|Dh|_{a(\tau)}\|_{\C\mapsto H^2}(\|\xs g_u\|_{L^\infty}+\|\xs g_{\bar
 u}\|_{L^\infty})\}\int_{t-\frac{1}{16}}^t\frac{C_p\|v\|_{L^{p'}}}{|\tau-s|^{N\ip}}d\tau\leq \frac{C\|v\|_{L^{p'}}}{|t-s|^{N\ip
}}.\nonumber
 \end{eqnarray} However to remove the non-integrable singularity at $\tau=t$ in the
remaining integral we need to plug in (\ref{sow}) in it:
\begin{eqnarray}
\lefteqn{\int_{t-\frac{1}{16}}^t \eht P_c \Dg[W(\tau)] d\tau}\nonumber \\
&=&\int_{t-\frac{1}{16}}^t \eht P_c g_u\Big(\int_s^\tau e^{-iH(\tau-\tau')} P_c (-i\Dg[e^{-iH(\tau'-s)}P_c v]+Dhi\<\po,\Dg[e^{-iH(\tau'-s)}P_c v]\>)d\tau'\nonumber \\
&&\quad\quad\quad\quad\quad\quad\quad\quad\quad\quad\quad\quad+\int_s^\tau e^{-iH(\tau-\tau')} P_c (-i\Dg[W(\tau')]+Dhi\<\po,\Dg[W(\tau')]\>)d\tau'\Big)d\tau \nonumber \\
&&+\int_{t-\frac{1}{16}}^t \eht P_c g_{\bar{u}}\Big(\int_s^\tau e^{iH(\tau-\tau')}P_c(\overline{-i \Dg[e^{-iH(\tau'-s)}P_c v]}+\overline{Dhi\<\po,\Dg[e^{-iH(\tau'-s)}P_c v]\>})d\tau' \nonumber \\
&&\quad\quad\quad\quad\quad\quad\quad\quad\quad\quad\quad\quad+\int_s^\tau
e^{iH(\tau-\tau')} P_c(\overline{-i
\Dg[W(\tau')]}+\overline{Dhi\<\po,\Dg[W(\tau')]\>})d\tau'\Big)d\tau
\nonumber \end{eqnarray} We will add $e^{iH(t-\tau)}$ and $\eht$
terms after $g_u$ and $g_{\bar{u}}$. Then  all the terms will be
similar to $L_1,\ L_2,\ \tilde L_1,$ and $\tilde L_2,$ see
$(\ref{l1})-(\ref{lt2}).$  After separating the the inside integrals
into pieces, we will estimate short time step integrals exactly the
same way we did for short time behavior by using estimate
\eqref{est:fm}, and the other integrals will be estimated using the
usual norms. For completeness we show below how each term is
treated:
\begin{itemize}
 \item For $X(\tau')=-ig_u e^{-iH(\tau'-s)}P_c v$ we have
\begin{align}
\|L_1\|_{L^p}&\leq\int_{t-\frac{1}{16}}^t \|\eht P_c g_u e^{iH(t-\tau)}\|_{L^p\rightarrow L^p}\nonumber \\
&\quad\quad\Big[\int_s^{s+\frac{1}{4}}\|e^{-iH(t-s)}P_c\|_{L^{p'}\rightarrow L^p} \|e^{iH(\tau'-s)}P_c g_u e^{-iH(\tau'-s)}P_c v\|_{L^{p'}}d\tau' \nonumber \\
&\quad\quad+\int_{s+\frac{1}{4}}^{t-\frac{1}{4}} \|e^{-iH(t-\tau')} P_c\|_{L^{p'}\rightarrow L^p}\|g_u e^{-iH(\tau'-s)}P_c v\|_{L^p}d\tau' \nonumber \\
&\quad\quad+\int_{t-\frac{1}{4}}^{\tau}\|e^{-iH(t-\tau')}P_c g_u e^{iH(t-\tau')}\|_{L^p\rightarrow L^p} \|e^{-iH(t-s)}P_c v\|_{L^{p}}d\tau'\Big]d\tau \nonumber \\
&\leq \frac{C\|v\|_{L^{p'}}}{|t-s|^{N\ip}} \nonumber \end{align}

\begin{align}
\|L_2\|_{L^p}&\leq\int_{t-\frac{1}{16}}^t \|\eht P_c g_{\bar{u}} e^{iH(t-\tau)}\|_{L^p\rightarrow L^p} \nonumber \\
&\quad\quad\quad\quad\times\Big[ \int_s^{s+\frac{1}{4}}\|e^{-iH(t+s-2\tau)}P_c\|_{L^{p'}\rightarrow L^p}\|e^{-iH(\tau'-s)}P_c g_{\bar{u}} e^{iH(\tau'-s)}P_c \bar{v}\|_{L^{p'}}d\tau' \nonumber \\
&\quad\quad\quad+\int_{s+\frac{1}{4}}^{t-\frac{1}{4}}\|e^{-iH(t-2\tau+\tau')}P_c\|_{L^{p'}\rightarrow L^p}\|g_{\bar{u}} e^{iH(\tau'-s)}P_c \bar{v}\|_{L^{p'}}d\tau' \nonumber \\
&\quad\quad\quad+\int_{t-\frac{1}{4}}^{\tau}\|e^{-iH(t+\tau'-2\tau)}P_c g_{\bar{u}}e^{iH(t+\tau'-2\tau)}\|_{L^{p'}\rightarrow L^p}\|e^{iH(t+s-2\tau)}P_c \bar{v}\|_{L^{p'}}d\tau'\Big]d\tau \nonumber \\
&\leq\frac{C\|v\|_{L^{p'}}}{|t-s|^{N\ip}} \nonumber \end{align}

\begin{align}
\|\tilde{L}_1\|_{L^p}&\leq\int_{t-\frac{1}{16}}^t \|\eht P_c g_u e^{iH(t-\tau)}\|_{L^p\rightarrow L^p}\nonumber \\
&\quad\quad\Big[\int_s^{s+\frac{1}{4}}\|e^{-iH(t-\tau')}P_c\|_{L^{p'}\mapsto L^p}\|Dh\|_{\C\mapsto L^{p'}}|\<e^{iH(\tau'-s)}\po, e^{iH(\tau'-s)} g_u e^{-iH(\tau'-s)}P_c v\>|d\tau' \nonumber \\
&\quad\quad+\int_{s+\frac{1}{4}}^{t-\frac{1}{4}} \|e^{-iH(t-\tau')}
P_c\|_{L^{p'}\mapsto L^p}\|Dh\|_{\C\mapsto L^{p'}}|\<\po,g_u
e^{-iH(\tau'-s)}P_c v\>|d\tau'\nonumber\\
&\quad\quad+\int_{t-\frac{1}{4}}^{\tau} \|e^{-iH(t-\tau')} P_c\|_{H^2\mapsto H^2}\|Dh\|_{\C\mapsto H^2}|\<\po,g_u e^{-iH(\tau'-s)}P_c v\>|d\tau'\Big]d\tau \nonumber \\
&\leq\int_{t-\frac{1}{16}}^t\|\widehat{g_u}\|_{L^1}\Big[\int_s^{s+\frac{1}{4}}\frac{C_p\|Dh\|}{|t-\tau'|^{N\ip}}\|\po\|_{L^p}\|e^{iH(\tau'-s)}P_c g_u e^{-iH(\tau'-s)}P_c v\|_{L^{p'}}d\tau' \nonumber \\
&\quad\quad\quad\quad\quad\qquad+\int_{s+\frac{1}{4}}^{t-\frac{1}{4}}\frac{C_p\|Dh\|}{|t-\tau'|^{N\ip}}\|\po\|_{L^2}\|g_u\|_{L^\beta}\| e^{-iH(\tau'-s)}P_c v\|_{L^{p}}d\tau' \nonumber \\
&\quad\quad\quad\quad\quad\qquad+\int_{t-\frac{1}{4}}^\tau
C\|Dh\|\|\po\|_{L^2}\|g_u\|_{L^\beta}\| e^{-iH(\tau'-s)}P_c
v\|_{L^p}d\tau'\Big]d\tau
\leq\frac{C\|v\|_{L^{p'}}}{|t-s|^{N\ip}}\nonumber
\end{align} $\tilde{L}_2$ is treated exactly the same as $\tilde{L}_1$ except that in the decomposition of the inside integral $4\tau-3t$ is used instead of $t-1/4/.$

\item For $X(\tau')=-ig_{\bar{u}} e^{iH(\tau'-s)}P_c \bar{v}$ we have
\begin{align}
\|L_1\|_{L^p}\leq&\int_{t-\frac{1}{16}}^t \|\eht P_c g_u e^{iH(t-\tau)}\|_{L^p\rightarrow L^p}\nonumber \\ &\times\Big[\int_{s}^{s+\frac{1}{4}} \|e^{-iH(t+s-2\tau')} P_c\|_{L^{p'}\rightarrow L^p}\|e^{-iH(\tau'-s)}P_c g_{\bar{u}} e^{iH(\tau'-s)}P_c \bar{v}\|_{L^{p'}}d\tau'\nonumber \\
&\quad\quad\quad+\int_{s+\frac{1}{4}}^{t-\frac{1}{4}}\|e^{-iH(t-\tau')}\|_{L^{p'}\rightarrow L^p}\|g_{\bar{u}} e^{iH(\tau'-s)}P_c \bar{v}\|_{L^{p'}}d\tau' \nonumber \\
&\quad\quad\quad+\int_{t-\frac{1}{4}}^{\tau}\|e^{-iH(t-\tau')} P_c g_{\bar{u}}e^{iH(t-\tau')}\|_{L^{p'}\rightarrow L^p}\|e^{-iH(t+s-2\tau')}P_c \bar{v}\|_{L^{p'}}d\tau'\Big]d\tau \nonumber \\
&\leq \frac{C\|v\|_{L^{p'}}}{|t-s|^{N\ip}} \nonumber\end{align}
\begin{align}
\|L_2\|_{L^p}&\leq\int_{t-\frac{1}{16}}^t \|\eht P_c g_{\bar{u}} e^{iH(t-\tau)}\|_{L^p\rightarrow L^p}\Big[\int_{s}^{s+\frac{1}{4}}\| e^{-iH(t-2\tau+2\tau'-s)}e^{iH(\tau'-s)} P_c g_u e^{-iH(\tau'-s)}P_c v\|_{L^p}d\tau' \nonumber \\
&\quad\quad\quad+\int_{s+\frac{1}{4}}^{t-\frac{1}{4}} \|e^{-iH(t-2\tau+\tau')} P_c\|_{L^{p'}\rightarrow L^{p'}} \|g_u e^{-iH(\tau'-s)}P_c v\|_{L^{p'}}d\tau'  \nonumber \\
&\quad\quad\quad+\int_{t-\frac{1}{4}}^\tau\| e^{-iH(t+\tau'-2\tau)}P_c g_u e^{iH(t+\tau'-2\tau)} e^{-iH(t-2\tau+2\tau'-s)}P_c v\|_{L^p}d\tau'\Big]d\tau \nonumber \\
&\leq \frac{C\|v\|_{L^{p'}}}{|t-s|^{N\ip}} \nonumber \end{align}
$\tilde L_1$ and $\tilde L_2$ are treated as in the previous case.
\item For $X(\tau')=-ig_u W(\tau')$ and $X(\tau')=-ig_{\bar{u}} \overline{W(\tau')}$ we will separate the $L_1$ term into three
integrals. For the first integral we will use short time $\Lsn$
estimate for $W$. Also note that one can obtain the same estimates
for $|t-s|\leq\frac{1}{4}$ and $|t-s|>\frac{1}{4}$ in the Theorem
\ref{th:lw} part $(iv).$  For the last integral we will change the
order of the integration:
\begin{align}
\|L_1\|_{L^p}&\leq\int_{t-\frac{1}{16}}^t \|\eht P_c g_u(\tau) e^{iH(t-\tau)}\|_{L^p\rightarrow L^p}\int_s^{\tau}\|e^{-iH(t-\tau')}P_c\|_{L^{p'}\rightarrow L^p} \|g_u W(\tau')\|_{L^{p'}}d\tau'd\tau \nonumber \\
&\leq\sup_{\tau\in\R}\{\|\widehat{g_u}(\tau)\|_{L^1}\}\int_{t-\frac{1}{16}}^t \int_s^{\tau}\frac{C_p}{|t-\tau'|^{N\ip}}\|\xs g_u\|_{L^\beta}\|W(\tau')\|_{\Lsn}d\tau' d\tau \nonumber \\
&\leq\int_{t-\frac{1}{16}}^t\Bigg[\int_s^{s+\frac{1}{4}} \frac{C}{|t-\tau'|^{N\ip}}\frac{C_p\|v\|_{L^{p'}}}{|\tau'-s|^{N\ip-1}}d\tau'+\int_{s+\frac{1}{4}}^{t-\frac{1}{16}}\frac{C}{|t-\tau'|^{N\ip}}\frac{C_p\|v\|_{L^{p'}}}{(1+|\tau'-s|)^{N\ip}}d\tau'\Bigg] d\tau \nonumber \\
&+\int_{t-\frac{1}{16}}^t\int_{\tau'}^t\frac{C}{|t-\tau'|^{N\ip}}\frac{C_p\|v\|_{L^{p'}}}{(1+|\tau'-s|)^{N\ip}}d\tau d\tau'\nonumber \\
&\leq \frac{C\|v\|_{L^{p'}}}{|t-s|^{N\ip}} \nonumber
\end{align}
Similar to $L_1$ we will split $L_2$ in three integrals. In the
first and last we use estimate \eqref{est:fm} and we also change the
order of integration in the last integral:
\begin{align}
\|L_2\|_{L^p}&\leq\int_{t-\frac{1}{16}}^t \|\eht P_c g_{\bar{u}} e^{iH(t-\tau)}\|_{L^p\rightarrow L^p}\int_s^{t-\frac{1}{4}}\|e^{-iH(t+\tau'-2\tau)}P_c\|_{L^{p'}\rightarrow L^p} \|g_u W(\tau')\|_{L^{p'}}d\tau'd\tau \nonumber \\
&+\int_{t-\frac{1}{16}}^t\int_{4\tau-3t}^\tau \|\eht P_c \|_{L^{p'}\rightarrow L^p}\|g_{\bar{u}}\|_{L^\beta}\|e^{iH(\tau-\tau')}P_c\|_{L^2\rightarrow L^2}\|\xs g_u\|_{L^\infty}\|W(\tau')\|_{\Lsn}d\tau' d\tau \nonumber \\
&+\int_{t-\frac{1}{4}}^t\int_{t-\frac{t-\tau'}{4}}^t \|\eht P_c g_{\bar{u}} e^{iH(t-\tau)}\|_{L^p\rightarrow L^p}\|e^{-iH(t+\tau'-2\tau)}P_c\|_{L^{p'}\rightarrow L^p} \|g_u W(\tau')\|_{L^{p'}}d\tau d\tau' \nonumber \\
&\leq\int_{t-\frac{1}{16}}^t \int_s^{s+\frac{1}{4}}\frac{C\|\widehat{g_{\bar{u}}}\|_{L^1}}{|t+\tau'-2\tau|^{N\ip}}\frac{\|\xs g_u\|_{L^\beta}\|v\|_{L^{p'}}}{|\tau'-s|^{N\ip-1}}d\tau' d\tau  \nonumber \\
&+\int_{t-\frac{1}{16}}^t
 \int_{s+\frac{1}{4}}^{t-\frac{1}{4}}\frac{C\|\widehat{g_{\bar{u}}}\|_{L^1}}{|t+\tau'-2\tau|^{N\ip}}\frac{\|\xs g_u\|_{L^{\beta}}\|v\|_{L^{p'}}}{(1+|\tau'-s|)^{N\ip}}
                                               d\tau' d\tau  \nonumber \\
&+\int_{t-\frac{1}{16}}^t\int_{4\tau-3t}^\tau \frac{C}{|t-\tau|^{N\ip}}\frac{\|v\|_{L^{p'}}}{(1+|\tau'-s|)^{N\ip}}d\tau' d\tau\nonumber \\
&+\int_{t-\frac{1}{4}}^t\int_{t-\frac{t-\tau'}{4}}^t \frac{C}{|t+\tau'-2\tau|^{N\ip}}\frac{\|v\|_{L^{p'}}}{(1+|\tau'-s|)^{N\ip}}d\tau d\tau' \nonumber \\
&\leq \frac{C\|v\|_{L^{p'}}}{|t-s|^{N\ip}} \nonumber
\end{align} $\tilde{L}_1$ and $\tilde{L}_2$ terms are estimated as
in the previous cases, more precisely:
\begin{align}
\|\tilde{L}_2\|_{L^p}&\leq\int_{t-\frac{1}{16}}^t \|\eht P_c g_u e^{iH(t-\tau)}\|_{L^p\rightarrow L^p}\nonumber \\
&\quad\quad\Big[\int_s^{4\tau-3t}\|e^{-iH(t-2\tau+\tau')}P_c\|_{L^{p'}\mapsto L^p}\|Dh\|_{\C\mapsto L^{p'}}|\<\po, g_u W(\tau')\>|d\tau' \nonumber \\
&\quad\quad+\int_{4\tau-3t}^{\tau} \|e^{-iH(t-2\tau+\tau')} P_c\|_{H^2\mapsto H^2}\|Dh\|_{\C\mapsto H^2}||\<\po, g_u W(\tau')\>|d\tau'\Big]d\tau \nonumber \\
&\leq\frac{C\|v\|_{L^{p'}}}{|t-s|^{N\ip}}\nonumber
\end{align} where we used $$|\<\po, g_u W(\tau')\>|\leq
\|\po\|_{L^2}\sup_{\tau'\in\R}\{\|\xs
g_u(\tau')\|_{L^\infty}\}\|W(\tau')\|_{\Lsn}\leq\left\{\begin{array}{ll}
\frac{C\|v\|_{L^{p'}}}{|\tau'-s|^{N\ip-1}} & \textrm{if }
|\tau'-s|\leq 1\\ \frac{C\|v\|_{L^{p'}}}{(1+|\tau'-s|)^{N\ip}} &
\textrm{if } |\tau'-s|>1\end{array}\right.$$
%\begin{align}
%\|\tilde{L}_1\|_{L^p}&\leq\int_{t-\frac{1}{4}}^t \|\eht P_c g_u e^{iH(t-\tau)}\|_{L^p\rightarrow L^p} \nonumber \\
%&\quad\quad\quad\times\Bigg[\int_s^{t-\frac{1}{16}}\|e^{-iH(t-\tau')}P_c\|_{L^{p'}\rightarrow L^p}\|Dh\|\|\po\|_{L^2} \|\xs g_u\|_{L^\infty} \|W(\tau')\|_{\Lsn}d\tau'\nonumber \\
%&+\int_{t-\frac{1}{16}}^\tau\|e^{-iH(t-\tau')}P_c\|_{H^2\rightarrow H^2}\|Dh\|\|\po\|_{L^2} \|\xs g_u\|_{L^\infty} \|W(\tau')\|_{\Lsn}d\tau'\Bigg]d\tau \nonumber \\
%&\leq\int_{t-\frac{1}{16}}^t \|\widehat{g_u}\|_{L^1}\int_s^{s+\frac{1}{4}}\frac{C}{|t-\tau'|^{N\ip}}\|Dh\|\|\po\|_{L^2}\|\xs g_u\|_{L^\infty}\frac{C\|v\|_{L^{p'}}}{|\tau'-s)^{N\ip-1}}d\tau' d\tau \nonumber \\
%&\leq\int_{t-\frac{1}{16}}^t\Bigg[\int_{s+\frac{1}{4}}^{t-\frac{1}{16}} \frac{C}{|t-\tau'|^{N\ip}}\frac{C\|v\|_{L^{p'}}}{(1+|\tau'-s|)^{N\ip}}d\tau \nonumber \\
%&\quad\quad\quad\quad+\int_{t-\frac{1}{16}}^\tau \frac{C\|v\|_{L^{p'}}}{(1+|\tau'-s|)^{N\ip}}d\tau\Bigg] d\tau' \nonumber \\
%&\leq \frac{C\|v\|_{L^{p'}}}{|t-s|^{N\ip}} \nonumber
%\end{align}

\end{itemize}

This finishes the proof of $(ii)$.

$(iii)$  The case $p=2$ has already been proven in part (i). It
remains to show the estimate for $p=\frac{2N}{N-2}$ since the ones
for $2<p<\frac{2N}{N-2}$ follow from Riesz-Thorin interpolation. We
will again use the definition \eqref{w} and expansion \eqref{sow}
together with notations \eqref{def:f}-\eqref{def:tildef}, see
\eqref{L} for the definition of $L(s).$ We will treat the $t\geq s$
case as the $t<s$ one can be treated similarly.

For the $f$ term let us first consider $s\leq t\leq 1.$ Recall that
$$\tilde{T}(t,s)=-i\int_s^{\min\{t,s+1\}}e^{-iH(t-\tau)}P_c
g_u(\tau)e^{-iH(\tau-s)}P_c .d\tau$$ Then for this time interval the
forcing term corresponding to $f$ of the operator
$T(t,s)-\tilde{T}(t,s)$ becomes
\begin{equation}
\tilde{I}_1=-i\int_s^t e^{-iH(t-\tau)}P_c g_{\bar{u}} e^{iH(\tau-s)}
P_c \bar{v} d\tau \label{Ishort} \end{equation} For fixed $t$ and
$s$ we have \begin{align}
\|\tilde{I}_1\|_{L^2}&=\|\int_s^t e^{-iH(t-\tau)}P_c g_{\bar{u}} e^{iH(\tau-s)} P_c \bar{v} d\tau\|_{L^2} \nonumber \\
&=\|e^{-iH(t+s)}P_c \int_s^t e^{2iH\tau}P_c \underbrace{e^{-iH(\tau-s)}P_c g_{\bar{u}} e^{iH(\tau-s)} P_c \bar{v}}_{q(\tau)} d\tau\|_{L^2} \nonumber \\
&\leq \|e^{-iH(t+s)}P_c\|_{L^2\rightarrow L^2} \|\int_s^t e^{2iH\tau}P_c q(\tau) d\tau\|_{L^2} \nonumber \\
&\leq C \|q(\tau)\|_{L^{2}([s,s+1], L^{p'})}\leq C\|v\|_{L^{p'}}
\nonumber \end{align} at the last step we used end point Stricharz
estimates \cite[Theorem 1.2]{kt:eps} and the fact that
$\|q(\tau)\|_{L^{p'}}\leq
C\|\widehat{g_{\bar{u}}}(\tau)\|_{L^1}\|v\|_{L^{p'}},$ see
\eqref{est:fm}.

For the long time we split $f$ as follows
$$f=\underbrace{\int_s^{s+1}}_{I_1}+\underbrace{\int_{s+1}^t-ie^{-iH(t-\tau)}P_c \Dg[e^{-iH(\tau-s)} P_c v] d\tau}_{I_2}$$
Then $I_1$ is estimated  exactly as \eqref{Ishort} above and for
$I_2$ we have via Stricharz estimates for the admissible pair
$(\gamma,\rho),$ $\gamma\geq 2$ fixed but $\gamma\not= \infty:$
\begin{align}
\|I_2\|_{L^2}&\leq C\Big(\int_{s+1}^\infty \|g_u \ehts P_c v\|_{L^{\rho'}}^{\gamma'}d\tau\Big)^{\frac{1}{\gamma'}}+C\Big(\int_{s+1}^\infty \|g_{\bar{u}} e^{iH(\tau-s)} P_c \bar{v}\|_{L^{\rho'}}^{\gamma'}d\tau\Big)^{\frac{1}{\gamma'}} \nonumber \\
&\leq C\Big(\int_{s+1}^\infty \|g_u\|_{L^{\beta}}^{\gamma'}\|\ehts P_c v\|_{L^p}^{\gamma'}d\tau\Big)^{\frac{1}{\gamma'}}+C\Big(\int_{s+1}^\infty \|g_{\bar{u}}\|_{L^{\beta}}^{\gamma'}\|e^{iH(\tau-s)} P_c \bar{v}\|_{L^p}^{\gamma'}d\tau\Big)^{\frac{1}{\gamma'}} \nonumber \\
&\leq C\Big(\int_{s+1}^\infty
\frac{d\tau}{|\tau-s|^{N(\frac{1}{2}-\frac{1}{p})\gamma'}}\Big)^{\frac{1}{\gamma'}}<\infty
\nonumber
\end{align}
where $1/\beta+1/p=1/\rho'$ and we used
$N(\frac{1}{2}-\frac{1}{p})\gamma'>1$ for $p=\frac{2N}{N-2}$ and
$\gamma\not= 2$.

Similarly for the other forcing term we
have:$$\tilde{f}(t)=\underbrace{\int_s^{s+1}}_{II_1}+\underbrace{\int_{s+1}^t\eht
P_c Dh|_{a(\tau)}i\<\po,\Dg[\ehts P_c v]\>d\tau.}_{II_2}$$ If $s\leq
t\leq s+1$ only $II_1$ appears with $s+1$ replaced by $t$ and it is
estimated as follows:
\begin{align}
\|II_1\|_{L^2}&\leq\int_s^{s+1} \|\eht
P_c \|_{H^2\mapsto H^2}\|Dh|_{a(\tau)}\|_{\C\mapsto H^2}\left(\left| \<\po,g_u e^{-iH(\tau-s)}P_c v\>\right|+\left| \<\po,g_{\bar{u}} e^{iH(\tau-s)}P_c \bar{v}\>\right|\right) d\tau \nonumber \\
&\leq\int_s^{s+1} C\|Dh|_{a(\tau)}\|_{\C\mapsto H^2}\Big(\left| \<e^{-iH(\tau-s)}\po,e^{iH(\tau-s)}g_u e^{-iH(\tau-s)}P_c v\>\right|\nonumber \\
&\quad\quad\quad\quad\quad\quad\quad\quad\quad\quad+\left| \<e^{iH(\tau-s)}P_c\po,e^{-iH(\tau-s)}g_{\bar{u}} e^{iH(\tau-s)}P_c \bar{v}\>\right|\Big) d\tau \nonumber \\
&\leq\int_s^{s+1} C\|Dh|_{a(\tau)}\| \|\po\|_{L^{p}}(\|e^{iH(\tau-s)}g_u e^{-iH(\tau-s)}P_c v\|_{L^{p'}}+\|e^{-iH(\tau-s)}g_{\bar u} e^{iH(\tau-s)}P_c v\|_{L^{p'}}) d\tau \nonumber \\
&\leq\int_s^{s+1} C\|Dh|_{a(\tau)}\|_{\C\mapsto H^2}
\|\po\|_{L^{p}}(\|\widehat{g_u}\|_{L^1}+\|\widehat{g_{\bar
u}}\|_{L^1}) \|v\|_{L^{p'}} d\tau\leq C\|v\|_{L^{p'}} \nonumber
\end{align} For the
$II_2$ term we have again via Stricharz estimates for the admissible
pair $(\gamma,\rho),$ $\gamma\geq 2$ fixed but $\gamma\not= \infty:$
\begin{align}
\|II_2\|_{L^2}&\leq C\Big(\int_{s+1}^\infty \|Dh\<\po,g_u \ehts P_c v\>\|_{L^{\rho'}}^{\gamma'}d\tau\Big)^{\frac{1}{\gamma'}}+C\Big(\int_{s+1}^\infty \|Dh\<\po,g_{\bar{u}} e^{iH(\tau-s)} P_c \bar{v}\>\|_{L^{\rho'}}^{\gamma'}d\tau\Big)^{\frac{1}{\gamma'}} \nonumber \\
&\leq C\Big(\int_{s+1}^\infty\|Dh\|_{\C\mapsto L^{\rho'}}^{\gamma'}\|\po\|_{L^2}^{\gamma'}\| g_u\|_{L^{\beta}}^{\gamma'}\|\ehts P_c v\|_{L^p}^{\gamma'}d\tau\Big)^{\frac{1}{\gamma'}}\nonumber \\
&+C\Big(\int_{s+1}^\infty \|Dh\|_{\C\mapsto L^{\rho'}}^{\gamma'}\|\po\|_{L^2}^{\gamma'}\|g_{\bar{u}}\|_{L^{\beta}}^{\gamma'}\|e^{iH(\tau-s)} P_c \bar{v}\|_{L^p}^{\gamma'}d\tau\Big)^{\frac{1}{\gamma'}} \nonumber \\
&\leq C\Big(\int_{s+1}^\infty
\frac{d\tau}{|\tau-s|^{N(\frac{1}{2}-\frac{1}{p})\gamma'}}\Big)^{\frac{1}{\gamma'}}<\infty
\nonumber
\end{align} where $1/\beta+1/p=1/\rho'.$

Similarly we can estimate $L(s)W:$
\begin{align}
\|L(s)W(t)\|_{L^2}&\leq C\Big(\int_s^\infty \|\Dg[W(\tau)]\|_{L^{\rho'}}^{\gamma'}d\tau\Big)^{\frac{1}{\gamma'}}+C\Big(\int_s^\infty \|Dh\<\po,\Dg[W(\tau)]\>\|_{L^{\rho'}}^{\gamma'}d\tau\Big)^{\frac{1}{\gamma'}} \nonumber \\
&\leq C\Big(\int_s^\infty \left(\|\xs g_u\|_{L^{\frac{3\gamma}{2}}}+\|g_{\bar{u}}\|_{L^{\frac{3\gamma}{2}}}\right)^{\gamma'}\|W\|_{\Lsn}^{\gamma'}d\tau\Big)^{\frac{1}{\gamma'}}\nonumber\\
&\qquad+C\Big(\int_s^\infty \|Dh\|_{\C\mapsto L^{\rho'}}^{\gamma'}\|\po\|_{L^2}^{\gamma'} \left(\|\xs g_u\|_{L^\infty}+\|g_{\bar{u}}\|_{L^\infty}\right)^{\gamma'}\|W\|_{\Lsn}^{\gamma'}d\tau\Big)^{\frac{1}{\gamma'}}\nonumber \\
&\leq C\Big(\int_s^\infty
\frac{d\tau}{(1+|\tau-s|)^{3(\frac{1}{2}-\frac{1}{p})\gamma'}}\Big)^{\frac{1}{\gamma'}}<\infty
\nonumber \end{align} Hence $T(t,s)-\tilde{T}(t,s):L^{p'}\rightarrow
L^2$ is bounded for $p=\frac{2N}{N-2}$ and by part (i) for $p=2.$ By
Riesz-Thorin interpolation it is bounded for any $2\leq p\leq
\frac{2N}{N-2}.$ This finishes the proof of part $(iii)$ and the
theorem. $\Box$

\bibliographystyle{plain}
\def\cprime{$'$}

\end{document}